\documentclass[a4paper]{amsart}
\usepackage{amsmath,amssymb,amsfonts,graphicx, array} %

\numberwithin{equation}{section}
\theoremstyle{plain}
 \newtheorem{thm}{Theorem}[section]
 \newtheorem{cor}[thm]{Corollary}
 \newtheorem{lem}[thm]{Lemma}
 \newtheorem{prop}[thm]{Proposition}
\theoremstyle{definition}
 \newtheorem{defn}[thm]{Definition}
 \newtheorem{exmp}[thm]{Example}
\theoremstyle{remark}
 \newtheorem{rem}[thm]{Remark}
\DeclareMathOperator{\Ad}{Ad}
\DeclareMathOperator{\ad}{ad}
\DeclareMathOperator{\Ann}{Ann}
\DeclareMathOperator{\trace}{Trace}
\DeclareMathOperator{\Hom}{Hom}
\DeclareMathOperator{\End}{End}
\DeclareMathOperator{\Ker}{Ker}
\DeclareMathOperator{\sign}{sgn}

\def\ang#1#2{\langle #1,#2\rangle}
\def\ve{\varepsilon}

\allowdisplaybreaks

\begin{document}

\title[Annihilators of generalized Verma modules]
{Minimal polynomials and annihilators of generalized Verma modules 
of the scalar type}%
\author{Hiroshi Oda}
\address{Faculty of Engineering, Takushoku University,
815-1, Tatemachi, Hachioji-shi, Tokyo 193-0985, Japan}
\email{hoda@la.takushoku-u.ac.jp}
\author{Toshio Oshima}
\address{Graduate School of Mathematical Sciences,
University of Tokyo, 7-3-1, Komaba, Meguro-ku, Tokyo 153-8914, Japan}
\email{oshima@ms.u-tokyo.ac.jp}
\begin{abstract}
We construct a generator system of the annihilator of a generalized Verma
module of a reductive Lie algebra induced from a character of a 
parabolic subalgebra as an analogue of the minimal polynomial of
a matrix.
\end{abstract}
\maketitle

\section{Introduction}

In the representation theory of a real reductive Lie group $G$ the center 
$Z(\mathfrak g)$ of the universal enveloping algebra $U(\mathfrak g)$ of the 
complexification $\mathfrak g$ of the Lie algebra of $G$ plays an important role.
For example, any irreducible admissible representation $\tau$ of $G$ realized in 
a subspace $E$ of sections of a certain $G$-homogeneous vector bundle is a 
simultaneous eigenspace of $Z(\mathfrak g)$ parameterized by the infinitesimal 
character of $\tau$.
The differential equations induced from $Z(\mathfrak g)$ are often used to
characterize the subspace $E$.

If the representation $\tau$ is small, we expect more differential equations
corresponding to the primitive ideal $I_\tau$, that is, the annihilator of $\tau$ in
$U(\mathfrak g)$.
For the study of $I_\tau$ and these differential equations it is interesting and 
important to get a good generator system of $I_\tau$.

Let $\mathfrak p_\Theta$ be a parabolic subalgebra containing
a Borel subalgebra $\mathfrak b$ of $\mathfrak g$ and let $\lambda$ be a character
of $\mathfrak p_\Theta$.  Then the generalized Verma module of the scalar type is
by definition
\begin{equation}
 M_\Theta(\lambda)=U(\mathfrak g)/J_\Theta(\lambda)\quad\text{with }
J_\Theta(\lambda) = 
 \sum_{X\in\mathfrak p_\Theta}U(\mathfrak g)\bigl(X-\lambda(X)\bigr).
\end{equation}
In this paper we construct generator systems of the annihilator 
$\Ann\bigl(M_\Theta(\lambda)\bigr)$ of the generalized Verma module 
$M_\Theta(\lambda)$ in a unified way.
If $\tau$ can be realized in a space $E$ of sections of a line bundle over a 
generalized flag manifold, the annihilator of the corresponding generalized 
Verma module kills $E$.

When $\mathfrak g = \mathfrak{gl}_n$,  \cite{O-Cap} and \cite{O-Mat} construct 
such a generator system by {\it generalized Capelli operators\/} defined through 
{\it quantized elementary divisors.\/}  
This is a good generator system and in fact it is used there to characterize the image of 
the Poisson integrals on various boundaries of the symmetric space and also to 
define generalized hypergeometric functions.
A similar generator system is studied by \cite{Oda} for $\mathfrak g=\mathfrak o_n$
but it is difficult to construct the corresponding generator system in the case of 
other general reductive Lie groups.  
On the other hand, in \cite{O-Cl} we give other generator systems as a 
quantization of {\it minimal polynomials\/} when $\mathfrak g$ is classical.

Associated to a faithful finite dimensional representation $\pi$ of $\mathfrak g$ 
and a $\mathfrak g$-module $M$, \cite{O-Cl} defines a minimal polynomial 
$q_{\pi,M}(x)$ as is quoted in ~Definition~\ref{def:min}
and Definition~\ref{def:canomin}.
If $\mathfrak g=\mathfrak {gl}_n$ and $\pi$ is a natural representation of 
$\mathfrak g$, $q_{\pi,M}(x)$ is characterized by the condition 
$q_{\pi,M}(F_\pi)M=0$.
Here $F_\pi=\Bigl(E_{ij}\Bigr)_{\substack{1\le i\le n\\1\le j\le n}}$ is the matrix
whose $(i,j)$-component is the fundamental matrix unit $E_{ij}$ and then $F_\pi$ is 
identified with a square matrix with components in $\mathfrak g\subset U(\mathfrak g)$.
In this case $q_{\pi,M_\Theta(\lambda)}(x)$ is naturally regarded as a 
quantization of the minimal polynomial which corresponds to the conjugacy class of 
matrices given by a {\it classical limit\/} of $M_\Theta(\lambda)$.
For example, if $\mathfrak p_\Theta$ is a maximal parabolic subalgebra of 
$\mathfrak{gl}_n$, the minimal polynomial $q_{\pi,M_\Theta(\lambda)}(x)$ is a 
polynomial of degree 2.

For general $\pi$ and $\mathfrak g$, the matrix $F_\pi$ is the image 
$\Bigl(p(E_{ij})\Bigr)$ of $\bigl(E_{ij}\bigr)$ under the contragredient map $p$ of 
$\pi$ and then $F_\pi$ is a square matrix of the size $\dim\pi$ with components in 
$\mathfrak g$.
For example, if $\pi$ is the natural representation of $\mathfrak o_n$, then the
$(i,j)$-component of $F_\pi$ equals $\frac12(E_{ij} - E_{ji})$.

In \cite{O-Cl} we calculate the minimal polynomial $q_{\pi,M_\Theta(\lambda)}(x)$
for the natural representation $\pi$ of each type of classical Lie algebra $\mathfrak g$ and
by putting
\begin{equation}
 I_{\pi,\Theta}(\lambda)=\sum_{i,j}U(\mathfrak g)q_{\pi,M_\Theta(\lambda)}(F_\pi)_{ij}
 + \sum_{\Delta\in Z(\mathfrak g)\cap \Ann M_\Theta(\lambda)}
    \!\!\!\!\!\! U(\mathfrak g) \Delta,
\end{equation}
it is shown that
\begin{equation}\label{eq:0gap}
 J_{\Theta}(\lambda) = I_{\pi,\Theta}(\lambda) + J(\lambda_\Theta)
 \quad\text{with }J(\lambda_\Theta) = 
 \sum_{X\in\mathfrak b}U(\mathfrak g)\bigl(X-\lambda(X)\bigr)
\end{equation}
for a generic $\lambda$.
This equality is essential because it shows that $q_{\pi,M_\Theta(\lambda)}(F_\pi)_{ij}$
give elements killing $M_{\Theta}(\lambda)$ which cannot be described 
by  $Z(\mathfrak g)$ and define differential equations characterizing the local 
sections of the corresponding line bundle of a generalized flag manifold.
Moreover \eqref{eq:0gap} assures that $I_{\pi,\Theta}(\lambda)$ %
equals $\Ann\bigl(M_\Theta(\lambda)\bigr)$ 
for a generic $\lambda$ (Proposition~\ref{prop:HCzero}).

In this paper, $\pi$ may be any faithful irreducible finite dimensional 
representation of a reductive Lie algebra $\mathfrak g$.
In Theorem~\ref{thm:min} we calculate a polynomial 
$q_{\pi,\Theta}(x;\lambda)$ which is divisible by the 
the minimal polynomial $q_{\pi,M_\Theta(\lambda)}(x)$
and it is shown in Theorem~\ref{thm:minimality} 
that the former polynomial equals the latter for a generic $\lambda$.
If $\mathfrak p_\Theta = \mathfrak b$, this result gives the {\it characteristic
polynomial\/} associated to $\pi$ as is stated in Theorem~\ref{thm:char}, 
which is studied by \cite{Go2}.
We prove Theorem~\ref{thm:min} in a similar way as in \cite{O-Cl} but in a more
generalized way and the proof is used to get the condition for \eqref{eq:0gap}.
Another proof which is similar as is given in \cite{Go2} is also possible and 
it is based on the decomposition of the tensor product of some finite dimensional 
representations of $\mathfrak g$ given by Proposition~\ref{prop:tensor}.
The proof of Theorem~\ref{thm:minimality} uses
{\it infinitesimal Mackey's tensor product theorem\/} which is explained
in Appendix~\ref{app:Mackey}.

In \S3 we examine \eqref{eq:0gap} and obtain a sufficient condition for 
\eqref{eq:0gap} by Theorem~\ref{thm:gap}.
Proposition~\ref{prop:goodweight}
and Proposition~\ref{prop:goodrep} assure that
a generic $\lambda$ satisfies this condition if $\pi$ is one of many proper representations
including minuscule representations, adjoint representations,
representations of multiplicity free, and representations with regular highest weights.
In such cases
the sufficient condition is satisfied if $\lambda$ is not in the union of 
a certain finite number of complex hypersurfaces in the parameter space, which 
are defined by the difference of certain weights of the representation $\pi$.
On the other hand, in Appendix~\ref{app:badcase},
we give counter examples for which our sufficient condition is never satisfied by any $\lambda$.
In Proposition~\ref{prop:geninv} we also study the element of 
$Z(\mathfrak g)$ contained in $I_{\pi,\Theta}(\lambda)$.

A corresponding problem in the {\it classical limit\/} is to construct a generator
system of the defining ideal of the coadjoint orbit of $\mathfrak g$ and in fact 
Theorem~\ref{thm:orbit} is considered to be the classical limit of 
Corollary~\ref{cor:ann}.

If $\pi$ is smaller, the two-sided ideal $I_{\pi,\Theta}(\lambda)$ is better in 
general and therefore in \S4 we give examples of the characteristic polynomials 
of some small $\pi$ for every simple $\mathfrak g$ and describe some minimal 
polynomials, especially in each case where $\mathfrak p_\Theta$ is maximal.
Note that the minimal polynomial is a divisor of the characteristic 
polynomial evaluated at the infinitesimal character.
In Proposition~\ref{prop:every} we present a two-sided ideal of $U(\mathfrak g)$ 
for every $(\mathfrak g, \mathfrak p_\Theta)$ and examine the condition 
\eqref{eq:0gap} for this ideal by applying Theorem~\ref{thm:gap}.
In particular, the condition is satisfied if the infinitesimal character of 
$M_\Theta(\lambda)$ is regular in the case when $\mathfrak g=\mathfrak{gl}_n$, 
$\mathfrak o_{2n+1}$, $\mathfrak{sp}_n$ or $G_2$.  The condition is also satisfied
if the infinitesimal character is in the positive Weyl chamber containing the 
infinitesimal characters of the Verma modules which have finite dimensional 
irreducible quotients.

Some applications of our results in this paper to the integral geometry will be found 
in \cite[\S5]{O-Cl} and
\cite{OSh}.

\subsection*{Acknowledgments}
The referee gave the authors many helpful comments
to make the paper more accessible.
The authors greatly appreciate it.

\section{Minimal Polynomials and Characteristic Polynomials}\label{sec:min}

For an associative algebra $\mathfrak A$ and a positive integer $N$, 
we denote by $M(N,\mathfrak A)$ the associative algebra of square matrices of 
size $N$ with components in $\mathfrak A$.
We use the standard notation $\mathfrak{gl}_n$, $\mathfrak o_n$ and 
$\mathfrak{sp}_n$ for classical Lie algebras over $\mathbb C$.
The exceptional simple Lie algebra is denoted by its type $E_6$, $E_7$, 
$E_8$, $F_4$ or $G_2$.

The Lie algebra $\mathfrak{gl}_N$ is identified with 
$M(N,\mathbb C)\simeq\End(\mathbb C^N)$ with the bracket
$[X,Y] = XY - YX$.
In general, if we fix a base $\{v_1,\ldots,v_N\}$ of an $N$-dimensional 
vector space $V$ over $\mathbb C$, we naturally identify an element 
$X=(X_{ij})$ of $M(N,\mathbb C)$ with an element of 
$\End(V)$ by $Xv_j = \sum_{i=1}^NX_{ij}v_i$.
Let 
$E_{ij} = 
  \Bigl(\delta_{\mu i}\delta_{\nu j}\Bigr)_{
  \substack{1\le\mu\le N\\1\le\nu\le N}
 }
  \in M(N,\mathbb C)
$ be the standard matrix units and put $E_{ij}^*=E_{ji}$.
Note that the symmetric bilinear form
\begin{equation}\label{eq:bilinear}
 \ang{X}{Y} = \trace XY
 \quad\text{for}\quad X, Y\in\mathfrak{gl}_N
\end{equation}
on $\mathfrak{gl}_N$ is non-degenerate and satisfies
\begin{equation}
\begin{gathered}
  \ang{E_{ij}}{E_{\mu\nu}} = \ang{E_{ij}}{E_{\nu\mu}^*} = \delta_{i\nu}\delta_{j\mu},\\
  X = \sum_{i,j}\ang{X}{E_{ji}}E_{ij},\\
  \ang{\Ad(g)X}{\Ad(g)Y}=\ang XY\quad\text{for }
  X,\,Y\in\mathfrak{gl}_N\text{ and }g\in GL(N,\mathbb C).
\end{gathered}
\end{equation}
In general, for a Lie algebra $\mathfrak g$ over $\mathbb C$, we denote by
$U(\mathfrak g)$ and $Z(\mathfrak g)$ the {\it universal enveloping algebra\/}
of $\mathfrak g$ and the {\it center\/} of $U(\mathfrak g)$, respectively.
Then we have the following lemma.
\def\TMPCite{\cite[Lemma~2.1]{O-Cl}}
\begin{lem}[\TMPCite]\label{lem:1}
Let $\mathfrak g$ be a Lie algebra over $\mathbb C$ and let 
$(\pi,\mathbb C^N)$ be a representation of $\mathfrak g$.
Let $p$ be a linear map of 
$\mathfrak {gl}_N$ to $U(\pi(\mathfrak g))$ satisfying
\begin{equation}
  p([X, Y]) = [X, p(Y)]
\quad\text{for $X\in\pi(\mathfrak g)$ and $Y\in\mathfrak{gl}_N$},
\label{eq:lem1proj}
\end{equation}
that is,
$p\in\Hom_{\mathfrak \pi(\mathfrak g)}
\bigl(\mathfrak{gl}_N,U(\pi(\mathfrak g))\bigr)$.

Fix 
$q(x)\in\mathbb C[x]$
and put
\begin{equation}
\begin{cases}
 \ F = \Bigl(p(E_{ij})\Bigr)_{\substack{1\le i\le N\\1\le j\le N}}
 \in M\bigl(N,U(\pi(\mathfrak g))\bigr),\\
 \ \Bigl(Q_{ij}\Bigr)_{\substack{1\le i\le N\\1\le j\le N}}
 = q(F) \in M\bigl(N,U(\pi(\mathfrak g))\bigr).
\end{cases}
\end{equation}
Then
\begin{equation}\label{eq:adjE}
 \Bigl(p(\Ad(g)E_{ij})\Bigr)_{\substack{1\le i\le N\\1\le j\le N}}
 = {}^t\!g\,F\,{}^t\!g^{-1}\quad\text{for \ }g\in GL(n,\mathbb C)
\end{equation}
and
\begin{multline}
 [X, Q_{ij}] = \sum_{\mu=1}^N X_{\mu i}Q_{\mu j}
                  - \sum_{\nu=1}^N X_{j \nu}Q_{i\nu}\\
             = \sum_{\mu=1}^N\ang{X}{E_{i\mu}}Q_{\mu j}
               - \sum_{\nu=1}^N Q_{i\nu}\ang{X}{E_{\nu j}}
 \quad\text{for $X=\Bigl(X_{\mu\nu}\Bigr)_{\substack{1\le \mu\le N\\1\le \nu\le N}}
 \in\pi(\mathfrak g)$}. \label{eq:lem1}
\end{multline}
Hence the linear map $\mathfrak{gl}_N \rightarrow U(\pi(\mathfrak g))$
defined by $E_{ij} \mapsto Q_{ij}$ is an element of
$\Hom_{\mathfrak \pi(\mathfrak g)}
\bigl(\mathfrak{gl}_N,U(\pi(\mathfrak g))\bigr)$.
In particular, $\sum_{i=1}^NQ_{ii}\in Z(\pi(\mathfrak g))$.
\end{lem}
\begin{rem}
The referee suggested that we should give the reader
the following conceptual explanation of Lemma~\ref{lem:1}:
Since $\left(\mathfrak{gl}_N\right)^*\simeq M(N,\mathbb C)^*$ is naturally identified with
$M(N,\mathbb C)$ via \eqref{eq:bilinear},
the linear $\pi(\mathfrak g)$-homomorphism $p : \mathfrak{gl}_N\rightarrow U(\pi(\mathfrak g))$
is considered as an element of
$\left(\mathfrak{gl}_N\right)^*\otimes U(\pi(\mathfrak g))
 \simeq M(N,\mathbb C) \otimes U(\pi(\mathfrak g))
 \simeq M(N,U(\pi(\mathfrak g)))$.
By this identification,
the image of $p$ equals ${}^tF$ and hence \eqref{eq:adjE} holds almost immediately. 
Furthermore, \eqref{eq:lem1} is equivalent to the fact that
$\left(M(N,\mathbb C) \otimes U(\pi(\mathfrak g))\right)^{\pi(\mathfrak g)}$
is a subalgebra of $M(N,\mathbb C) \otimes U(\pi(\mathfrak g))
\simeq M(N,U(\pi(\mathfrak g)))$.
\end{rem}

Now we introduce the minimal polynomial defined by \cite{O-Cl}, which will be
studied in this section.

\begin{defn}[characteristic polynomials and minimal polynomials]\label{def:min}
Given a Lie algebra $\mathfrak g$, a faithful finite dimensional 
representation $(\pi, \mathbb C^N)$
and a $\mathfrak g$-homomorphism $p$ of 
$\End(\mathbb C^N)\simeq\mathfrak{gl}_N$ to $U(\mathfrak g)$.
Here we identify $\mathfrak g$ as a subalgebra of $\mathfrak{gl}_N$ through 
$\mathfrak \pi$.
Let $\Hat{Z}(\mathfrak g)$ denote the quotient field of $Z(\mathfrak g)$.
(Recall $Z(\mathfrak g)$ is an integral domain.)
Put $F=\Bigl(p(E_{ij})\Bigr)\in M(N,U(\mathfrak g))$.
We say $q_F(x)\in \Hat{Z}(\mathfrak g)[x]$
is the {\it characteristic polynomial\/} of $F$ if it is the monic polynomial
with the minimal degree which satisfies
\[
  q_F(F) = 0
\]
in $M\bigl(N,\hat Z(\mathfrak g)\otimes_{Z(\mathfrak g)} U(\mathfrak g)\bigr)$.
Suppose moreover a $\mathfrak g$-module $M$ is given.
Then we say $q_{F,M}(x)\in\mathbb C[x]$ is the
{\it minimal polynomial\/} of the pair $(F,M)$ if
it is the monic polynomial with the minimal degree which satisfies
\[
 q_{F,M}(F)M = 0.
\]
\end{defn}

\begin{rem}
The uniqueness of the characteristic (or minimal) polynomial is clear
if it exists.
Suppose $\mathfrak g$ is reductive.
Then the characteristic polynomial actually exists by \cite[Theorem~2.6]{O-Cl}.
The same theorem assures the existence of the minimal polynomial
if $M$ has a finite length or an infinitesimal character.
\end{rem}

\begin{defn}\label{def:canomin}
If the symmetric bilinear form \eqref{eq:bilinear} is non-degenerate
on $\pi(\mathfrak g)$, the orthogonal projection of 
$\mathfrak{gl}_N$ onto $\pi(\mathfrak g)$ satisfies the assumption
for $p$ in Lemma~\ref{lem:1}, which we call the 
{\it canonical projection\/} of $\mathfrak{gl}_N$ to $\pi(\mathfrak g)\simeq\mathfrak{g}$.
In this case we put $F_\pi= \Bigl(p(E_{ij})\Bigr)$.
Then we call $q_{F_\pi}(x)$ (resp.~$q_{F_\pi,M}(x)$)
in Definition~\ref{def:min} the characteristic polynomial of 
$\pi$ (resp.~the minimal polynomial of the pair $(\pi,M)$)
and denote it by $q_{\pi}(x)$ (resp.~$q_{\pi,M}(x)$).
\end{defn}

\begin{rem}
For a given involutive automorphism $\sigma$ of
$\mathfrak{gl}_N$,
put $$\mathfrak g= \{X\in\mathfrak{gl}_N;\ \sigma(X) = X\}$$
and let $\pi$ be the inclusion map of 
$\mathfrak g\subset\mathfrak{gl}_N$.
Then $p(X) = \frac{X + \sigma(X)}2$.
\end{rem}

Hereafter in the general theory of minimal polynomials
which we shall study,
we restrict our attention to a fixed finite dimensional representation 
$(\pi,V)$ of $\mathfrak g$ such that
\begin{equation}\label{eq:setting}
\begin{cases}
  \mathfrak g\text{ is a reductive Lie algebra over }\mathbb C,\\
  \pi\text{ is faithful and irreducible}.
\end{cases}
\end{equation}
Moreover we put $N=\dim V$ and 
identify $V$ with $\mathbb C^N$ through some basis of $V$. 
The assumption of Definition~\ref{def:canomin} is then satisfied.
\begin{rem}\label{rem:1} i)
The dimension of the center of 
$\mathfrak g$ is at most one.

\noindent
ii)
Fix $g\in GL(V)$.
If we replace $(\pi,V)$ by $(\pi^g,V)$ with 
$\pi^g(X)=\Ad(g)\pi(X)$ for $X\in\mathfrak g$ in Lemma~\ref{lem:1}, 
$F_\pi\in M(N,\mathfrak g)$ is naturally changed into 
${}^t\!g^{-1}F_\pi\,{}^t\!g$ 
under the fixed identification $V\simeq\mathbb C^N$.
This is clear from Lemma~\ref{lem:1} (cf.~\cite[Remark~2.7~ii)]{O-Cl}).
iii)
Exceptionally
the condition \eqref{eq:setting} will not be assumed
in Definition~\ref{defn:minuscule} and Proposition~\ref{prop:minuscule}.
\end{rem}

\begin{defn}[root system]\label{def:roots}
We fix a Cartan subalgebra $\mathfrak a$ of $\mathfrak g$ and
let $\Sigma(\mathfrak g)$ be a root system for the pair 
$(\mathfrak g,\mathfrak a)$.
We choose an order in $\Sigma(\mathfrak g)$ and denote by
$\Sigma(\mathfrak g)^+$ and $\Psi(\mathfrak g)$ the set of the positive roots
and the fundamental system, respectively.
For each root $\alpha\in\Sigma(\mathfrak g)$ we fix a root vector 
$X_\alpha\in\mathfrak g$.
Let $\mathfrak g=\mathfrak n\oplus\mathfrak a\oplus\bar{\mathfrak n}$
be the triangular decomposition of $\mathfrak g$ so that $\mathfrak n$ is 
spanned by $X_\alpha$ with $\alpha\in\Sigma(\mathfrak g)^+$.
We say $\mu\in\mathfrak a^*$ is {\it dominant\/} if and only if
\[
2\frac{\ang{\mu}{\alpha}}{\ang{\alpha}{\alpha}} \notin \{-1,-2,\ldots\}
\quad\text{for any }\alpha\in \Sigma(\mathfrak g)^+.
\] 
\end{defn}
Let us prepare some lemmas and definitions.
\begin{lem}\label{lem:proj}
Let $U$ be a $k$-dimensional subspace of $\mathfrak{gl}_N$
such that $\ang{\ }{\ }|_U$ is non-degenerate.
Let $p_U$ be the orthogonal projection of $\mathfrak{gl}_N$ to $U$ and
let $\{v_1,\ldots,v_k\}$ be a basis of $U$ with $\ang{v_2}{v_j}=0$ for 
$2\le j\le k$.
Suppose that $u\in\mathfrak{gl}_N$ satisfies $\ang{u}{v_j}=0$ for 
$2\le j\le k$. 
Then $p_U(u) = \frac{\ang{u}{v_1}}{\ang{v_1}{v_2}}v_2$.
\end{lem}
The proof of this lemma is easy and we omit it.
\begin{lem}\label{lem:2} 
Choose a base $\{v_i;\, i=1,\ldots,N\}$ of $V$
for the identification $V\simeq\mathbb C^N$
so that $v_i$ are weight vectors with 
weights $\varpi_i\in\mathfrak a^*$, respectively.
We identify $\mathfrak g$ with the subalgebra $\pi(\mathfrak g)$ of
$\mathfrak{gl}_N\simeq M(N,\mathbb C)$ and put 
$\mathfrak a_N=\sum_{i=1}^N\mathbb CE_{ii}$.
For $F_\pi = \Bigl(F_{ij}\Bigr)_{
  \substack{1\le i \le N\\1\le j \le N}
 }$
we have
\begin{equation}\label{eq:F1}
\begin{aligned}
{}  &F_{ii} = \varpi_i = \sum_{j=1}^N\varpi_i(E_{jj})E_{jj},\\
  &\ad(H)(F_{ij})=(\varpi_i-\varpi_j)(H)F_{ij}\ (\forall H\in\mathfrak a),\\
  &\ang{F_{ij}}{E_{\mu\nu}}\ne 0\text{\ with $i\ne j$ implies }
  \varpi_i - \varpi_j = \varpi_\nu - \varpi_\mu\in\Sigma(\mathfrak g),\\
 &\mathfrak a =\sum_{i=1}^N\mathbb CF_{ii}\subset\mathfrak a_N,\ 
 \mathfrak n =\sum_{\varpi_i-\varpi_j\in\Sigma(\mathfrak g)^+}
    \mathbb CF_{ij},\ 
  \bar{\mathfrak n} =\sum_{\varpi_j-\varpi_i\in\Sigma(\mathfrak g)^+}
    \mathbb CF_{ij}\ 
\end{aligned}
\end{equation}
under the identification 
$\mathfrak a^*\simeq\mathfrak a\subset\mathfrak a_N\simeq 
\mathfrak a_N^*$ by the bilinear form \eqref{eq:bilinear}.
\end{lem}
\begin{proof}
Note that
$H\in\mathfrak a$ is identified with 
$\sum_{j=1}^N\varpi_j(H)E_{jj}\in\mathfrak a_N\subset\mathfrak{gl}_N$.
Hence $\ad(H)(E_{ij})=(\varpi_i-\varpi_j)(H)E_{ij}$ and therefore
$\ad(H)(F_{ij})=(\varpi_i-\varpi_j)(H)F_{ij}.$
In particular we have $F_{ii}\in\mathfrak{a}$.
Since
\[
\ang H{F_{ii}} =
\ang H{E_{ii}} = 
\ang{\sum_{j=1}^N\varpi_j(H)E_{jj}}{E_{ii}} = \varpi_i(H)
\quad(\forall H\in\mathfrak a),
\]
we get $F_{ii} = \varpi_i$.

For each root $\alpha$,
the condition $(X_\alpha)_{ij}=\ang{X_\alpha}{E_{ji}}\ne 0$ 
means $\varpi_i - \varpi_j=\alpha$.
Hence if $i\ne j$ and 
$X\in\mathfrak a + \sum_{\alpha\in\Sigma(\mathfrak g),\,\alpha\ne\varpi_j-\varpi_i}
\mathbb C X_\alpha$, then $\ang{E_{ij}}{X}=0$ and therefore $\ang{F_{ij}}{X}=0$.
Hence $F_{ij}=0$ if $i\ne j$ and $\varpi_j-\varpi_i\not\in\Sigma(\mathfrak g)$.
On the other hand, if $\varpi_j-\varpi_i\in\Sigma(\mathfrak g)$, 
we can easily get
$F_{ij}=CX_{\varpi_i-\varpi_j}$ for some $C\in\mathbb C$.
Hence $\ang{F_{ij}}{E_{\mu\nu}}=0$ if $\varpi_i-\varpi_j\ne\varpi_\nu-\varpi_\mu$.
\end{proof}

Through the identification of $\mathfrak{a}^* \simeq \mathfrak{a}\subset\mathfrak{a}_N$
in the lemma,
we introduce the symmetric bilinear form $\ang{\ }{\ }$
on $\mathfrak{a}^*$.
We note this bilinear form is real-valued and positive definite on 
$\sum_{\alpha\in\Psi(\mathfrak g)}\mathbb{R}\alpha$.

Now we take a subset $\Theta\subset \Psi(\mathfrak g)$ with $\Theta\ne \Psi(\mathfrak g)$
and fix it.

\begin{defn}[generalized Verma module]\label{def:gV}
Put
\begin{align*}
 \mathfrak a_\Theta &= \{H\in\mathfrak a;\, \alpha(H)=0,
  \quad\forall\alpha\in\Theta\},\\ 
 \mathfrak g_\Theta &= \{X\in\mathfrak g;\, [X,H] = 0,
  \quad\forall H\in\mathfrak a_\Theta\},\\
 \mathfrak m_\Theta &= \{X\in\mathfrak g_\Theta;\, \ang{X}{H}=0,\quad
  \forall H\in\mathfrak a_\Theta\},\\
 \Sigma(\mathfrak g)^- &= \{\alpha;\, -\alpha\in\Sigma(\mathfrak g)^+\},\\
 \Sigma(\mathfrak g_\Theta) &= \{\alpha\in\Sigma(\mathfrak g);\, \alpha(H)=0,
  \quad\forall H\in\mathfrak a_\Theta\},\\
 \Sigma(\mathfrak g_\Theta)^+ &= \Sigma(\mathfrak g_\Theta)\cap \Sigma(\mathfrak g)^+,\quad
 \Sigma(\mathfrak g_\Theta)^- = \{-\alpha;\,\alpha\in\Sigma(\mathfrak g_\Theta)^+\},\\
 \mathfrak n_\Theta &= \sum_{\alpha\in \Sigma(\mathfrak g)^+
    \setminus\Sigma(\mathfrak g_\Theta)}\mathbb C X_\alpha,
    \quad \bar{\mathfrak n}_\Theta = 
  \sum_{\alpha\in \Sigma(\mathfrak g)^-\setminus\Sigma(\mathfrak g_\Theta)}
    \mathbb C X_\alpha,\\
  \mathfrak b &= \mathfrak a + \mathfrak n,\quad
  \mathfrak p_\Theta = \mathfrak g_\Theta + \mathfrak n_\Theta,\\
  \rho &= \frac12\sum_{\alpha\in\Sigma(\mathfrak g)^+}\alpha,\quad
  \rho(\Theta) = \frac12\sum_{\alpha\in\Sigma(\mathfrak g_\Theta)^+}\alpha,\quad
  \rho_\Theta = \rho - \rho(\Theta).
\end{align*}
For $\Lambda\in\mathfrak a^*$ which satisfies 
$2\frac{\ang{\Lambda}{\alpha}}{\ang{\alpha}{\alpha}} \in \{0, 1, 2, \ldots \}$
for $\alpha\in\Theta$,
let $U_{(\Theta, \Lambda)}$ denote
the finite dimensional irreducible $\mathfrak g_\Theta$-module
with highest weight $\Lambda$.
By the trivial action of $\mathfrak n_\Theta$,
we consider $U_{(\Theta, \Lambda)}$ to be a $\mathfrak p_\Theta$-module.
Put 
\begin{equation}
M_{(\Theta, \Lambda)} = U(\mathfrak g) \otimes_{U(\mathfrak p_\Theta)} U_{(\Theta, \Lambda)}.
\end{equation}
Then $M_{(\Theta, \Lambda)}$ is called a 
{\it generalized Verma module of the finite type\/}.
\end{defn}

\begin{rem}\label{rem:gV}
i)
$\mathfrak p_\Theta$ is a parabolic subalgebra containing the Borel subalgebra
$\mathfrak b$.
$\mathfrak p_\Theta=\mathfrak m_\Theta + \mathfrak a_\Theta + \mathfrak n_\Theta$
gives its direct sum decomposition.

\noindent
ii)
Every finite dimensional irreducible $\mathfrak p_\Theta$-module
is isomorphic to $U_{(\Theta, \Lambda)}$
with a suitable choice of $\Lambda$.
 
\noindent
iii)
$M_{(\emptyset, \Lambda)}$ is nothing but the Verma module
for the highest weight $\Lambda\in\mathfrak a^*$.

\noindent
iv)
Let $u_{\Lambda}$ be a highest weight vector of $U_{(\Theta,\Lambda)}$. 
Then $1\otimes u_{\Lambda}$ is a highest weight vector of $M_{(\Theta,\Lambda)}$.
Moreover $1\otimes u_{\Lambda}$ generates $M_{(\Theta,\Lambda)}$ because
\[
\begin{aligned}
M_{(\Theta,\Lambda)}&=U(\mathfrak g)\otimes_{U(\mathfrak p_\Theta)}U_{(\Theta,\Lambda)}
=U(\bar{\mathfrak n}_\Theta)\otimes_{\mathbb C}U(\mathfrak p_\Theta)\otimes_{U(\mathfrak p_\Theta)}U_{(\Theta,\Lambda)} \\
&=U(\bar{\mathfrak n}_\Theta)\otimes_{\mathbb C}U_{(\Theta,\Lambda)}
=U(\bar{\mathfrak n}_\Theta)\otimes_{\mathbb C}U(\bar{\mathfrak n}\cap\mathfrak g_\Theta) u_{\Lambda}
=U(\bar{\mathfrak n}) (1\otimes u_{\Lambda}).
\end{aligned}
\]
Hence $M_{(\Theta,\Lambda)}$ is a highest weight module
and is therefore a quotient of the Verma module $M_{(\emptyset, \Lambda)}$.

\noindent
v)
If $\ang{\Lambda}{\alpha}=0$ for each $\alpha\in\Theta$,
then $\dim U_{(\Theta, \Lambda)}=1$
and we have the character $\lambda_\Theta$ of $\mathfrak p_\Theta$
such that $X u_\Lambda = \lambda_\Theta(X) u_\Lambda$ for $X\in \mathfrak p_\Theta$.
Since $$U(\mathfrak g)=U(\bar{\mathfrak n}_\Theta)\oplus
\sum_{X\in\mathfrak p_\Theta}U(\mathfrak g)\bigl(X-\lambda_\Theta(X)\bigr)$$
is a direct sum and
$M_{(\Theta, \Lambda)} = U(\bar{\mathfrak n}_\Theta)\otimes_{\mathbb C}\mathbb C u_\Lambda$,
we have the kernel of the surjective $U(\mathfrak g)$-homomorphism
$U(\mathfrak g) \rightarrow M_{(\Theta, \Lambda)}$
defined by $D\mapsto D (1\otimes u_\Lambda)$
equals $\sum_{X\in\mathfrak p_\Theta}U(\mathfrak g)\bigl(X-\lambda_\Theta(X)\bigr)$.
\end{rem}

\begin{defn}[generalized Verma module of the scalar type]\label{def:gVs}
For $\lambda\in\mathfrak a_\Theta^*$ 
define a character $\lambda_\Theta$ of $\mathfrak p_\Theta$
by $\lambda_\Theta(X+H) = \lambda(H)$ for $X\in\mathfrak m_\Theta + \mathfrak n_\Theta$ and
$H\in\mathfrak a_\Theta$.
Put
\begin{equation}
\begin{aligned}
   J_\Theta(\lambda) &= \sum_{X\in\mathfrak p_\Theta}U(\mathfrak g)
     \bigl(X-\lambda_\Theta(X)\bigr),\\
   J(\lambda_\Theta) &= \sum_{X\in\mathfrak b}U(\mathfrak g)
     \bigl(X-\lambda_\Theta(X)\bigr),\\
   M_\Theta(\lambda) &= U(\mathfrak g)/J_\Theta(\lambda),\quad
   M(\lambda_\Theta) = U(\mathfrak g)/J(\lambda_\Theta).
\end{aligned}
\end{equation}
Then $M_\Theta(\lambda)$ is 
isomorphic to $M_{(\Theta, \lambda_\Theta)}$,
which is called a {\it generalized Verma module of the scalar type\/}.
If $\Theta=\emptyset$, we denote $J_\emptyset(\lambda)$ and
$M_\emptyset(\lambda)$ by $J(\lambda)$ and $M(\lambda)$, respectively.
\end{defn}

\begin{defn}[Weyl group]
Let $W$ denote the Weyl group of $\Sigma(\mathfrak g)$, which
is generated by the reflections 
$w_\alpha:\mathfrak a^*\ni\mu\mapsto
\mu - 2\frac{\ang{\mu}{\alpha}}{\ang{\alpha}{\alpha}}\alpha\in\mathfrak a^*$
with respect to $\alpha\in\Psi(\mathfrak g)$.
Put
\begin{equation}
\begin{aligned}
W_\Theta &= \{w\in W;\,w\bigl(\Sigma(\mathfrak g)^+\setminus
  \Sigma(\mathfrak g_\Theta)\bigr) =
  \Sigma(\mathfrak g)^+\setminus\Sigma(\mathfrak g_\Theta)\},\\
W(\Theta) &= \{w\in W;\,w\bigl(\Sigma(\mathfrak g_\Theta)^+\bigr)
  \subset\Sigma(\mathfrak g)^+\}.
\end{aligned}
\end{equation}
Then 
each element $w\in W(\Theta)$ is a unique element with the smallest length in
the right coset $w W_\Theta$ and
the map $W(\Theta)\times W_\Theta\ni(w_1,w_2)\mapsto w_1w_2\in W$ is a bijection.

For $w\in W$ and $\mu\in\mathfrak a^*$, define
\begin{equation}
 w.\mu = w(\mu+\rho)-\rho.
\end{equation}
\end{defn}

Here we note that $W_\Theta$ is generated by the reflections $w_\alpha$
with $\alpha\in\Theta$ and 
\begin{equation}
  \ang{\rho_\Theta}{\alpha} = 0
   \quad\text{for \ }\alpha\in\Sigma(\mathfrak g_\Theta).
\end{equation}

\begin{defn}[infinitesimal character]
Let $D\in U(\mathfrak g)$.
We denote by $D_\mathfrak a$ the element of $U(\mathfrak a)$
which satisfies $D-D_\mathfrak a\in\bar{\mathfrak n}U(\mathfrak g)+U(\mathfrak g)\mathfrak n$
and identify $D_\mathfrak a \in U(\mathfrak a) \simeq S(\mathfrak a)$ with a polynomial function on $\mathfrak a^*$.
Then
$\Delta_\mathfrak a(\mu)=\Delta_\mathfrak a(w.\mu)$
for $\Delta\in Z(\mathfrak g)$,
$\mu\in\mathfrak a^*$, and $w\in W$.

Let $\mu\in\mathfrak a^*$.
We say
a $\mathfrak g$-module $M$ has {\it infinitesimal character\/} $\mu$
if each $\Delta\in Z(\mathfrak g)$ operates by the scalar $\Delta_\mathfrak a(\mu)$ in $M$.
We say an infinitesimal character $\mu$ is regular
if $\ang{\mu+\rho}{\alpha}\ne0$ for any $\alpha\in\Sigma(\mathfrak g)$.
\end{defn}

\begin{rem}\label{rem:inf_char}
The generalized Verma module
$M_{(\Theta,\Lambda)}$ in Definition~\ref{def:gV}
has infinitesimal character $\Lambda$.
It is clear by Remark~\ref{rem:gV} iv).
\end{rem}

\begin{defn}[Casimir operator]\label{def:Casimir}
Let $\{X_i;\,i=1,\ldots,\omega\}$ be a basis of $\mathfrak g$.
Then put
\[
 \Delta_\pi = \sum_{i=1}^\omega X_iX_i^*
\]
with the dual basis $\{X_i^*\}$ of $\{X_i\}$ with respect to the 
symmetric bilinear form \eqref{eq:bilinear} under the identification
$\mathfrak g\subset\mathfrak{gl}_N$ through $\pi$ and call $\Delta_\pi$
the {\it Casimir operator\/} of $\mathfrak g$ for $\pi$.
\end{defn}

\begin{rem}
As is well-known, $\Delta_\pi\in Z(\mathfrak g)$ and
$\Delta_\pi$ does not depend on the choice of $\{X_i\}$.
\end{rem}

We may assume 
in Definition~\ref{def:Casimir}
that $\{X_1,\ldots,X_{\omega'}\}$ and $\{X_{\omega'+1},\ldots,X_\omega\}$
be bases of $\mathfrak g_\Theta$ and 
$\bar{\mathfrak n}_\Theta+\mathfrak n_\Theta$,
respectively.
Then
$X_i^*\in\mathfrak g_\Theta$ for $i=1,\ldots,\omega'$ and
\begin{equation}
   \Delta_\pi^\Theta = \sum_{i=1}^{\omega'}X_iX_i^*
\end{equation}
is the Casimir operator of $\mathfrak g_\Theta$ for $\pi$.
\begin{lem}\label{lem:casimir}
Fix a basis $\{H_1,\ldots,H_r\}$ of the Cartan 
subalgebra $\mathfrak a$ of $\mathfrak g$. 

\noindent
{\rm i)}
Let $\{H_1^*,\ldots,H_r^*\}$ be the dual basis of $\{H_1,\ldots,H_r\}$.
Put $H_\alpha=[X_\alpha,X_{-\alpha}]$.
Then
\begin{align*}
\Delta_\pi
   &= 
   \sum_{\alpha\in\Sigma(\mathfrak g)}
     \frac{X_\alpha X_{-\alpha}}{\ang{X_\alpha}{X_{-\alpha}}}
     + \sum_{i=1}^r H_iH_i^*\\
   &= \sum_{i=1}^r H_iH_i^*
      + \sum_{\alpha\in\Sigma(\mathfrak g)^+}
      \left(
      \frac{2X_{-\alpha} X_{\alpha}}{\ang{X_\alpha}{X_{-\alpha}}}
     +  \frac{\alpha(H_\alpha)H_\alpha}{\ang{H_\alpha}{H_\alpha}}\right) \\
   &= \Delta_\pi^\Theta
     + \sum_{\alpha\in\Sigma(\mathfrak g)^+\setminus\Sigma(\mathfrak g_\Theta)}
      \left(
      \frac{2X_{\alpha} X_{-\alpha}}{\ang{X_\alpha}{X_{-\alpha}}}
     -  \frac{\alpha(H_\alpha)H_\alpha}{\ang{H_\alpha}{H_\alpha}}\right).
\end{align*}

\noindent
{\rm ii)}
Let $M$ be a highest weight module of $\mathfrak g$
with highest weight $\mu \in \mathfrak a^*$.
Then $\Delta_\pi v = \ang{\mu}{\mu+2\rho} v$
for any $v\in M$.

\noindent
{\rm iii)}
Let $v$ be a weight vector of $\pi$ belonging to an irreducible
representation of $\mathfrak g_\Theta$ realized as a subrepresentation 
of $\pi|_{\mathfrak g_\Theta}$ and let $\varpi$ denote the lowest weight
of the irreducible subrepresentation.
Then
\begin{align*}
  \Delta_\pi v &= 
    \ang{\bar\pi}{\bar\pi - 2\rho} v,\\
  \Delta_\pi^\Theta v &= 
    \ang{\varpi}{\varpi - 2\rho(\Theta)} v,\\
  \sum_{\alpha\in\Sigma(\mathfrak g)^+\setminus\Sigma(\mathfrak g_\Theta)}
      \frac{X_{\alpha} X_{-\alpha}}{\ang{X_\alpha}{X_{-\alpha}}}v
   &= \frac12\ang{\bar\pi-\varpi}{\bar\pi+\varpi - 2\rho}v.
\end{align*}
Here $\bar\pi$ denotes the lowest weight of $\pi$.

\noindent
{\rm iv)}
Fix $\beta\in\Sigma(\mathfrak g)^+$ and put 
$\mathfrak g(\beta) = 
\mathbb C X_\beta + \mathbb C X_{-\beta}
 + \sum_{i=1}^r\mathbb CH_i$.
Let $v$ be a weight vector of $\pi$ belonging to an irreducible representation
of $\mathfrak g(\beta)$ realized as a subrepresentation of 
$\pi|_{\mathfrak g(\beta)}$ and let $\varpi$ denote the lowest weight
of the irreducible subrepresentation.
Let $\varpi+\ell\beta$ be the weight of $v$.
Then
\begin{equation}\label{eq:beta}
 \frac{X_{\beta}X_{-\beta}}{\ang{X_\beta}{X_{-\beta}}}v = 
  -\Bigl(\ell\ang{\varpi+\frac{\ell-1}2\beta}\beta\Bigr)v.
\end{equation}

\noindent
{\rm v)} Suppose $\mathfrak g$ is simple.
Let $\alpha_{\max}$ is the maximal root of $\Sigma(\mathfrak g)^+$ and
let $B(\ ,\ )$ be the Killing form of $\mathfrak g$. Then
\[
  B(\alpha_{\max},\alpha_{\max} + 2\rho) = 1.
\]
\end{lem}
\begin{proof} i)
Note that
\begin{equation}
\ang{H_\alpha}{H_\alpha}
 = \ang{H_\alpha}{[X_\alpha,X_{-\alpha}]}
 =\ang{[H_\alpha,X_\alpha]}{X_{-\alpha}}
 = \alpha(H_\alpha)\ang{X_\alpha}{X_{-\alpha}}
\end{equation}
Since the dual base of 
$\{X_\alpha,\ H_i;\,\alpha\in\Sigma(\mathfrak g),\ i=1,\ldots,r\}$
equals
$\{\frac{X_{-\alpha}}{{\ang{X_\alpha}{X_{-\alpha}}}},\ H_i^*;\,
\alpha\in\Sigma(\mathfrak g),\ i=1,\ldots,r\}$, the claim is clear.

ii)
Let $v_\mu$ be a highest weight vector of $M$.
Then
\begin{align*}
  \Delta_\pi v_\mu &=
  \sum_{i=1}^rH_iH_i^* v_\mu
  +  \sum_{\alpha\in\Sigma(\mathfrak g)^+}
     \frac{\alpha(H_\alpha)H_\alpha}
     {\ang{H_\alpha}{H_\alpha}} v_\mu\\
  &= \sum_{i=1}^r\mu(H_i)\mu(H_i^*) v_\mu
   + \sum_{\alpha\in\Sigma(\mathfrak g)^+}
  \frac{\alpha(H_\alpha)\mu(H_\alpha)}{\ang{H_\alpha}{H_\alpha}}v_\mu.
\end{align*}
Hence $\Delta_\pi v_\mu = \ang{\mu}{\mu + 2\rho} v_\mu$
because $H_\alpha$ is
a non-zero constant multiple of $\alpha$ with the identification 
$\mathfrak a^* \simeq \mathfrak a$ by $\ang{\ }{\ }$ and therefore
$\Delta_\pi v = \ang{\mu}{\mu + 2\rho}v$ because
$M$ is generated by $v_\mu$.  

iii) 
Let $v_{\bar\pi}$ be a lowest weight vector of $\pi$.
Then we have $\Delta_\pi v_{\bar\pi} = \ang{\bar\pi}{\bar\pi - 2\rho} v_{\bar\pi}$
and therefore
$\Delta_\pi v = \ang{\bar\pi}{\bar\pi - 2\rho}v$.
Similarly we have $\Delta_\pi^\Theta v = \ang{\varpi}{\varpi-2\rho(\Theta)}v$.

Let $\varpi'$ be the weight of $v$.  Then we have
\begin{align*}
  \sum_{\alpha\in\Sigma(\mathfrak g)^+
  \setminus\Sigma(\mathfrak g_\Theta)}
      \frac{X_{\alpha} X_{-\alpha}}{\ang{X_\alpha}{X_{-\alpha}}}v
 &= \frac12 \Delta_\pi v - \frac12 \Delta_\pi^\Theta v
   + \ang{\varpi'}{\rho_\Theta}v\\
 &= \frac12\ang{\bar\pi-\varpi}{\bar\pi+\varpi - 2\rho}v.
\end{align*}
Here we note that 
$\ang{\varpi'}{\rho_\Theta}=\ang{\varpi}{\rho_\Theta}$.

iv) By the same argument as above we have
\[
 \frac{2X_{\beta}X_{-\beta}}{\ang{X_\beta}{X_{-\beta}}}v 
 + \sum_{i=1}^rH_iH_i^*v 
 - \frac{\beta(H_\beta)H_\beta}{\ang{H_\beta}{H_\beta}}v
 = \ang{\varpi}{\varpi - \beta}v.
\]
Hence
\begin{align*}
 \frac{2X_{\beta}X_{-\beta}}{\ang{X_\beta}{X_{-\beta}}}v &=
 \ang{\varpi}{\varpi - \beta}v - \ang{\varpi + \ell\beta}{\varpi + \ell\beta}v
 + \ang{\beta}{\varpi + \ell\beta}v\\
 &= -\bigl(2\ell\ang{\varpi}{\beta}+\ell(\ell-1)\ang{\beta}{\beta}\bigr)v
\end{align*}

v)
Suppose $\pi$ is the adjoint representation of the simple Lie algebra
$\mathfrak g$.
Then for $H\in\mathfrak a$ we have
\begin{align*}
\ang{\pi(\Delta_\pi)(H)}{H} &= \sum_{\alpha\in\Sigma(\mathfrak g)}
      \frac{\ang{[X_\alpha,[X_{-\alpha},H]]}H}{\ang{X_\alpha}{X_{-\alpha}}}\\
      &= \sum_{\alpha\in\Sigma(\mathfrak g)}
      \frac{-\ang{[X_{-\alpha},H]}{[X_\alpha,H]}}
           {\ang{X_\alpha}{X_{-\alpha}}}\\
      &= \sum_{\alpha\in\Sigma(\mathfrak g)}\alpha(H)^2\\
      &= \ang H H.
\end{align*}
Hence $\pi(\Delta_\pi)(H) = H$ and $B(\alpha_{\max},\alpha_{\max}+2\rho) = 
B(-\alpha_{\max},-\alpha_{\max}-2\rho)=1$.
\end{proof}
\begin{defn}[weights]\label{def:shift}
Let $\mathcal W(\pi)$ denote the set of the weights of the finite dimensional 
irreducible representation $\pi$ of $\mathfrak g$.
For $\varpi\in\mathcal W(\pi)$ define a real constant
\begin{equation}
D_\pi(\varpi) = \frac12\ang{\bar\pi-\varpi}{\bar\pi+\varpi-2\rho}.
\end{equation}
Here $\bar\pi$ is the lowest weight of $\pi$.
Put $R_+ = \{\sum_{\alpha\in\Psi(\mathfrak g)}m_\alpha\alpha;\,
m_\alpha \in\{0,1,2,\ldots\}\}$.
We define a partial order among the elements of $\mathcal W(\pi)$
so that $\varpi\le\varpi'$ if and only if $\varpi'-\varpi\in R_+$.

Moreover we put
\begin{equation}\begin{aligned} 
 \mathcal W_\Theta(\pi) &= \{
  \varpi\text{ are the highest weights of the irreducible components of }
  \pi|_{\mathfrak g_\Theta}\},\\
 \overline{\mathcal W}_\Theta(\pi) &= \{
  \varpi\text{ are the lowest weights of the irreducible components of }
  \pi|_{\mathfrak g_\Theta}\},\\
 \mathcal W(\pi)|_{\mathfrak a_\Theta} &= \{ 
  \varpi|_{\mathfrak a_\Theta};\,
  \varpi \in \mathcal W(\pi) \}. 
\end{aligned}\end{equation}
Let $\mu$ and $\mu' \in \mathcal W(\pi)|_{\mathfrak a_\Theta}$.
Then we define $\mu\le_\Theta\mu'$ if and only if
$\mu'-\mu\in 
\{\sum_{\alpha\in\Psi(\mathfrak g) \setminus\Theta}m_\alpha
\alpha|_{\mathfrak a_\Theta};\,
m_\alpha \in\{0,1,2,\ldots\}\}$.
\end{defn}

\begin{rem} i)
$\mathcal W_\emptyset(\pi)=\overline{\mathcal W}_\emptyset(\pi) 
= \mathcal W(\pi)$ and 
$\overline{\mathcal W}_\Theta(\pi) = -\mathcal W_\Theta(\pi^*)$.
Here $(\pi^*,V^*)$ denotes the contragredient representation of $(\pi,V)$
defined by 
\begin{equation}
(\pi^*(X)v^*)(v)=-v^*(\pi(X)v) \text{ for } X\in\mathfrak g,\ 
v^*\in V^* \text{ and } v\in V.
\end{equation}

\noindent
ii)
$\mathcal W(\pi)|_{\mathfrak a_\Theta}
=\{ \varpi|_{\mathfrak a_\Theta};\,
\varpi \in \mathcal{W}_\Theta(\pi) \}
=\{ \varpi|_{\mathfrak a_\Theta};\,
\varpi \in \overline{\mathcal W}_\Theta(\pi) \}$. 

\noindent
iii)
Suppose $\varpi$ and $\varpi'\in\mathcal W(\pi)$ 
and put $\varpi'-\varpi=\sum_{\alpha\in\Psi(\mathfrak g)}m_\alpha\alpha$.
Then $\varpi|_{\mathfrak a_\Theta}
\le_\Theta
\varpi'|_{\mathfrak a_\Theta}$
if and only if $m_\alpha\ge 0$ for
any $\alpha\in\Psi(\mathfrak g) \setminus \Theta$.
Hence $\bar\pi|_{\mathfrak a_\Theta}$ is the smallest element
of $\mathcal W(\pi)|_{\mathfrak a_\Theta}$.
Note that $\varpi\le\varpi'$ if and only if $\varpi\le_\emptyset\varpi'$.
\end{rem}

\begin{lem}\label{lem:order}
Let $\varpi$ and $\varpi'\in\mathcal W(\pi)$.

\noindent
{\rm i)} If $\alpha=\varpi'-\varpi\in\Psi(\mathfrak g)$, then
$D_\pi(\varpi)-D_\pi(\varpi')=\ang{\varpi}{\varpi'-\varpi}$.

\noindent
{\rm ii)} Suppose $\varpi'\in\overline{\mathcal W}_\Theta(\pi)$, 
$\varpi<\varpi'$ and 
$\varpi|_{\mathfrak a_\Theta}=\varpi'|_{\mathfrak a_\Theta}$.
Then $D_\pi(\varpi)<D_\pi(\varpi')$.
\end{lem}
\begin{proof}
ii) Note that
\[
 D_\pi(\varpi)-D_\pi(\varpi') = 
  \frac12\ang{\varpi'-\varpi}{\varpi+\varpi'-2\rho}.
\]
The assumption in ii) implies 
$\varpi'-\varpi = \sum_{\alpha\in\Theta}m_\alpha\alpha$ with
$m_\alpha\ge 0$. Here at least one of $m_\alpha$ is positive.
Hence $\ang{\varpi'-\varpi}{\rho}>0$.
Since $\varpi'$ are the lowest weights of irreducible
representations of $\mathfrak g_\Theta$,
$\ang{\alpha}{\varpi'}\le0$ for $\alpha\in\Theta$.
Thus we have 
$\ang{\sum_{\alpha\in\Theta}
m_\alpha\alpha}{2\varpi'-\sum_{\alpha\in\Theta}m_\alpha\alpha-2\rho}<0$.

i)
Put $\alpha=\varpi'-\varpi$.
Then
\[D_\pi(\varpi)-D_\pi(\varpi') - \ang{\varpi}{\varpi'-\varpi}
 = -\frac12\ang{\alpha}{2\rho - \alpha},
\]
which equals 0 if $\alpha\in\Psi(\mathfrak g)$
because 
$w_\alpha(\Sigma(\mathfrak g)^+\setminus\{\alpha\})=
\Sigma(\mathfrak g)^+\setminus\{\alpha\}$.
\end{proof}

Now we give a key lemma which is used to calculate our minimal polynomial.
\begin{lem}\label{lem:zen}
Fix an irreducible decomposition $\bigoplus_{i=1}^\kappa(\pi_i,V_i)$ of 
$(\pi|_{\mathfrak g_\Theta}, V)$ and a basis $\{v_{i,1},\ldots,v_{i,m_i}\}$
of $V_i$ so that $v_{i,j}$ are weight vectors for $\mathfrak a$.
Let $\varpi_{i,j}$ and $\varpi_i$ be the weight of $v_{i,j}$ and the lowest weight 
of the representation $\pi_i$, respectively.

Suppose $\varpi_{i,j} = \varpi_{i',j'}$.
Then for a positive integer $k$ with $k\ge 2$ and complex numbers 
$\mu_1,\ldots,\mu_k$
\begin{multline*}
 \Big(\prod_{\nu=1}^k(F_\pi-\mu_\nu)\Big)_{(i',j')(i,j)}
 \equiv
 \Big(\prod_{\nu=1}^{k-1}(F_\pi-\mu_\nu)\Big)_{(i',j')(i,j)}
 \bigl(\varpi_i - \mu_k
 + D_\pi(\varpi_i)\bigr)\\
 \mod U(\mathfrak g)(\mathfrak m_\Theta + \mathfrak n_\Theta) + 
       \!\!\!\!\!\! 
       \sum_{\substack{(s, t), (s'', t'');\\
       \varpi_{s}|_{\mathfrak a_\Theta}
       <_\Theta\varpi_{i}|_{\mathfrak a_\Theta} \\
       \varpi_{s, t} = \varpi_{s'', t''}}}
       \!\!\!\!\!\!
       \mathbb C\ %
       \Big(\prod_{\nu=1}^{k-1}
       (F_\pi-\mu_\nu)\Big)_{(s'',t'')(s,t)}.
\end{multline*}
\end{lem}
\begin{proof}
Note that $\varpi_{i,j}\equiv\varpi_i\mod U(\mathfrak g)\mathfrak m_\Theta$.
It follows from Lemma~\ref{lem:2} that
\[F_{(s, t)(i, j)} \equiv \delta_{si}\delta_{tj}\varpi_i
\mod U(\mathfrak g)(\mathfrak m_\Theta+\mathfrak n_\Theta)
\]
if $\varpi_{s,t} - \varpi_{i,j}\not\in
\Sigma(\mathfrak g)^-\setminus\Sigma(\mathfrak g_\Theta)$.

Put $F^\ell = \prod_{\nu=1}^\ell(F_\pi-\mu_\nu)$. Then Lemma~\ref{lem:2} implies
\begin{equation}\label{eq:zenka}
\begin{aligned}
{}   &F_{(i',j')(i,j)}^k
   - F_{(i',j')(i,j)}^{k-1}\bigl(\varpi_i - \mu_k\bigr)\\
   &\quad\equiv \sum_{\varpi_{s,t} - 
    \varpi_{i,j}\in\Sigma(\mathfrak g)^-\setminus\Sigma(\mathfrak g_\Theta)}
   [F_{(i',j')(s, t)}^{k-1}, F_{(s, t)(i,j)}]
    \mod U(\mathfrak g)(\mathfrak m_\Theta+\mathfrak n_\Theta)
    +\bar{\mathfrak n}_\Theta U(\mathfrak g)\\
   &\quad= \sum_{
      \substack{\alpha
      \in\Sigma(\mathfrak g)^+\setminus\Sigma(\mathfrak g_\Theta)\\
      \varpi_{s,t}=\varpi_{i,j}-\alpha}
      }
   \frac{\ang{E_{(s,t)(i,j)}}{X_{\alpha}}}
        {\ang{X_\alpha}{X_{-\alpha}}}
        [F_{(i',j')(s,t)}^{k-1}, X_{-\alpha}]\\
   &\quad=  \sum_{\substack
       {\alpha\in\Sigma(\mathfrak g)^+\setminus\Sigma(\mathfrak g_\Theta)\\
       \varpi_{s,t}=\varpi_{i,j}-\alpha\\
       \varpi_{s',t'}-\varpi_{s'',t''}=\alpha}
     }
   \frac{\ang{E_{(s,t)(i,j)}}{X_{\alpha}}
         \ang{E_{(s',t')(s'',t'')}}{X_{-\alpha}}}
        {\ang{X_\alpha}{X_{-\alpha}}}
  \\
    &\quad\quad
     \cdot\bigl(\delta_{ss''}\delta_{tt''}F_{(i',j')(s', t')}^{k-1}
      - \delta_{i's'}\delta_{j't'}F_{(s'',t'')(s,t)}^{k-1}\bigr)\\
    &\quad=
     \frac12\ang{\bar\pi-\varpi_i}{\bar\pi+\varpi_i - 2\rho}
     F_{(i',j')(i,j)}^{k-1}\\
    &\quad\quad - 
      \sum_{\substack
         {\alpha\in\Sigma(\mathfrak g)^+\setminus\Sigma(\mathfrak g_\Theta)\\
         \varpi_{s,t}=\varpi_{s'',t''}=\varpi_{i,j}-\alpha}
       }
       \frac{\ang{E_{(s,t)(i,j)}}{X_{\alpha}}
       \ang{E_{(i',j')(s'',t'')}}{X_{-\alpha}}}
       {\ang{X_\alpha}{X_{-\alpha}}}
        F_{(s'',t'')(s,t)}^{k-1}.
\end{aligned}
\end{equation}
In the above the second equality follows from \eqref{eq:F1} and Lemma~\ref{lem:proj} 
with $U=\mathfrak g$. 
The third equality follows from Lemma~\ref{lem:1} with
\[
 X_{-\alpha} = \sum_{\varpi_{s',t'}-\varpi_{s'',t''}=\alpha}
   \ang{E_{(s',t')(s'',t'')}}{X_{-\alpha}}E_{(s'',t'')(s',t')}
\]
which follows from the identification $\mathfrak g\subset\mathfrak{gl}_N$
together with the property of $\ang{\ }{\ }$.

Put $X^\vee=-{}^tX$ for $X\in M(N,\mathbb C)\simeq\mathfrak{gl}_N$.
Let $\{v_{i,j}^*\}$ be the dual base of $\{v_{i,j}\}$
and consider the contragredient representation $\pi^*$ of $\pi$.
Then $\pi^*(X)=X^\vee$ for $X\in\mathfrak g$ with respect to these basis.
Then $\ang{X}{Y}=\ang{X^\vee}{Y^\vee}$ for $X$, $Y\in\mathfrak g$ and
\[
 \sum_{\alpha\in\Sigma(\mathfrak g)^+\setminus\Sigma(\mathfrak g_\Theta)}
 \frac{X_{-\alpha}^\vee X_{\alpha}^\vee}{\ang{X_{-\alpha}^\vee}{X_{\alpha}^\vee}}
  v_{i,j}^*
 = \!\!\!
   \sum_{\substack
       {\alpha\in\Sigma(\mathfrak g)^+\setminus\Sigma(\mathfrak g_\Theta)\\
       \varpi_{s,t}=\varpi_{i,j}-\alpha\\
       \varpi_{s',t'}=\varpi_{i,j}}
     }
     \!\!\!\!
       \frac{\ang{E_{(s,t)(s',t')}}{X_{-\alpha}^\vee}
       \ang{E_{(i,j)(s,t)}}{X_{\alpha}^\vee}}
       {\ang{X_{-\alpha}^\vee}{X_{\alpha}^\vee}}
       v_{s',t'}^*,
\]
which is proved to be equal to $D_\pi(\varpi_i)v_{i,j}^*$ by 
Lemma~\ref{lem:casimir}~iii) because $(\bar\pi,\varpi_i,\rho)$ for $\pi$ changes into 
$(-\bar\pi,-\varpi_i,-\rho)$ in the dual $\pi^*$ with the reversed order of roots.
This implies the last equality in \eqref{eq:zenka}.

Note that if 
$D\in \bar{\mathfrak n}_\Theta U(\mathfrak g)
+U(\mathfrak g)(\mathfrak m_\Theta+\mathfrak n_\Theta)$ satisfies $[H,D]=0$
for all $H\in\mathfrak a_\Theta$, then 
$D\in U(\mathfrak g)(\mathfrak m_\Theta+\mathfrak n_\Theta)$.
Since the condition 
$\varpi_{i,j} - \varpi_{s,t}\in\Sigma(\mathfrak g)^+
\setminus\Sigma(\mathfrak g_\Theta)$ implies 
$\varpi_s|_{\mathfrak a_\Theta}<_\Theta\varpi_i|_{\mathfrak a_\Theta}$, 
we have the lemma.
\end{proof}

\begin{thm}\label{thm:min}
Retain the notation in\/ {\rm Definition~\ref{def:shift}}.
For $\varpi\in\mathfrak a^*$ we identify $\varpi|_{\mathfrak a_\Theta}$ with a linear
function on $\mathfrak a_\Theta^*$ by $\varpi|_{\mathfrak a_\Theta}(\lambda) = 
\ang{\lambda_\Theta}\varpi$ for $\lambda\in\mathfrak a_\Theta^*$. 
Put
\begin{equation}
\begin{aligned}
  \Omega_{\pi,\Theta} &= \{\bigl(\varpi|_{\mathfrak a_\Theta},D_\pi(\varpi)\bigr);\,
     \varpi\in\overline{\mathcal W}_\Theta(\pi)\},\\
  q_{\pi,\Theta}(x;\lambda) &= 
     \prod_{(\mu,C)\in\Omega_{\pi,\Theta}}\bigl(x - \mu(\lambda)-C\bigr).
\end{aligned}
\end{equation}
Then $q_{\pi,\Theta}(F_\pi;\lambda)M_\Theta(\lambda) = 0$ for any
$\lambda\in\mathfrak a_\Theta^*$.
\end{thm}
\begin{proof}
For any $D\in U(\mathfrak g)$ there exists a unique constant 
$T(D)\in\mathbb C$ satisfying
\[
  T(D)\equiv D\mod\bar{\mathfrak n}U(\mathfrak g) + J_\Theta(\lambda)
\]
because the dimension of the space
$M_\Theta(\lambda)/\bar{\mathfrak n}M_\Theta(\lambda)$ equals 1.
Notice that 
\[
J_\Theta(\lambda) = \sum_{H \in \mathfrak a_\Theta}
U(\mathfrak g)(H-\lambda(H)) + U(\mathfrak g)(\mathfrak m_\Theta + \mathfrak n_\Theta).
\]

Use the notation in Lemma~\ref{lem:zen}.
Since 
\[
 \ad(H)q_{\pi,\Theta}(F_\pi;\lambda)_{(i',j')(i,j)} = 
  (\varpi_{i',j'}-\varpi_{i,j})(H)q_{\pi,\Theta}(F_\pi;\lambda)_{(i',j')(i,j)}
 \quad\text{for \ }H\in\mathfrak a,
\]
$T(q_{\pi,\Theta}(F_\pi;\lambda)_{(i',j')(i,j)}) = 0$ if 
$\varpi_{i,j}\ne\varpi_{i',j'}$.  

Next assume $\varpi_{i,j}=\varpi_{i',j'}$ and put
\begin{align*}
  \Omega_{\pi,\Theta,i} &= \{(\mu,C) \in \Omega_{\pi,\Theta};\,
      \mu\le_\Theta \varpi_i|_{\mathfrak a_\Theta} \},\\
    q_{\pi,\Theta,i}(x;\lambda) &= 
     \prod_{(\mu,C)\in\Omega_{\pi,\Theta,i}}\bigl(x - \mu(\lambda)-C\bigr).
\end{align*}
Then $q(F_\pi)_{(i',j')(i,j)} \in J_\Theta(\lambda)$
for any $q(x) \in \mathbb C[x]$ which is a multiple of 
$q_{\pi,\Theta,i}(x;\lambda)$.
It is proved by the induction on 
$\varpi_i|_{\mathfrak a_\Theta}$
with the partial order $\le_\Theta$.
Take $i_0 \in \{ 1, \ldots, \kappa \}$ so that $\varpi_{i_0}=\bar\pi$.
If $i=i_0$
then Lemma~\ref{lem:2} and Lemma~\ref{lem:zen}
with $D_\pi(\varpi_i)=D_\pi(\bar\pi)=0$
imply our claim.
If $i \ne i_0$
then $\bar\pi|_{\mathfrak a_\Theta} <_\Theta \varpi_i|_{\mathfrak a_\Theta}$
and therefore
$\deg_x q_{\pi,\Theta,i}(x;\lambda) \ge 2$.
Hence we can use Lemma~\ref{lem:zen} again
to prove our claim inductively.

Thus we get the condition
\begin{equation}\label{eq:T_vanish}
T(q_{\pi,\Theta}(F_\pi;\lambda)_{(i',j')(i,j)})=0
\text{ for any }(i,j)\text{ and }(i',j').
\end{equation}
Let $\mathbf{V}(\lambda)$ denote the $\mathbf{C}$-subspace of $U(\mathfrak g)$
spanned by $q_{\pi,\Theta}(F_\pi;\lambda)_{(i',j')(i,j)}$.
Then $\mathbf{V}(\lambda)$ is $\ad(\mathfrak g)$-stable
by Lemma~\ref{lem:1}.
The $\mathfrak g$-module
\[
M_\lambda = \mathbf{V}(\lambda)M_\Theta(\lambda)
\]
is contained in $\bar{\mathfrak n}M_\Theta(\lambda)$
because
putting $u_\lambda=1 \bmod J_\Theta(\lambda)$,
\[
M_\lambda = \mathbf{V}(\lambda)U(\bar{\mathfrak n})u_\lambda
 = U(\bar{\mathfrak n})\mathbf{V}(\lambda)u_\lambda
 \subset U(\bar{\mathfrak n})\bar{\mathfrak n}U(\mathfrak g)u_\lambda
 =\bar{\mathfrak n}M_\Theta(\lambda).
\]

On the other hand, since $M_\Theta(\lambda)$ is irreducible if 
$\lambda$ belongs to a suitable open subset of 
$\mathfrak a_\Theta^*$, $M_\lambda=\{0\}$ in the open set.
If we fix a base
$\{Y_1,\ldots,Y_m\}$ of $\bar{\mathfrak n}_\Theta$, we have the unique
expression
\[
 q_{\pi,\Theta}(F_\pi;\lambda)_{(i',j')(i,j)}
 \equiv \sum_{\nu} Q_\nu(\lambda) Y_1^{\nu_1}\cdots Y_m^{\nu_m}
 \mod J_\Theta(\lambda)
\]
with polynomial functions $Q_\nu(\lambda)$.
All these $Q_\nu(\lambda)$ vanish on the open set and therefore they are
identically zero and we have $\mathbf{V}(\lambda)\subset J_\Theta(\lambda)$
for any $\lambda$.
We have then for any $\lambda$
\[
M_\lambda = \mathbf{V}(\lambda)U(\mathfrak g)u_\lambda
= U(\mathfrak g)\mathbf{V}(\lambda)u_\lambda = \{0\}.
\qedhere
\]
\end{proof}

Theorem~\ref{thm:min} is one of our central results
since $q_{\pi,\Theta}(x;\lambda)=q_{\pi,M_\Theta(\lambda)}(x)$
for a generic $\lambda\in{\mathfrak a}_\Theta^*$.
Before showing this minimality,
which will be done in Theorem~\ref{thm:minimality},
we mention the possibility of other approaches to Theorem~\ref{thm:min}.
In fact we have three different proofs.
The first one given above has the importance that the calculation
in the proof is also used in \S\ref{sec:ideal}
to study the properties of the two-sided ideal
of $U(\mathfrak g)$ generated by $q_{\pi,\Theta}(F_\pi;\lambda)_{ij}$.
The second one comes from a straight expansion of the method in \cite{Go1} and \cite{Go2}
to construct characteristic polynomials.
In the following we first discuss it.
The third one is based on {\it infinitesimal Mackey's tensor product theorem\/}
which we explain in Appendix~\ref{app:Mackey}. 
With this method we shall get the sufficient condition for 
the minimality of $q_{\pi,\Theta}(x;\lambda)$ (Theorem~\ref{thm:minimality})
and slightly strengthen the result of Theorem~\ref{thm:min} (Theorem~\ref{thm:minmore}).

\begin{defn}
Let $(\pi^*,V^*)$ be the contragredient representation of $(\pi,V)$
and $\{v^*_1,\ldots,v^*_N\}$ the dual base of the base $\{v_1,\ldots,v_N\}$ of $V$.
For a $\mathfrak g$-module $M$ 
define the homomorphism
\[h_{(\pi,M)}:
M(N,U(\mathfrak g))\rightarrow
\End\left(M \otimes V^*\right)\] 
of associative algebras by
\begin{equation}\label{eqn:h}
\bigl(
h_{(\pi,M)}(Q)
\bigr)
\bigl(\sum_{j=1}^N u_j\otimes v_j^*\bigr)
=\sum_{i=1}^N\sum_{j=1}^N
\bigl(Q_{ij}u_j\bigr)\otimes v_i^*
\end{equation}
for $u_j \in M$ and
$Q=\bigl(Q_{ij}\bigr) \in M(N,U(\mathfrak g))$.
Then $Q M=0$, namely, $Q_{ij}\in\Ann(M)$ for any $i, j$
if and only if $h_{(\pi,M)}(Q)=0$.
\end{defn}

The following lemma is considered in \cite{Go1} and \cite{Go2}.
\begin{lem}\label{lem:tcasimir}
Let $M$ be a $\mathfrak g$-module.
For an element $\sum_{j=1}^Nu_j\otimes v^*_j$ of
$M\otimes V^*$ with $u_j\in M$, we have
\[
 2 h_{(\pi,M)}(F_\pi)\bigl(\sum_{j=1}^Nu_j\otimes v^*_j\bigr)
 =  \sum_{j=1}^N\Delta_\pi(u_j)\otimes v^*_j
    + \sum_{j=1}^Nu_j\otimes \Delta_\pi(v^*_j)
    - \Delta_\pi\bigl(\sum_{j=1}^Nu_j\otimes v^*_j\bigr).
\]
In particular $h_{(\pi,M)}(F_\pi)\in\End_{\mathfrak g}(M\otimes V^*)$.
\end{lem}
\begin{proof}
Let $\{X_1,\ldots,X_\omega\}$ be a base of $\mathfrak g$ and let
$\{X_1^*,\ldots,X_\omega^*\}$ be its dual base with respect to $\ang{\ }{\ }$.
Then
\begin{align*}
 \sum_{j=1}^N&\Delta_\pi(u_j)\otimes v^*_j
    + \sum_{j=1}^Nu_j\otimes \Delta_\pi(v^*_j)
    - \Delta_\pi\bigl(\sum_{j=1}^Nu_j\otimes v^*_j\bigr)\\
 &= -\sum_{j=1}^N\sum_{\nu=1}^\omega X_\nu^*u_j\otimes X_\nu v^*_j
 - \sum_{j=1}^N\sum_{\nu=1}^\omega X_\nu u_j\otimes X_\nu^* v^*_j\\
 &= \sum_{j=1}^N\sum_{\nu=1}^\omega\Bigl(X_\nu^*u_j\otimes %
     \sum_{i=1}^N\ang{X_\nu}{E_{ij}}v^*_i
   + X_\nu u_j\otimes
     \sum_{i=1}^N\ang{X_\nu^*}{E_{ij}} v^*_i\Bigr)\\
 &= 2\sum_{i=1}^N\sum_{j=1}^N(p(E_{ij})u_j)\otimes v^*_i.
\end{align*}
Here we use the fact that $X v_j^* = -\sum_{i=1}^N\ang{X}{E_{ij}}v_i^*$ for 
$X\in\mathfrak g$
because  $X v_j = \sum_{i=1}^N\ang{X}{E_{ji}}v_i$.
\end{proof}

Now we examine the tensor product $M\otimes V^*$ in the preceding lemma when
$M$ is realized as a finite dimensional quotient of a generalized Verma module
$M_\Theta(\lambda)$.

\begin{prop}[a character identity for a tensor product]\label{prop:tensor}
Put
\[
 \chi_\Lambda = 
 \frac{\sum_{w\in W}\sign(w)e^{w(\Lambda+\rho)}}
  {\prod_{\alpha\in\Sigma(\mathfrak g)^+}(e^{\frac\alpha2}-e^{-\frac\alpha2})}
\]
for $\Lambda\in\mathfrak a^*$.
If $\ang{\Lambda}{\alpha}=0$ for any $\alpha\in\Theta$, then
\begin{equation}\label{eq:decomp}
 \chi_{\pi^*}\chi_{\Lambda} =
 \sum_{\varpi\in\mathcal W(\pi^*)}m_{\pi^*,\Theta}(\varpi)\chi_{\Lambda+\varpi}
\end{equation}
by denoting
\[
 m_{\pi^*,\Theta}(\varpi) = \dim\{v^*\in V^*;\,Hv^*=\varpi(H)v^*\ (\forall H\in\mathfrak a),\ 
 Xv^*=0\ (\forall X\in\mathfrak g_\Theta\cap\mathfrak n)\}.
\]
Here $\chi_{\pi^*}$ is the character of the representation $(\pi^*,V^*)$ and
for $\mu\in\mathfrak a^*$,
$e^\mu$ denotes the function on $\mathfrak a$ which takes the value
$e^{\mu(H)}$ at $H\in\mathfrak a$.
\end{prop}
\begin{proof}
It is sufficient to prove \eqref{eq:decomp} under the condition that
$\ang{\Lambda}{\alpha}$ is a sufficiently large real number for any 
$\alpha\in\Psi(\mathfrak g)\setminus\Theta$ because both hand sides of
\eqref{eq:decomp} are holomorphic with respect to $\Lambda\in\mathfrak a^*$.
Put
\begin{align*}
 \mathfrak a^*_0 &= \{\mu\in\mathfrak a^*;\, \ang{\mu}{\alpha}\in
  \mathbb R\quad\bigl(\forall\alpha\in\Sigma(\mathfrak g)\bigr)\},\\
 \mathfrak a^*_+ &= \{\mu\in\mathfrak a^*_0;\, \ang{\mu}{\alpha}\ge 0\
  \quad\bigl(\forall\alpha\in\Sigma(\mathfrak g)^+\bigr)\},\\
 \chi_\Lambda^+ &= \frac{\sum_{w\in W_\Theta}\sign(w)e^{w(\Lambda+\rho)}}
  {\prod_{\alpha\in\Sigma(\mathfrak g)^+}(e^{\frac\alpha2}-e^{-\frac\alpha2})},\\
 \bar\chi_\varpi &= \frac{\sum_{w'\in W_\Theta}\sign(w')e^{w'(\varpi+\rho(\Theta))}}
  {\prod_{\alpha\in\Sigma(\mathfrak g_\Theta)^+}(e^{\frac\alpha2}-e^{-\frac\alpha2})}.
\end{align*}
Then $\chi_{\pi^*} = \sum_{\varpi\in\mathcal W(\pi^*)}m_{\pi^*,\Theta}(\varpi)\bar\chi_\varpi$ by
Weyl's character formula and if $\varpi\in\mathcal W(\pi^*)$ satisfies 
$m_{\pi^*,\Theta}(\varpi) > 0$, then $\Lambda+\varpi\in\mathfrak a^*_+$ and
\begin{align*}
  \bar\chi_\varpi\chi_\Lambda^+ {\prod_{\alpha\in\Sigma(\mathfrak g)^+}(e^{\frac\alpha2}-e^{-\frac\alpha2})} &= 
  \frac{\sum_{w\in W_\Theta}\sign(w)e^{w(\varpi+\rho(\Theta))}}
  {\prod_{\alpha\in\Sigma(\mathfrak g_\Theta)^+}(e^{\frac\alpha2}-e^{-\frac\alpha2})}
  e^{\Lambda+\rho_\Theta}\sum_{w'\in W_\Theta}\sign(w')e^{w'\rho(\Theta)}\\
 &= \sum_{w\in W_\Theta}\sign(w)e^{w(\Lambda+\varpi+\rho)}\\
 &\equiv
  e^{\Lambda+\varpi+\rho}
  \mod\sum_{\mu\in\mathfrak a^*_0\setminus\mathfrak a^*_+}\mathbb Ze^\mu.
\end{align*}
For any $w\in W\setminus W_\Theta$ there exists 
$\alpha\in\Sigma(\mathfrak g)^-\setminus\Sigma(\mathfrak g_\Theta)$ with
$w\alpha\in\Sigma(\mathfrak g)^+$ and then the value
$-\ang{w(\Lambda+\rho)}{w\alpha}=-\ang{(\Lambda+\rho)}{\alpha}$ 
is sufficiently large and 
therefore \[\bar\chi_\varpi(\chi_\Lambda-\chi_\Lambda^+)
{\prod_{\alpha\in\Sigma(\mathfrak g)^+}(e^{\frac\alpha2}-e^{-\frac\alpha2})} \in
\sum_{\mu\in\mathfrak a^*_0\setminus\mathfrak a^*_+}\mathbb Ze^\mu.\]
Hence
\[
 \chi_{\pi^*}\chi_\Lambda 
 \prod_{\alpha\in\Sigma(\mathfrak g)^+}(e^{\frac\alpha2}-e^{-\frac\alpha2})
 \equiv\sum_{\varpi\in\mathcal W_\Theta(\pi^*)}
 m_{\pi^*,\Theta}(\varpi)e^{\Lambda+\varpi+\rho}
 \mod\sum_{\mu\in\mathfrak a^*_0\setminus\mathfrak a^*_+}\mathbb Ze^\mu
\]
and we have the proposition because
$\chi_{\pi^*}\chi_\Lambda 
\prod_{\alpha\in\Sigma(\mathfrak g)^+}(e^{\frac\alpha2}-e^{-\frac\alpha2})$
is an odd function under $W$.
\end{proof}

\begin{lem}[eigenvalue]\label{lem:eigen}
Let $(\pi_\Lambda,V_\Lambda)$ be an irreducible finite dimensional 
representation of $\mathfrak g$ with highest weight $\Lambda$.
Suppose $\ang{\Lambda}{\alpha}=0$ for $\alpha\in\Theta$ and
$\ang{\Lambda+\varpi}{\alpha}\ge0$ for $\varpi\in\mathcal W_\Theta(\pi^*)$
and $\alpha\in\Psi(\mathfrak g)\setminus\Theta$.
Then %
the set of the eigenvalues 
of $h_{\pi, V_\Lambda}(F_\pi)\in\End(V_\Lambda\otimes V^*)$ without counting their 
multiplicities equals
\begin{align*}
&\{-\ang{\Lambda}{\varpi}+\frac12\ang{\pi^*-\varpi}{\pi^*+\varpi+2\rho};\,
\varpi\in\mathcal W_\Theta(\pi^*)\}\\
&=\{\ang{\Lambda}{\varpi}+\frac12\ang{\bar\pi-\varpi}{\bar\pi+\varpi-2\rho};\,
\varpi\in\overline{\mathcal W}_\Theta(\pi)\}.
\end{align*}
Here we identify $\pi^*$ with the highest weight of $(\pi^*,V^*)$.
\end{lem}
\begin{proof}
The assumption of the lemma and Proposition~\ref{prop:tensor} imply
$$\pi^*\otimes\pi_\Lambda = 
  \sum_{\varpi\in\mathcal W_\Theta(\pi^*)}m_{\pi^*,\Theta}(\varpi)\pi_{\Lambda+\varpi}$$
and hence by Lemma~\ref{lem:casimir} ii) and Lemma~\ref{lem:tcasimir} the eigenvalues of
$2h_{\pi, V_\Lambda}(F_\pi)$ are
\[
\ang{\Lambda}{\Lambda+2\rho}
+\ang{\pi^*}{\pi^*+2\rho}-\ang{\Lambda+\varpi}{\Lambda+\varpi+2\rho}
=-2\ang{\Lambda}{\varpi}+\ang{\pi^*-\varpi}{\pi^*+\varpi+2\rho}
\]
with $\varpi\in\mathcal W_\Theta(\pi^*)$.
Since $\overline{\mathcal W}_\Theta(\pi)= -\mathcal W_\Theta(\pi^*)$, we have
the lemma.
\end{proof}

\begin{proof}[Proof of\/ {\rm Theorem \ref{thm:min}} -- the 2nd version]
This proof differs from the previous one in how to deduce
the condition~\eqref{eq:T_vanish}.
The rests of two proofs are the same.

Note that for fixed $(i,j)$ and $(i',j')$ the value $T(q_{\pi,\Theta}(F_\pi;\lambda)_{(i',j')(i,j)})$
depends algebraically on the parameter $\lambda\in{\mathfrak a}_\Theta^*$.
Since the set
$$S=\{\lambda\in{\mathfrak a}_\Theta^*;\,
\ang{\lambda_\Theta+\varpi}{\alpha}\in\{0,1,2,\ldots \}\text{ for }
\varpi\in\mathcal W_\Theta(\pi^*)\cup\{0\}\text{ and }
\alpha\in\Psi(\mathfrak g)\setminus\Theta\}$$
is Zariski dense in ${\mathfrak a}_\Theta^*$,
we have only to show \eqref{eq:T_vanish} for $\lambda\in S$.
In this case we have
from Lemma~\ref{lem:eigen} and the definition of $q_{\pi,\Theta}(x;\lambda)$,
\[
h_{\pi, V_{\lambda_\Theta}}(q_{\pi,\Theta}(F_\pi;\lambda))
=q_{\pi,\Theta}(h_{\pi, V_{\lambda_\Theta}}(F_\pi);\lambda)
=0.
\]
Hence $q_{\pi,\Theta}(F_\pi;\lambda)_{(i',j')(i,j)}\in\Ann(V_{\lambda_\Theta})$
for any $(i,j)$ and $(i',j')$. 
On the other hand,
if we take a highest weight vector $v_\lambda$ of $V_{\lambda_\Theta}$,
we get
\[
q_{\pi,\Theta}(F_\pi;\lambda)_{(i',j')(i,j)}v_\lambda
\in T(q_{\pi,\Theta}(F_\pi;\lambda)_{(i',j')(i,j)})v_\lambda
+\bar{\mathfrak n}V_{\lambda_\Theta}
\] 
and therefore $T(q_{\pi,\Theta}(F_\pi;\lambda)_{(i',j')(i,j)})=0$.
\end{proof}

\begin{thm}[minimality]\label{thm:minimality}
Let $\lambda\in\mathfrak{a}_\Theta^*$.

\noindent
{\rm i)}
The set of the roots of $q_{\pi, M_\Theta(\lambda)}(x)$
equals 
$\{ \ang{\lambda_\Theta}{\varpi}+D_\pi(\varpi);\,
\varpi\in\overline{\mathcal W}_\Theta(\pi) \}.$

\noindent
{\rm ii)}
If each root of $q_{\pi,\Theta}(x;\lambda)$ is simple,
then $q_{\pi,\Theta}(x;\lambda)=q_{\pi, M_\Theta(\lambda)}(x)$. 
Hence we call $q_{\pi,\Theta}(x;\lambda)$ the {\it global minimal polynomial\/}
of the pair $(\pi,M_\Theta(\lambda))$.
\end{thm}
\begin{proof}
i)
Fix an irreducible decomposition $\bigoplus_{i=1}^\kappa U_i$ of 
the $\mathfrak g_\Theta$-module $V^*|_{\mathfrak g_\Theta}$.
Let $\varpi_i\in\mathfrak a^*$ be the highest weight of $U_i$.
With a suitable change of indices
we may assume 
$\varpi_i|_{\mathfrak a_\Theta} <_\Theta \varpi_j|_{\mathfrak a_\Theta}$
implies $i > j$.
Then putting $V_i=\bigoplus_{\nu=1}^iU_\nu$
we get a $\mathfrak p_\Theta$-stable filtration 
\[ \{ 0 \} = V_0 \subsetneq V_1
\subsetneq \cdots \subsetneq V_\kappa = V^*|_{\mathfrak p_\Theta}.
\]
Note that $V_i / V_{i-1} \simeq U_i$
is an irreducible $\mathfrak p_\Theta$-module
on which $\mathfrak n_\Theta$ acts trivially.

Recall $M_\Theta(\lambda) \simeq M_{(\Theta, \lambda_\Theta)}
= U(\mathfrak g) \otimes_{U(\mathfrak p_\Theta)} U_{(\Theta, \lambda_\Theta)}$
and $\dim U_{(\Theta, \lambda_\Theta)}=1$.
Hence writing $\mathbb C_\lambda$ instead of $U_{(\Theta, \lambda_\Theta)}$
we get by Theorem~\ref{thm:Mackey} of Appendix~\ref{app:Mackey}
\[
M_\Theta(\lambda)\otimes V^*
= \left( U(\mathfrak g) \otimes_{U(\mathfrak p_\Theta)} \mathbb C_\lambda \right)
\otimes V^*
\simeq U(\mathfrak g) \otimes_{U(\mathfrak p_\Theta)}
\left(\mathbb C_\lambda \otimes V^*|_{\mathfrak p_\Theta} \right).
\]
Since $\mathbb C_\lambda\otimes_{\mathbb C}\cdot$
and
$U(\mathfrak g) \otimes_{U(\mathfrak p_\Theta)} \cdot
= U(\bar{\mathfrak n}_\Theta) \otimes_{\mathbb C} \cdot$
are exact functors,
putting $M_i=U(\mathfrak g) \otimes_{U(\mathfrak p_\Theta)}
\left(\mathbb C_\lambda \otimes V_i \right)$
we get a $\mathfrak g$-stable filtration 
\[
\{ 0 \} = M_0 \subsetneq M_1
\subsetneq \cdots \subsetneq M_\kappa = M_\Theta(\lambda) \otimes V^*
\]
with
\begin{equation}\label{eqn:filter_q}
M_i/M_{i-1}\simeq U(\mathfrak g)\otimes_{U(\mathfrak p_\Theta)}
\left(\mathbb C_\lambda \otimes U_i \right)=M_{(\Theta,\lambda_\Theta+\varpi_i)}.
\end{equation}

Now as a subalgebra of $\End\left(M_\Theta(\lambda) \otimes V^*\right)$
we take 
\[A=
\{D;\,
D M_i \subset M_i\ \text{for}\ i=1, \ldots, \kappa\}.
\]
Then by Lemma~\ref{lem:tcasimir} and Lemma~\ref{lem:casimir}~ii)
we have $h_{(\pi,M_\Theta(\lambda))}
\left( q(F_\pi) \right) \in A$
for any polynomial $q(x)\in\mathbb C[x]$.
Let $\eta_i:A\rightarrow
\End\left(M_i/M_{i-1}\right)
\simeq\End\left(M_{(\Theta,\lambda_\Theta+\varpi_i)}\right)$
be a natural algebra homomorphism.
Then using Lemma~\ref{lem:tcasimir} and Lemma~\ref{lem:casimir}~ii)
again we get
\begin{equation}\label{eqn:M_i_scalar}
\begin{aligned}
\eta_i\bigl( &h_{(\pi,M_\Theta(\lambda))}(F_\pi)\bigr) \\
&= \frac12\ang{\lambda_\Theta}{\lambda_\Theta+2\rho}
+ \frac12\ang{-\bar\pi}{-\bar\pi+2\rho}
- \frac12\ang{\lambda_\Theta+\varpi_i}{\lambda_\Theta+\varpi_i+2\rho}\\
&= \ang{\lambda_\Theta}{-\varpi_i}+D_\pi(-\varpi_i).
\end{aligned}
\end{equation}
and therefore
\[
\begin{aligned}
q_{\pi,M_\Theta(\lambda)}\left(\ang{\lambda_\Theta}{-\varpi_i}+D_\pi(-\varpi_i)\right) &=
q_{\pi,M_\Theta(\lambda)}\left(\eta_i\left( h_{(\pi,M_\Theta(\lambda))}(F_\pi)\right)\right) \\
&= \eta_i\left( h_{(\pi,M_\Theta(\lambda))}\left(q_{\pi,M_\Theta(\lambda)}(F_\pi)\right)\right) \\
&=0.
\end{aligned}
\]
Since $\{\varpi_i\}=\mathcal W_\Theta(\pi^*)=-\overline{\mathcal W}_\Theta(\pi)$
we can conclude $\ang{\lambda_\Theta}{\varpi}+D_\pi(\varpi)$ is 
a root of the minimal polynomial for each $\varpi\in\overline{\mathcal W}_\Theta(\pi)$.
Conversely Theorem~\ref{thm:min} assures any other roots do not exist.

ii)
The claim immediately follows from i) and the definition of $q_{\pi,\Theta}(x;\lambda)$.
\end{proof}

\begin{rem}
In general 
it may happen for a certain $\lambda$
that $q_{\pi,\Theta}(x;\lambda)\ne q_{\pi, M_\Theta(\lambda)}(x)$.
Such example is
shown in \cite{O-Cl} when $\mathfrak g$ is $\mathfrak o_{2n}$ and 
$\lambda$ is invariant under an outer automorphism 
of $\mathfrak g$, which is related to the following theorem.
It gives more precise information
on our minimal polynomials.
\end{rem}

\begin{thm}\label{thm:minmore}
Let $\lambda\in\mathfrak{a}_\Theta^*$.
Let $\overline{\mathcal W}_\Theta(\pi)=
\overline{\mathcal W}_\lambda^1\sqcup\overline{\mathcal W}_\lambda^2\cdots
\sqcup\overline{\mathcal W}_\lambda^{m_\lambda}$
be a division of $\overline{\mathcal W}_\Theta(\pi)$
into non-empty subsets $\overline{\mathcal W}_\lambda^\ell$  
such that
the relation
$\lambda_\Theta-\varpi\in\{w.(\lambda_\Theta-\varpi');\,w\in W\}$
holds for $\varpi, \varpi' \in\overline{\mathcal W}_\Theta(\pi)$
if and only if $\varpi, \varpi' \in \overline{\mathcal W}_\lambda^\ell$ for some $\ell$.
For each $\ell$
we denote by $\kappa_\ell$
the maximal length of 
sequences $\{\varpi, \varpi', \ldots, \varpi'' \}$
of weights in $\overline{\mathcal W}_\lambda^\ell$
such that the restriction of each weight to ${\mathfrak{a}_\Theta}$ gives both strictly and linearly ordered sequences:
$$\varpi|_{\mathfrak{a}_\Theta} <_\Theta \varpi'|_{\mathfrak{a}_\Theta} <_\Theta \cdots <_\Theta \varpi''|_{\mathfrak{a}_\Theta}.$$

\noindent
{\rm i)}
$\ang{\lambda_\Theta}{\varpi}+D_\pi(\varpi)=\ang{\lambda_\Theta}{\varpi'}+D_\pi(\varpi')$
if $\varpi, \varpi' \in \overline{\mathcal W}_\lambda^\ell$ for some $\ell$.

\noindent
{\rm ii)}
Let $q(x)\in\mathbb{C}[x]$
and suppose for each $\ell=1,\ldots,m_\lambda$,
$q(x)$ is a multiple of 
$(x-\ang{\lambda_\Theta}{\varpi}-D_\pi(\varpi))^{\kappa_\ell}$
with $\varpi\in\overline{\mathcal W}_\lambda^\ell$.
Then $q(F_\pi)M_\Theta(\lambda)=0$.
\end{thm}

\begin{proof}
i)
By the $W$-invariance of $\ang{\ }{\ }$ and the assumption,
we have
$$\ang{\lambda_\Theta+\rho-\varpi}{\lambda_\Theta+\rho-\varpi}
=\ang{\lambda_\Theta+\rho-\varpi'}{\lambda_\Theta+\rho-\varpi'},$$
which implies the claim.

ii)
Use the notation in the proof of Theorem~\ref{thm:minimality}.
Let $M$ be a $\mathfrak g$-module
and $\mu\in\mathfrak a^*$.
We say that a non-zero vector $v$ in $M$
is a {\it generalized weight vector for the generalized infinitesimal character\/} $\mu$ 
if for any $\Delta\in Z(\mathfrak g)$
there exists a positive integer $k$ such that
$(\Delta-\Delta_{\mathfrak a}(\mu))^kv=0$.
We denote by $(M)_{(\mu)}$
the submodule of $M$
spanned by the generalized weight vectors for the generalized infinitesimal character $\mu$.
Note that $(M)_{(\mu)}=(M)_{(\mu')}$
if and only if $\mu=w.\mu'$ for some $w\in W$.
By virtue of \eqref{eqn:filter_q} and Remark~\ref{rem:inf_char},
$M_\Theta(\lambda)\otimes V^*$ is uniquely decomposed
as a direct sum of submodules in
$\{(M_\Theta(\lambda)\otimes V^*)_{(\lambda_\Theta+\varpi_\nu)};\,\nu=1, \ldots, i\}$.

For $i=1,\ldots,\kappa$
using a $\mathfrak{p}_\Theta$-module
\[
V_{[i]}=U_i \oplus 
\bigoplus_{\nu;\,\varpi_i|_{\mathfrak{a}_\Theta}<_\Theta \varpi_\nu|_{\mathfrak{a}_\Theta}} U_\nu
\subset V_i,
\]
define 
\[
M_{[i]} = U(\mathfrak g)\otimes_{\mathfrak{p}_\Theta}(\mathbb{C}_\lambda\otimes V_{[i]})
= U(\bar{\mathfrak{n}}_\Theta)\otimes\mathbb{C}_\lambda\otimes V_{[i]}.
\]
It is naturally considered as a $\mathfrak{g}$-submodule of 
$M_i=U(\bar{\mathfrak{n}}_\Theta)\otimes\mathbb{C}_\lambda\otimes V_i$.
If we define the surjective homomorphism
\[
\tau_{[i]} : M_{[i]}\hookrightarrow M_i \rightarrow M_i/M_{i-1} \simeq M_{(\Theta,\lambda_\Theta+\varpi_i)},
\]
then
\begin{equation}\label{eq:M_kernel}
\Ker \tau_{[i]} = 
\sum_{\nu;\,\varpi_i|_{\mathfrak{a}_\Theta}<_\Theta \varpi_\nu|_{\mathfrak{a}_\Theta}}
M_{[\nu]}.
\end{equation}
Since $M_{(\Theta,\lambda_\Theta+\varpi_i)}$ has
infinitesimal character $\lambda_\Theta+\varpi_i$
we get
\[
M_{[i]}=(M_{[i]})_{(\lambda_\Theta+\varpi_i)}
+ \sum_{\nu;\,\varpi_i|_{\mathfrak{a}_\Theta}<_\Theta \varpi_\nu|_{\mathfrak{a}_\Theta}}
M_{[\nu]}.
\]
Therefore we get inductively
\begin{equation}\label{eq:geninfdcp}
M_{[i]}=(M_{[i]})_{(\lambda_\Theta+\varpi_i)}
+ \sum_{\nu;\,\varpi_i|_{\mathfrak{a}_\Theta}<_\Theta \varpi_\nu|_{\mathfrak{a}_\Theta}}
(M_{[\nu]})_{(\lambda_\Theta+\varpi_\nu)}.
\end{equation}

Notice that
the $\mathfrak{g}$-homomorphism $h_{(\pi,M_\Theta(\lambda))}(F_\pi)$
leaves any $\mathfrak{g}$-submodule of $M_\Theta(\lambda)\otimes V^*$ stable.
Then from \eqref{eqn:M_i_scalar} and \eqref{eq:M_kernel}
\[\begin{aligned}
\Bigl(h_{(\pi,M_\Theta(\lambda))}(F_\pi)-\ang{\lambda_\Theta}{-\varpi_i}&-D_\pi(-\varpi_i)\Bigr) 
(M_{[i]})_{(\lambda_\Theta+\varpi_i)} \\
&\subset\left(
\sum_{\nu;\,\varpi_i|_{\mathfrak{a}_\Theta}<_\Theta \varpi_\nu|_{\mathfrak{a}_\Theta}}
\!\!\!\!\! M_{[\nu]}
\right)_{(\lambda_\Theta+\varpi_i)} \\
&=\left(
\sum_{\nu;\,\varpi_i|_{\mathfrak{a}_\Theta}<_\Theta \varpi_\nu|_{\mathfrak{a}_\Theta}} 
\!\!\!\!\! (M_{[\nu]})_{(\lambda_\Theta+\varpi_\nu)}
\right)_{(\lambda_\Theta+\varpi_i)} \\
&=
\sum_{\substack{\nu;\,
\varpi_i|_{\mathfrak{a}_\Theta}<_\Theta \varpi_\nu|_{\mathfrak{a}_\Theta}, \cr
\lambda_\Theta+\varpi_\nu
\in \{ w.(\lambda_\Theta+\varpi_i);\, w\in W \}
}}
\!\!\!\!\! (M_{[\nu]})_{(\lambda_\Theta+\varpi_\nu)}.
\end{aligned}\]
By the relation $\{\varpi_i\}=\mathcal W_\Theta(\pi^*)=-\overline{\mathcal W}_\Theta(\pi)$
and the assumption of ii)
we get inductively
\[
h_{(\pi,M_\Theta(\lambda))}(q(F_\pi))(M_{[i]})_{(\lambda_\Theta+\varpi_i)}=
q(h_{(\pi,M_\Theta(\lambda))}(F_\pi))(M_{[i]})_{(\lambda_\Theta+\varpi_i)}=\{0\} 
\]
for $i=1,\ldots,\kappa$.
Now our claim is clear because by \eqref{eq:geninfdcp} we have
\[
M_\Theta(\lambda)\otimes V^*=\sum_{i=1}^\kappa M_{[i]}
=\sum_{i=1}^\kappa (M_{[i]})_{(\lambda_\Theta+\varpi_i)}.\qedhere
\]
\end{proof}

\begin{cor}\label{cor:involution}
Let $\tau$ be an involutive automorphism of $\mathfrak g$ which corresponds to
an automorphism of the Dynkin diagram of $\mathfrak g$.
Then $\tau(\mathfrak a) = \mathfrak a$ and $\tau(\mathfrak n)=\mathfrak n$.
Furthermore we suppose 
$\tau(\mathfrak{p}_\Theta)=\mathfrak{p}_\Theta$,
or equivalently, $\tau(\mathfrak{a}_\Theta)=\mathfrak{a}_\Theta$.
For $\varpi\in\mathfrak{a}^*$
we identify $\varpi|_{(\mathfrak{a}_\Theta)^\tau}$
as a linear function on $(\mathfrak a_\Theta^*)^\tau$
by $\varpi|_{(\mathfrak{a}_\Theta)^\tau}(\lambda)=\ang{\lambda_\Theta}{\varpi}$
for $\lambda\in(\mathfrak a_\Theta^*)^\tau$.
Put
\[
\begin{aligned}
  \Omega_{\pi,\Theta,\tau} &= \{\bigl(\varpi|_{(\mathfrak a_\Theta)^\tau},D_\pi(\varpi)\bigr);\,
     \varpi\in\overline{\mathcal W}_\Theta(\pi)\},\\
  q_{\pi,\Theta,\tau}(x;\lambda) &= 
     \prod_{(\mu,C)\in\Omega_{\pi,\Theta,\tau}}\bigl(x - \mu(\lambda)-C\bigr).
\end{aligned}
\]
Then for $\lambda\in\mathfrak (a_\Theta^*)^\tau$
we have the following.

\noindent
{\rm i)}
$q_{\pi,\Theta,\tau}(F_\pi;\lambda)M_\Theta(\lambda) = 0$.

\noindent
{\rm ii)}
If each root of $q_{\pi,\Theta,\tau}(x;\lambda)$ is simple,
then $q_{\pi,\Theta,\tau}(x;\lambda)=q_{\pi,M_\Theta(\lambda)}(x)$. 
\end{cor}

\begin{proof}
We naturally identify
$\rho_\Theta$ with an element in
$(\mathfrak{a}_\Theta^*)^\tau$.
For a given pair of weights $\varpi, \varpi' \in \overline{\mathcal W}_\Theta(\pi)$ with 
$\varpi|_{\mathfrak a_\Theta} <_\Theta \varpi'|_{\mathfrak a_\Theta}$,
choose the non-negative integers $\{m_\alpha;\,\alpha\in\Psi(\mathfrak g)\setminus\Theta\}$
so that $\varpi'|_{\mathfrak a_\Theta}-\varpi|_{\mathfrak a_\Theta}
=\sum_{\alpha\in\Psi(\mathfrak g)\setminus\Theta}m_\alpha \alpha|_{\mathfrak a_\Theta}$.
Then $\varpi'|_{\mathfrak a_\Theta}(\rho_\Theta)
-\varpi|_{\mathfrak a_\Theta}(\rho_\Theta)=
\sum_{\alpha\in\Psi(\mathfrak g)\setminus\Theta}m_\alpha 
\ang{\alpha}{\rho_\Theta}>0$.
It simply shows
\[
\bigl(\varpi|_{(\mathfrak a_\Theta)^\tau},D_\pi(\varpi)\bigr)
\ne
\bigl(\varpi'|_{(\mathfrak a_\Theta)^\tau},D_\pi(\varpi')\bigr).
\]
Hence from Theorem~\ref{thm:minmore} we get i).
Now ii) is clear from Theorem~\ref{thm:minimality}.
\end{proof}

We will shift $\mathfrak a^*$ by $\rho$ so that
the action $w.\mu=w(\mu+\rho) - \rho$ for $\mu\in\mathfrak a^*$ and 
$w\in W$ changes into the natural action of $W$ and then we can give the 
characteristic polynomial
as a special case of the global minimal polynomials.
The result itself is not new and it has already been studied in \cite{Go2}.
\begin{thm}[Cayley-Hamilton \cite{Go2}]\label{thm:char}
The characteristic polynomial $q_\pi(x)$ of $\pi$ is given by
\begin{equation}\label{eq:char}
  q_\pi(x) = \prod_{\varpi\in\mathcal W(\pi)}
  \Bigl(
    x - \varpi - 
   \frac{\ang\pi{\pi + 2\rho} - \ang\varpi\varpi}2
  \Bigr)
\end{equation}
under the identification %
$\mathbb C[x]\otimes S(\mathfrak a^*)^W\simeq
 \mathbb C[x]\otimes S(\mathfrak a)^W\simeq
Z(\mathfrak g)[x]$
by the symmetric bilinear form $\ang{\ }{\ }$ and
the Harish-Chandra isomorphism:
\begin{align*}
 Z(\mathfrak g)\simeq U(\mathfrak a)^{W};\
 &\Delta\mapsto \Upsilon(\Delta), \\
 &\Upsilon(\Delta)(\mu) = \Delta_{\mathfrak{a}}(\mu-\rho)
 \text{ for }\mu\in\mathfrak{a}^*.
\end{align*}
Here $\pi$ is identified with its highest weight.
In particular $q_\pi(x)\in Z(\mathfrak g)[x]$.
\end{thm}
\begin{proof}
Note that $\ang{\pi}{\pi+2\rho}=\ang{\bar\pi}{\bar\pi-2\rho}$.
Let $\tilde{q}_\pi(x)$  be the element
of $Z(\mathfrak g)[x]$ identified with the right-hand side of \eqref{eq:char}.
Put $\mathbf{V}=\sum_{i,\,j}\mathbb C\tilde q_\pi(F_\pi)_{ij}$ and 
$\mathbf{V}_\mathfrak a=\{D_\mathfrak a;\, D\in \mathbf{V}\}$.
Then Theorem~\ref{thm:min} with $\Theta=\emptyset$ shows $Q(\mu) = 0$ for any
$\mu\in\mathfrak a^*$ and $Q\in \mathbf{V}_\mathfrak a$, 
which implies $\mathbf{V}_{\mathfrak a}=\{0\}$.
Since $\mathbf{V}$ is $\ad(\mathfrak g)$-stable, we have $\mathbf{V}=\{0\}$ as is shown in 
\cite[Lemma~2.12]{O-BV}. 
Since the minimality of $\tilde q_\pi(x)$ follows from Theorem~\ref{thm:minimality},
we get $q_\pi(x)=\tilde q_\pi(x)$.
\end{proof}

\begin{cor}
{\rm i)}
Let $\mathfrak g$ be a simple Lie algebra.
Then the characteristic polynomial of the adjoint representation of 
$\mathfrak g$ is given by
\[
  q_{\alpha_{\max}}(x) = %
  \prod_{\alpha\in\Sigma(\mathfrak g)\cup\{0\}}
  \Bigl(x -  \alpha - 
   \frac{1- B(\alpha,\alpha)}2\Bigr).
\]
Here $B(\ , \ )$ denotes the Killing form of $\mathfrak g$.

{\rm ii)}
Suppose that the representation $\pi$ is {\it minuscule}, that is,
$\mathcal W(\pi)$ is a single $W$-orbit.
Then
\[
  q_\pi(x) = \prod_{\varpi\in\mathcal W(\pi)}
   \bigl(x-\varpi- \ang{\pi}{\rho}\bigr).
\]
\end{cor}
\begin{proof}
This is a direct consequence of Theorem~\ref{thm:char} and 
Lemma~\ref{lem:casimir}~v).
\end{proof}

\begin{cor}
Put $q_\pi(x) = x^m + \Delta_1x^{m-1} + \cdots + \Delta_{m-1}x + \Delta_m$
with $\Delta_j\in Z(\mathfrak g)$ and define
\[
\tilde F_\pi = -F_\pi^{m-1} - \Delta_1F_\pi^{m-2} - \cdots - \Delta_{m-1}I_N.
\]
Then
\[
    F_\pi\tilde F_\pi = \tilde F_\pi F_\pi = \Delta_mI_N = \prod_{\varpi\in\mathcal W(\pi)}
    \Bigl(-\varpi -\frac{\ang\pi{\pi + 2\rho} - \ang\varpi\varpi}2\Bigr)I_N,
\]
In particular, $F_\pi$ is invertible in 
$M\bigl(N,\hat Z(\mathfrak g)\otimes_{Z(\mathfrak g)} U(\mathfrak g)\bigr)$
with the quotient field $\hat Z(\mathfrak g)$ of $Z(\mathfrak g)$.
\end{cor}

In the next definition and the subsequent proposition,
we do not assume \eqref{eq:setting}.
Namely,
$\frak g$ is a general reductive Lie algebra
and $(\pi, V)$ denotes a finite dimensional irreducible representation
which is not necessarily faithful. 
Moreover we use the symbol $\ang{\ }{\ }$ for the symmetric bilinear form
on $\mathfrak a^*$ defined by the restriction of the Killing form
of $\mathfrak g$.

\begin{defn}[dominant minuscule weight]\label{defn:minuscule}
We say a weight $\pi_{\min}$ of $\pi$ is {\it dominant\/} and {\it minuscule\/} if
\begin{align*}
  \ang{\pi_{\min}}{\alpha}&\ge 0\quad\text{for all }\alpha\in\Sigma(\mathfrak g)^+
  \\ 
  \intertext{and}
  \ang{\pi_{\min}}{\pi_{\min}}&\le \ang{\varpi}{\varpi} \quad\text{for all }
  \varpi\in\mathcal W(\pi).
\end{align*}
If the highest weight of $\pi$ is dominant and minuscule, then $(\pi,V)$ is called
a {\it minuscule representation.}
\end{defn}

\begin{prop}\label{prop:minuscule}
Put $\Psi(\mathfrak g)=\{\alpha_1,\ldots,\alpha_r\}$
and define $\alpha^\vee = \frac{2\alpha}{\ang\alpha\alpha}$
for $\alpha\in\Sigma(\mathfrak g)$.
Let $(\pi, V)$ be a finite dimensional irreducible representation of $\mathfrak g$.
Let $\pi_{\min}$ be a dominant minuscule weight of $\pi$.

\noindent
{\rm i)}  If the highest weight of $\pi$ is in the root lattice, then $\pi_{\min} = 0$.

\noindent
{\rm ii)} $\pi_{\min}$ is uniquely determined by $\pi$.  Moreover if $(\pi',V')$ is a
finite dimensional irreducible representation of $\mathfrak g$ such that the
difference of the highest weight of $\pi'$ and that of $\pi$ is in the root lattice
of $\Sigma(\mathfrak g)$, then $\pi_{\min}=\pi'_{\min}$.

\noindent
{\rm iii)} $\varpi\in\mathcal W(\pi)$ is a dominant minuscule weight if and only if
\begin{equation}\label{eq:minuscule}
 \ang{\varpi}{\alpha^\vee}\in\{0,1\}
 \quad\text{for all }\alpha\in\Sigma(\mathfrak g)^+.
\end{equation}

\noindent
{\rm iv) }  If $\pi$ is a minuscule representation, then
$\mathcal W(\pi)=W\pi_{\min}$.

\noindent
{\rm v)} Suppose $\mathfrak g$ is simple. Let 
$\Sigma(\mathfrak g)^\vee:=\{\alpha^\vee;\,\alpha\in\Sigma(\mathfrak g)\}$
be the dual root system of\/ $\Sigma(\mathfrak g)$.
Let $\beta$ be the maximal root of $\Sigma(\mathfrak g)^\vee$ and put 
$\beta=\sum_{i=1}^rn_i\alpha_i^\vee$.
Define the fundamental 
weights $\Lambda_i$ by
$\ang{\Lambda_i}{\alpha_j^\vee}=\delta_{ij}$.
Then $\pi$ is a minuscule representation if and only if its highest weight is\/ $0$ 
or $\Lambda_i$ with $n_i=1$.
\end{prop}
\begin{proof} 
For $\alpha\in\Sigma(\mathfrak g)$ we denote by $\mathfrak g^\alpha$ the Lie
algebra generated by the root vectors corresponding to $\alpha$ and $-\alpha$.
Note that $\mathfrak g^\alpha$ is isomorphic to $\mathfrak{sl}_2$.

i) Suppose the highest weight of $\pi$ is in the root lattice.
Put $\varpi=\sum_{i=1}^r m_i(\varpi)\alpha_i$ for $\varpi\in\mathcal W(\pi)$.
Note that $m_i(\varpi)$ are integers.
Let $\varpi_0\in W(\pi)$ such that $m_i(\varpi_0)\ge 0$ and
$\sum_{i=1}^r m_i(\varpi_0)\le\sum_{i=1}^r m_i(\varpi)$ for
$\varpi\in W(\pi)$ satisfying $m_i(\varpi)\ge 0$ for $i=1,\ldots,r$.
The existence of $\varpi_0$ is clear because $m_i(\pi)\ge 0$ for
$i=1,\ldots,r$.  
Suppose $\varpi_0\ne 0$.
Since $0<\ang{\varpi_0}{\varpi_0}=\sum_{i=1}^r m_i(\varpi_0)\ang{\varpi_0}{\alpha_i}$,
there exists an index $k$ such that $\ang{\varpi_0}{\alpha_k}>0$ and $m_k(\varpi_0)>0$.
Hence $\varpi_0 - \alpha_k\in\mathcal W(\pi)$
by the representation $\pi|_{\mathfrak g^{\alpha_k}}$, which
contradicts the assumption for $\varpi_0$.
Thus $0=\varpi_0\in\mathcal W(\pi)$ and $\pi_{\min}=0$.

ii) -- iv)
Suppose the existence of $\alpha\in\Sigma(\mathfrak g)^+$
with $\ang{\pi_{\min}}{\alpha^\vee}>1$.
Then it follows from the representation $\pi|_{\mathfrak{g}^\alpha}$ that
$\pi_{\min}-\alpha\in\mathcal W(\pi)$
and
$\ang{\pi_{\min}}{\pi_{\min}}-\ang{\pi_{\min}-\alpha}{\pi_{\min}-\alpha} =
2\ang{\pi_{\min}}{\alpha} - \ang{\alpha}{\alpha} > 0$, which contradicts the
assumption of $\pi_{\min}$.
Thus we have \eqref{eq:minuscule} for $\varpi=\pi_{\min}$.

Suppose $\pi$ is an irreducible representation of $\mathfrak g$ with the
highest weight $\varpi$ satisfying \eqref{eq:minuscule}.
Suppose $\mathcal{W}(\pi)\ne W\varpi$. 
Then there exist $\mu\in W\varpi$ and 
$\mu'\in \mathcal{W}(\pi)$ such that $\mu'\notin W\varpi$ with 
$\alpha:=\mu - \mu'\in\Sigma(\mathfrak g)$.
By the $W$-invariance
we may assume $\mu=\varpi$ and therefore
$\mu'=\varpi-\alpha$ with $\alpha\in\Sigma(\mathfrak g)^+$.
Then by the representation $\pi_{\mathfrak{g}^\alpha}$
together with the condition \eqref{eq:minuscule}
we have $\ang{\varpi}{\alpha^\vee}=1$ and  
$\mu'=w_\alpha \varpi$,
which is a contradiction.
Thus we have iv).

Let $\varpi$ and $\varpi'$ be the elements of $\mathfrak a^*$
satisfying the condition \eqref{eq:minuscule}.
Then $\varpi'':=\varpi-\varpi'$ satisfies
$\ang{\varpi''}{\alpha^\vee}\in \{-1,0,1\}$ for $\alpha\in\Sigma(\mathfrak g)$.
Suppose that $\varpi''$ is in the root lattice.  
Let $\varpi_0\in W\varpi''$ such that $\ang{\varpi_0}{\alpha}\ge 0$
for $\alpha\in\Sigma(\mathfrak g)^+$.
Since $\varpi_0$ also satisfies \eqref{eq:minuscule},
the finite dimensional irreducible representation $\pi_0$ with
the highest weight $\varpi_0$ is minuscule by the argument above.
Since $\varpi_0$ is in the root
lattice, $\varpi_0=0$ by i) and hence $\varpi=\varpi'$.
Thus we obtain
ii) and iii).

v) 
Let $\alpha\in\Sigma(\mathfrak g)^+$.
If we denote $\alpha^\vee = \sum_{i=1}^r n_i(\alpha)\alpha_i^\vee$,
then $n_i(\alpha)\le n_i$ for $i=1,\ldots,r$.
Hence the claim is clear.
\end{proof}
\begin{rem}
Equivalent contents of Proposition~\ref{prop:minuscule} are found in
exercises of \cite{Bo1},~Ch.~VI.%
\end{rem}
Restore the previous setting \eqref{eq:setting} on $\mathfrak g$ and $(\pi,V)$.

\begin{prop}\label{prop:free}
{\rm i) }
Let $V_{\varpi}$ denote the weight space of $V$ with weight $\varpi\in\mathcal W(\pi)$.  
Define the projection map 
$\bar p_\Theta:\mathcal W(\pi)\to\mathcal W(\pi)|_{\mathfrak a_\Theta}$ by
$\bar p_\Theta(\varpi)=\varpi|_{\mathfrak a_\Theta}$
and put $V(\Lambda)=\sum_{\varpi\in\bar p^{-1}_\Theta(\Lambda)}V_\varpi$
for $\Lambda\in\mathcal W(\pi)|_{\mathfrak a_\Theta}$.
Then
\begin{equation}\label{eq:rest1}
   V = \bigoplus_{\Lambda\in\mathcal W(\pi)|_{\mathfrak a_\Theta}}V(\Lambda)
\end{equation}
is a direct sum decomposition of the $\mathfrak g_\Theta$-module $V$.

Let $V(\Lambda)=V(\Lambda)_1\oplus\cdots\oplus V(\Lambda)_{k_\Lambda}$ be a 
decomposition into irreducible $\mathfrak g_\Theta$-modules. 
We denote by $\varpi_\Lambda$ the dominant minuscule weight of 
$(\pi|_{\mathfrak g_\Theta}, V(\Lambda)_1)$.
Then
\begin{equation}\label{eq:rest2}
  V_{\varpi_\Lambda} = 
  \bigoplus_{i=1}^{k_\Lambda} V_{\varpi_\Lambda}\cap V(\Lambda)_i
  \qquad \text{with } \dim V_{\varpi_\Lambda}\cap V(\Lambda)_i>0.
\end{equation}
In particular, $V(\Lambda)$ is an irreducible 
$\mathfrak g_\Theta$-module if $\dim V_{\varpi_\Lambda}=1$.

\noindent
{\rm ii)} 
Put $\Psi(\mathfrak g)=\{\alpha_1,\cdots,\alpha_r\}$ and %
put $\Psi(\mathfrak g)\setminus\Theta
=\{\alpha_{i_1},\ldots,\alpha_{i_k}\}$, 
define the map
\[
  \begin{matrix}
   p_\Theta:&\Sigma(\mathfrak g)&\to&\mathbb Z^k\\
            &\alpha=\sum m_i\alpha_i&\mapsto&(m_{i_1},\ldots,m_{i_k})
  \end{matrix}
\]
and put 
\begin{align*}
  L_\Theta &= \{0\}\cup\{p_\Theta(\alpha);\,\alpha\in\Sigma(\mathfrak g)\},\\
  V(\mathbf m) &=
  \begin{cases}
    \displaystyle\sum_{\alpha\in p_\Theta^{-1}(\mathbf m)}\mathbb CX_\alpha
    &\text{if \ }\mathbf m\ne 0,\\
    \displaystyle\mathfrak a + 
     \sum_{\alpha\in p_\Theta^{-1}(\mathbf m)}\mathbb CX_\alpha
    &\text{if \ }\mathbf m= 0
  \end{cases}
\end{align*}
for\/ $\mathbf m\in L_\Theta$.
Then 
\begin{equation}\label{eq:adj}
  \mathfrak g=\bigoplus_{\mathbf m\in L_\Theta}V(\mathbf m)
\end{equation}
is a decomposition of the $\mathfrak g_\Theta$-module $\mathfrak g$.
If\/ $\mathbf m\ne0$, then\/ $V(\mathbf m)$ is an irreducible\/
$\mathfrak g_\Theta$-module.
On the other hand,\/ $V(0)=\mathfrak g_\Theta$ is isomorphic
to the adjoint representation of\/ $\mathfrak g_\Theta %
=\mathfrak{a}_\Theta \oplus \mathfrak{m}_\Theta$.
Let $\Theta=\Theta_1 \sqcup \Theta_2 \sqcup \cdots \sqcup \Theta_\ell$
be the division of $\Theta$ into the connected parts of vertexes
in the Dynkin diagram of $\Psi(\mathfrak g)$.
Then $\mathfrak{m}_\Theta=\mathfrak{m}_{\Theta_1}\oplus\mathfrak{m}_{\Theta_2}\oplus
\cdots\oplus\mathfrak{m}_{\Theta_\ell}$
gives a decomposition into irreducible $\mathfrak{g}_\Theta$-modules.

\noindent
{\rm iii)}
Suppose that the representation $(\pi, V)$ is minuscule.
Put $W^\pi=\{ w\in W;\, w\pi=\pi\}$.
Here we identify $\pi$ with its highest weight.
Let $\{w_1,\ldots,w_k\}$ be a representative system
of $W^\pi \backslash W / W_\Theta$ such that
$w_i \in W(\Theta)$.
Then with the notation in\/ {\rm i)}
\begin{equation}\label{eq:mindeco}
V = \bigoplus_{i=1}^k V(w_i^{-1}\pi|_{\mathfrak{a}_\Theta})
\end{equation}
gives a decomposition into irreducible $\mathfrak{g}_\Theta$-modules.
Moreover the $\mathfrak{g}_\Theta$-submodule $V(w_i^{-1}\pi|_{\mathfrak{a}_\Theta})$
has highest weight $w_i^{-1}\pi$.
\end{prop}

\begin{proof} i)
Since $\alpha|_{\mathfrak a_\Theta} = 0$ for $\alpha\in\Theta$, 
\eqref{eq:rest1} is a decomposition into $\mathfrak g_\Theta$-modules.
Then Proposition~\ref{prop:minuscule}~ii) implies that $\varpi_\Lambda$ is the minuscule 
weight for any $(\pi|_{\mathfrak g_\Theta}, V(\Lambda)_i)$ and therefore the other statements
in i) are clear.

ii)
Note that 
$\alpha_{i_k}|_{\mathfrak a_\Theta},\ldots,\alpha_{i_1}|_{\mathfrak a_\Theta}$
are linearly independent and $\mathfrak g_\Theta = V(0)$.
Then the statements in ii) follows from i).

iii)
From i)
each $V(w_i^{-1}\pi|_{\mathfrak{a}_\Theta})$ is an irreducible $\mathfrak{g}_\Theta$-module
and
\[
V(w_i^{-1}\pi|_{\mathfrak{a}_\Theta}) \supset 
\sum\{V_{w^{-1}\pi};\,w\in W^\pi w_i W_\Theta \}.
\]
Since $w_i\in W(\Theta)$ we have $w_i^{-1}\pi + \alpha \notin \mathcal{W}(\pi)$
for $\alpha\in\Sigma(\mathfrak{g}_\Theta)^+$.
It shows the highest weight of $V(w_i^{-1}\pi|_{\mathfrak{a}_\Theta})$ is $w_i^{-1}\pi$. 
Since $w_i^{-1}\pi\ne w_j^{-1}\pi$ if $i\ne j$
we have \eqref{eq:mindeco}.
\end{proof}

We give the minimal polynomials for some representations in the following 
proposition as a corollary of Lemma~\ref{lem:casimir}~v) and
Proposition~\ref{prop:free}.
\begin{prop}\label{prop:minfree}%
Retain the notation in\/ {\rm Theorem~\ref{thm:min}} and\/ 
{\rm Proposition~\ref{prop:free}}.

\noindent
{\rm i) (multiplicity free representation)}
Suppose $\dim V_\varpi=1$ for any $\varpi\in\mathcal W(\pi)$.
Let $\bar\Lambda$ be the lowest weight of $(\pi|_{\mathfrak g_\Theta}, V(\Lambda))$
for $\Lambda\in\mathcal W(\pi)|_{\mathfrak a_\Theta}$.
Then
\begin{equation}
\begin{aligned}
 q_{\pi,\Theta}(x;\lambda)
 &= \prod_{\Lambda\in\mathcal W(\pi)|_{\mathfrak a_\Theta}}
  \left(
   x - \ang{\lambda_\Theta}{\bar\Lambda} 
    - \frac12\ang{\bar\pi-\bar\Lambda}{\bar\pi+\bar\Lambda-2\rho}
  \right)\\
 &=  \prod_{\Lambda\in\mathcal W(\pi)|_{\mathfrak a_\Theta}}
  \left(
   x - \ang{\lambda_\Theta+\rho}{\bar\Lambda} + \ang{\bar\pi}{\rho}
    - \frac{\ang{\bar\pi}{\bar\pi}-\ang{\bar\Lambda}{\bar\Lambda}}2
  \right).
\end{aligned}
\end{equation}
\noindent
{\rm ii) (adjoint representation)}
Suppose $\mathfrak g$ is simple and $\Theta\ne\emptyset$.
Let $\Theta=\Theta_1\sqcup\cdots\sqcup\Theta_\ell$ be the division
in\/ {\rm Proposition~~\ref{prop:free}~ii)}.
Let $\alpha_{\max}^i$ denote the maximal root of the simple Lie algebra 
$\mathfrak{m}_{\Theta_i}$ for $i=1,\ldots,\ell$.
Put
\[\Omega_\Theta=
\{ B(\alpha_{\max}^1,\alpha_{\max}^1+2\rho(\Theta_1)), \ldots, 
B(\alpha_{\max}^\ell,\alpha_{\max}^\ell+2\rho(\Theta_\ell)) \}.\]
Let $\alpha_{\mathbf m}$ be the smallest root in $p_\Theta^{-1}(\mathbf m)$
for $\mathbf m\in L_\Theta\setminus\{0\}$ under the order in\/ {\rm Definition~\ref{def:shift}}.
Then for the adjoint representation of $\mathfrak g$,
\begin{multline}\label{eq:admin}
  q_{\alpha_{\max},\Theta}(x;\lambda) = 
  \left(x-\frac12\right)
  \prod_{C\in\Omega_\Theta}\left(x-\frac{1-C}2\right)\\
  \cdot \prod_{\mathbf m\in L_\Theta\setminus\{0\}}\left(x - B(\lambda_\Theta+\rho,\alpha_{\mathbf m})
     -\frac{1 - B(\alpha_{\mathbf m},\alpha_{\mathbf m})}2\right).
\end{multline}

\noindent
{\rm iii) (minuscule representation)}
Suppose $(\pi, V)$ is minuscule.
Then with $w_1, \ldots, w_k$ in\/ 
{\rm Proposition~\ref{prop:free}~iii)},
\begin{equation}
q_{\pi,\Theta}(x;\lambda)=
\prod_{i=1}^k\Bigl(x-\ang{w_i \bigl(\lambda_\Theta+\rho_\Theta-\rho(\Theta)\bigr)+\rho}{\pi}\Bigr).
\end{equation}
\end{prop}
\begin{proof}
It is easy to get i) and ii).

iii)
Let $\bar{w}_{\Theta}$ denote the longest element in $W_\Theta$.
Then the $\mathfrak{g}_\Theta$-module $V(w_i^{-1}\pi|_{\mathfrak{a}_\Theta})$
has lowest weight $\bar{w}_\Theta w_i^{-1}\pi$.
The claim follows from the next calculation:
\begin{align*}
\ang{\lambda_\Theta}{\bar{w}_\Theta w_i^{-1}\pi}
&+\frac12\ang{\bar{\pi}-\bar{w}_\Theta w_i^{-1}\pi}{\bar{\pi}+\bar{w}_\Theta w_i^{-1}\pi-2\rho}\\
&=\ang{\lambda_\Theta+\rho}{\bar{w}_\Theta w_i^{-1}\pi}
+\ang{\rho}{\pi}\\
&=\ang{w_i\bar{w}_\Theta(\lambda_\Theta+\rho)+\rho}{\pi}
=\ang{w_i\bigl(\lambda_\Theta+\rho_\Theta-\rho(\Theta)\bigr)+\rho}{\pi}.
\qedhere\end{align*}
\end{proof}

\section{Two-sided ideals}\label{sec:ideal}
Our main concern in this paper is the following two-sided ideal.

\begin{defn}[gap]
Let $\lambda\in\mathfrak{a}_\Theta^*$.
If a two-sided ideal $I_\Theta(\lambda)$ of $U(\mathfrak g)$
satisfies
\begin{equation}\label{eq:gapdef}
 J_\Theta(\lambda) = I_\Theta(\lambda) + J(\lambda_\Theta),
\end{equation}
then we say that $I_\Theta(\lambda)$ {\it describes the gap\/} between
the generalized Verma module $M_\Theta(\lambda)$ and 
the Verma module $M(\lambda_\Theta)$.
\end{defn}

It is clear that there exists a two-sided ideal $I_\Theta(\lambda)$ satisfying
\eqref{eq:gapdef} 
if and only if
\begin{equation}\label{eq:gapann}
J_\Theta(\lambda)=\Ann\bigl(M_\Theta(\lambda)\bigr) + J(\lambda_\Theta).
\end{equation}
This condition depends on $\lambda$
but such an ideal exists and is essentially unique for a generic
$\lambda$ (cf.~Proposition~\ref{prop:HCzero}, Theorem~\ref{thm:gapexist}, Remark~\ref{rm:existence}).
The main purpose in this paper is to construct a good generator system of the ideal
from a minimal polynomial.

\begin{defn}[two-sided ideal]
Using the global minimal polynomial defined in the last section, 
we define a two-sided ideal of $U(\mathfrak g)$:
\begin{equation}
  I_{\pi,\Theta}(\lambda)
  = \sum_{i,j}U(\mathfrak g)
    q_{\pi,\Theta}(F_\pi;\lambda)_{ij}
   + \sum_{\Delta\in Z(\mathfrak g)}
    U(\mathfrak g)\bigl(\Delta - \Delta_{\mathfrak a}(\lambda_\Theta)\bigr).
\end{equation}
\end{defn}

From Theorem~\ref{thm:min} and Remark~\ref{rem:inf_char}
this ideal satisfies
\begin{equation}
  I_{\pi,\Theta}(\lambda)\subset J_\Theta(\lambda).
\end{equation}
In this section we will examine the condition so that
\begin{equation}\label{eq:gap}
  J_\Theta(\lambda) = I_{\pi,\Theta}(\lambda) + J(\lambda_\Theta).
\end{equation}

\begin{prop}[invariant differential operators]\label{prop:geninv}
For $\Delta\in Z(\mathfrak g)$ and a non-negative integer $k$ we denote by 
$\Delta_{\mathfrak a}^{(k)}$ the homogeneous part of $\Delta_{\mathfrak a}$ with degree $k$
and put
\begin{equation}
 T_{\pi}^{(k)}= \sum_{\varpi\in\mathcal W(\pi)}m_\pi(\varpi)\varpi^k.
\end{equation}
Here $m_\pi(\varpi)$ is the multiplicity of the weight $\varpi$ of $\pi$ and 
we use the identification 
$\varpi\in\mathfrak a^*\simeq\mathfrak a\subset U(\mathfrak a)$.
Let $\{\Delta_1,\ldots,\Delta_r\}$ be a system of generators of $Z(\mathfrak g)$
as an algebra over $\mathbb C$ and let $d_i$ be the degree of 
$(\Delta_i)_{\mathfrak a}$ for $i=1,\ldots,r$. 
We assume that $(\Delta_1)_{\mathfrak a}^{(d_1)},\ldots,
(\Delta_r)_{\mathfrak a}^{(d_r)}$ are algebraically independent.
Suppose a subset $A$ of $\{1,\ldots,r\}$ satisfies
\begin{equation}
 \begin{cases}
  d_k\ge\deg_x q_{\pi,\Theta}(x,\lambda)\quad\text{if }k\in\{1,\ldots,r\}\setminus A,\\
  \mathbb C[(\Delta_1)_{\mathfrak a}^{(d_1)},\ldots,(\Delta_r)_{\mathfrak a}^{(d_r)}]
  = \mathbb C[(\Delta_i)_{\mathfrak a}^{(d_i)}, T_{\pi}^{(d_k)};
  i\in A,\ k\in \{1,\ldots,r\}\setminus A].
 \end{cases}
\end{equation}
Then
\begin{equation}
  I_{\pi,\Theta}(\lambda) = \sum_{i,j}U(\mathfrak g)
    q_{\pi,\Theta}(F_\pi;\lambda)_{ij}
   + \sum_{i\in A}
    U(\mathfrak g)\bigl( \Delta_i - (\Delta_i)_{\mathfrak a}(\lambda_\Theta)\bigr).
\end{equation}
\end{prop}
\begin{proof}
Note that $\sum_{i,j}U(\mathfrak g)q_{\pi,\Theta}(F_\pi;\lambda)_{ij}
\ni\trace\bigl(F_\pi^\nu q_{\pi,\Theta}(F_\pi;\lambda)\bigr)$ if $\nu\ge 0$.
On the other hand, 
since $\trace\bigl(F_\pi^{\ell_k}
q_{\pi,\Theta}\bigl(F_\pi;\lambda)\bigr)^{(d_k)}_{\mathfrak a} =
T_{\pi}^{(d_k)}$ by Lemma~\ref{lem:zen} with $\Theta=\emptyset$ if the integer
$\ell_k=d_k -\deg_x\bigl(q_{\pi,\Theta}(F_\pi;\lambda)\bigr)$ is non-negative, 
the assumption implies that for $k\not\in A$,
$\Delta_k$ may be replaced by 
$\trace\bigl(F_\pi^{\ell_k}
q_{\pi,\Theta}\bigl(F_\pi;\lambda)\bigr)$,
which implies the proposition.
\end{proof}

\begin{lem}\label{lem:gap}
Let $\mathbf V$ be an $\ad(\mathfrak g)$-stable subspace of 
$U(\mathfrak g)$ and let $\mathbf V=\bigoplus_{\varpi}\mathbf V_\varpi$ be
the decomposition of $\mathbf V$ into the weight spaces $\mathbf V_\varpi$ with
weight $\varpi\in\mathfrak a^*$.
Suppose $D_{\mathfrak a}(\lambda_\Theta)=0$ for $D\in \mathbf V_0$.
Then the following three conditions are equivalent.

\noindent
{\rm i)}\quad
 $J_\Theta(\lambda) \subset U(\mathfrak g)\mathbf V + J(\lambda_\Theta)$.

\noindent
{\rm ii)}
 For any $\alpha\in\Theta$ there exists $D\in \mathbf V_{-\alpha}$ such that
   $D - X_{-\alpha}\in J(\lambda_\Theta)$.

\noindent
{\rm iii)}
 For any $\alpha\in\Theta$ there exists $D\in\mathbf V_0$  such that
  $D_{\mathfrak a}(\lambda_\Theta-\alpha)\ne0$.
\end{lem}
\begin{proof}
Let $U(\mathfrak g)=\bigoplus_{\varpi}U(\mathfrak g)_\varpi$ be
the decomposition of $U(\mathfrak g)$ into the weight spaces $U(\mathfrak g)_\varpi$ with
weight $\varpi\in\mathfrak a^*$.
Let $\mu\in\mathfrak a^*$.
Since $U(\mathfrak g) %
= U(\bar{\mathfrak n})\oplus J(\mu)$, 
to $D\in U(\mathfrak g)$,
there corresponds %
a unique $D^\mu\in U(\bar{\mathfrak n})$
such that %
$D-D^\mu\in J(\mu)$.
Here we note that 
$D\in U(\mathfrak g)_\varpi$ implies
$D^\mu\in U(\bar{\mathfrak n})_\varpi$
and that $D^\mu=D_{\mathfrak a}(\mu)\in\mathbb C$
whenever $D\in U(\mathfrak g)_0$.

Put ${\mathbf V}^\mu=\{D^\mu;\,D\in \mathbf V\}$.
Since $\ad(X)\mathbf V\subset\mathbf V$ for $X\in\mathfrak b$, 
we have $PD\in \mathbf V +J(\mu)$ and therefore 
$(PD)^\mu\in{\mathbf V}^\mu$ for every $P\in U(\mathfrak b)$ and $D\in\mathbf V$.
Owing to $U(\mathfrak g)=U(\bar{\mathfrak n})\otimes
U(\mathfrak b)$, we have
\begin{equation}\label{eq:ideal}
  \{D^\mu; D\in U(\mathfrak g)\mathbf V\} = 
  U(\bar{\mathfrak n}){\mathbf V}^\mu.
\end{equation}
Note that 
\begin{equation}\label{eq:Vimage}
{\mathbf V}^\mu=\bigoplus\{(\mathbf V_\varpi)^\mu; 
\varpi=-\sum_{\gamma\in\Psi(\mathfrak g)}n_\gamma \gamma 
\text{ for some non-negative integers } n_\gamma \}.
\end{equation}

Suppose i).  Let $\alpha\in\Theta$.  
Since $X_{-\alpha}\in J_\Theta(\lambda)\setminus J(\lambda_\Theta)$, there exists 
$D\in U(\mathfrak g)\mathbf V$ with $D^{\lambda_\Theta}=X_{-\alpha}$.
On the other hand,
we can deduce $\left( U(\bar{\mathfrak n})\mathbf V^{\lambda_\Theta} \right)_{-\alpha}
=(\mathbf V_{-\alpha})^{\lambda_\Theta}$ from \eqref{eq:Vimage}
because the assumption of the lemma
assures $(\mathbf V_0)^{\lambda_\Theta}=0$.
Hence from \eqref{eq:ideal} 
we may assume $D\in\mathbf V_{-\alpha}$.
Thus we have ii).

It is clear that ii) implies i) because 
$J_\Theta(\lambda) = J(\lambda_\Theta) +
\sum_{\alpha\in\Theta}U(\mathfrak g)X_{-\alpha}$.

Let $\alpha\in\Theta$.
Since $\ad(H)X_{-\alpha} = -\alpha(H)X_{-\alpha}$ for $H\in\mathfrak a$,
we have %
$H_1\cdots H_k X_{-\alpha} = 
X_{-\alpha}(H_1-\alpha(H_1))\cdots (H_k-\alpha(H_k))$ for 
$H_1,\ldots,H_k\in\mathfrak a$.
We also have $X_\gamma X_{-\alpha}\in J(\lambda_\Theta)$
for $\gamma\in\Sigma(\mathfrak g)^+$
because $\lambda_\Theta([X_\alpha, X_{-\alpha}])=0$
and $[X_\gamma, X_{-\alpha}] \in \mathfrak n$ if $\gamma \ne \alpha$.
Hence for any $D\in U(\mathfrak g)_0$,
\begin{equation}\label{eq:aton}
 (\ad(X_{-\alpha})D)^{\lambda_\Theta} = 
[X_{-\alpha}, D_{\mathfrak a}]^{\lambda_\Theta} = 
  \bigl(D_{\mathfrak a}(\lambda_\Theta) - 
  D_{\mathfrak a}(\lambda_\Theta-\alpha)\bigr)X_{-\alpha}.
\end{equation}
Now it is clear that iii) implies ii).

Conversely suppose ii).
Let $\alpha\in\Theta$.
Since $\mathbf V_{-\alpha} = \ad(X_{-\alpha})\mathbf V_0$,
there exists $D\in\mathbf V_0$ with
$(\ad(X_{-\alpha})D)^{\lambda_\Theta}=X_{-\alpha}$
and we have iii) from \eqref{eq:aton}.
\end{proof}

\begin{rem}\label{rem:gapalpha}
In the above lemma $\lambda_\Theta-\alpha = w_\alpha.\lambda_\Theta$ for 
$\alpha\in\Theta$ because $\ang{\lambda_\Theta}{\alpha}=0.$
\end{rem}

By the {\it Duflo theorem\/} (\cite{Du}),
$\Ann\bigl(M(\mu)\bigr)=\sum_{\Delta\in Z(\mathfrak g)}U(\mathfrak g)\bigl(\Delta-\Delta_{\mathfrak a}(\mu)\bigr)$
for any $\mu\in\mathfrak a^*$. 
Then, by the following theorem, 
each $\Ann\bigl(M(\mu)\bigr)$ has the same $\ad(\mathfrak g)$-module structure.

\begin{thm}[the Kostant theorem \cite{Ko1}]\label{thm:Kostant}
There exists an $\ad(\mathfrak g)$-submodule $\mathcal{H}$ of $U(\mathfrak g)$
such that %
$U(\mathfrak g)$
is naturally isomorphic to $Z(\mathfrak g)\otimes\mathcal{H}$ by the multiplication.
For any finite dimensional $\mathfrak g$-module $\mathbf{V}$,
$\dim\Hom_{\mathfrak g}\left(\mathbf{V},\mathcal H\right)=\dim \mathbf{V}_0$.
\end{thm}

Similarly on the annihilators of generalized Verma modules we have

\begin{prop}\label{prop:ann_str}
Suppose $\lambda_\Theta+\rho$
is dominant.
Then for any finite dimensional $\mathfrak g$-module $\mathbf{V}$ and $\mathcal{H}$
in\/ {\rm Theorem~\ref{thm:Kostant}},
\begin{multline*}
\dim \Hom_{\mathfrak g}\left(\mathbf{V}, \Ann\bigl(M_\Theta(\lambda)\bigr)/\Ann\bigl(M(\lambda_\Theta)\bigr) \right)\\
=\dim \Hom_{\mathfrak g}\left(\mathbf{V}, \mathcal{H}\cap\Ann\bigl(M_\Theta(\lambda)\bigr)\right)
=\dim \mathbf{V}_0- \dim \mathbf{V}^{\mathfrak g_\Theta}
\end{multline*}
where $\mathbf{V}^{\mathfrak g_\Theta}
=\{ v\in\mathbf V;\, Xv=0\ (\forall X\in \mathfrak g_\Theta)\}$.
\end{prop}

Before proving the proposition,
we accumulate some necessary facts from %
\cite{BGG}, \cite{BG} and \cite{Jo0}.

\begin{defn}[category $\mathcal O$ \cite{BGG}]
Let $\mathcal O$ be the abelian category consisting of
the $\mathfrak g$-modules which are finitely generated, $\mathfrak{a}^*$-diagonalizable
and $U(\mathfrak n)$-finite.
All subquotients of Verma modules are objects of $\mathcal O$.
For $\mu\in\mathfrak{a}^*$ we denote by $L(\mu)$
the unique irreducible quotient of the Verma module $M(\mu)$.
There exists a unique indecomposable projective object $P(\mu)\in\mathcal O$ such that 
$\Hom_{\mathfrak g}(P(\mu),L(\mu))\ne0$.
\end{defn}

\begin{prop}[\cite{BGG}, \cite{BG}]\label{prop:Ofacts}
{\rm i)}
If $\mu+\rho$ is dominant,
then $P(\mu)=M(\mu)$
and 
\[
\dim\Hom_{\mathfrak g}(M(\mu),M(\mu'))
=\begin{cases}
1 & \text{if }\mu'=\mu, \\
0 & \text{if }\mu'\ne\mu.
\end{cases}\]

\noindent {\rm ii)}
For any $\mu, \mu'\in\mathfrak a^*$
\[
\dim\Hom_{\mathfrak g}(P(\mu),L(\mu'))
=\begin{cases}
1 & \text{if }\mu'=\mu, \\
0 & \text{if }\mu'\ne\mu.
\end{cases}\]

\noindent {\rm iii)}
For any finite dimensional $\mathfrak g$-module $\mathbf V$
and $\mu\in\mathfrak a^*$,
$\mathbf{V}\otimes P(\mu)$ is a projective object in $\mathcal O$.
\end{prop}

\begin{prop}[\cite{BG}, \cite{Jo0}]\label{prop:2toleft}
Suppose $\mu\in\mathfrak a^*$ and $\mu+\rho$ is dominant.
Then the map
\begin{equation}\label{eq:2toleft}
\{I\subset U(\mathfrak g);\,\text{two-sided ideal, }I \supset \Ann\bigl(M(\mu)\bigr)\}
\rightarrow
\{M\subset M(\mu);\,\text{submodule}\}
\end{equation}
defined by $I\mapsto IM(\mu)$ is injective
and hence $\Ann\bigl(M(\mu)/IM(\mu)\bigr)=I$ for any two-sided ideal $I$
with $I\supset \Ann\bigl(M(\mu)\bigr)$.
The image of the map~\eqref{eq:2toleft} consists of the submodules
which are isomorphic to quotients of direct sums of $P(\mu')$ with
\begin{equation}\label{eq:proj_rep}
2\frac{\ang{\mu'+\rho}{\beta}}{\ang{\beta}{\beta}}\in\{0,-1,-2,\ldots\}
\text{ for any }\beta\in\Sigma(\mathfrak g)^+\text{ such that }\ang{\mu+\rho}{\beta}=0.
\end{equation}
\end{prop}

\begin{proof}[Proof of\/ {\rm Proposition \ref{prop:ann_str}}]
We first show the map
\begin{equation}\label{eq:VMM}
\Hom_{\mathfrak g}\left(\mathbf{V}, \mathcal{H}\right)\ni\varphi
\mapsto\Phi\in
\Hom_{\mathfrak g}\left(\mathbf{V}\otimes M(\lambda_\Theta), M(\lambda_\Theta)\right)
\end{equation}
defined by $\Phi(v\otimes u)=\varphi(v)u$
is a linear isomorphism.
Since $U(\mathfrak g)=\mathcal{H}\oplus \Ann\bigl(M(\lambda_\Theta)\bigr)$
the map is injective.
To show the surjectivity we calculate the dimensions of both spaces.
By Theorem~\ref{thm:Kostant}
$\dim \Hom_{\mathfrak g}\left(\mathbf{V}, \mathcal{H}\right)
=\dim \mathbf{V}_0$.
On the other hand, note that 
\[\Hom_{\mathfrak g}\left(\mathbf{V}\otimes M(\lambda_\Theta), M(\lambda_\Theta)\right)
 \simeq
 \Hom_{\mathfrak g}\left(M(\lambda_\Theta), M(\lambda_\Theta) \otimes \mathbf{V}^*\right)
\]
and there exist a sequence $\{\mu_1,\ldots,\mu_\ell\}\subset\mathfrak a^*$
and a $\mathfrak g$-stable filtration
\[
\{0\}=M_0\subsetneq M_1 \subsetneq \cdots \subsetneq M_\ell = M(\lambda_\Theta) \otimes \mathbf{V}^*\\
\]
such that $M_i/M_{i-1} \simeq M(\mu_i)$ for $i=1,\ldots,\ell$.
Here the number of appearances of $\lambda_\Theta$ in the sequence $\{\mu_1,\ldots,\mu_\ell\}$
equals $\dim \mathbf{V}^*_0=\dim \mathbf{V}_0$ (cf.~the proof of Theorem~\ref{thm:minimality}).
Since $\lambda_\Theta+\rho$ is dominant, 
it follows from Proposition~\ref{prop:Ofacts}~i) that
$\dim \Hom_{\mathfrak g}\left(M(\lambda_\Theta), M(\lambda_\Theta) \otimes \mathbf{V}^*\right)
=\dim \mathbf{V}_0$.
Thus \eqref{eq:VMM} is isomorphism.

Secondly, consider the exact sequence
\[
0\rightarrow J_\Theta(\lambda)/J(\lambda_\Theta)
\rightarrow M(\lambda_\Theta) \rightarrow M_\Theta(\lambda)
\rightarrow 0.
\]
It is clear that under the isomorphism~\eqref{eq:VMM}
the subspace
\[
 \Hom_{\mathfrak g}\left(\mathbf{V}, \mathcal{H}\cap\Ann\bigl(M_\Theta(\lambda)\bigr)\right)
 \subset
 \Hom_{\mathfrak g}\left(\mathbf{V}, \mathcal{H}\right) 
\]
corresponds to the subspace
\[
 \Hom_{\mathfrak g}\left(\mathbf{V}\otimes M(\lambda_\Theta), J_\Theta(\lambda)/J(\lambda_\Theta)\right)
 \subset
 \Hom_{\mathfrak g}\left(\mathbf{V}\otimes M(\lambda_\Theta), M(\lambda_\Theta)\right).
\]
Let us calculate the dimension of the latter space.
By Proposition~\ref{prop:Ofacts}~i) and iii),
$\mathbf{V}\otimes M(\lambda_\Theta)$ is projective and therefore
\begin{multline*}
\dim \Hom_{\mathfrak g}\left(\mathbf{V}\otimes M(\lambda_\Theta), J_\Theta(\lambda)/J(\lambda_\Theta)\right)\\
=\dim \Hom_{\mathfrak g}\left(\mathbf{V}\otimes M(\lambda_\Theta), M(\lambda_\Theta)\right)
-\dim \Hom_{\mathfrak g}\left(\mathbf{V}\otimes M(\lambda_\Theta), M_\Theta(\lambda)\right).
\end{multline*}
Here we know
\[\Hom_{\mathfrak g}\left(\mathbf{V}\otimes M(\lambda_\Theta), M_\Theta(\lambda)\right)
 \simeq
 \Hom_{\mathfrak g}\left(M(\lambda_\Theta), M_\Theta(\lambda) \otimes \mathbf{V}^*\right)
\]
and there exist a sequence $\{\mu_1,\ldots,\mu_{\ell'}\}\subset\mathfrak a^*$
and a $\mathfrak g$-stable filtration
\[
\{0\}=M_0\subsetneq M_1 \subsetneq \cdots \subsetneq M_{\ell'} = M_\Theta(\lambda) \otimes \mathbf{V}^*\\
\]
such that $M_i/M_{i-1} \simeq M_{(\Theta,\mu_i)}$ for $i=1,\ldots,\ell'$.
The number of appearances of $\lambda_\Theta$ in the sequence $\{\mu_1,\ldots,\mu_{\ell'}\}$
equals $\dim (\mathbf{V}^*)^{\mathfrak g_\Theta}=\dim \mathbf{V}^{\mathfrak g_\Theta}$ (cf.~the proof of Theorem~\ref{thm:minimality}).
Since the generalized Verma module $M_{(\Theta,\mu_i)}$ is a quotient of $M(\mu_i)$,
Proposition~\ref{prop:Ofacts}~i) implies
$\dim \Hom_{\mathfrak g}\left(M(\lambda_\Theta), M_\Theta(\lambda) \otimes \mathbf{V}^*\right)
 =\dim \mathbf{V}^{\mathfrak g_\Theta}$.
Thus the proposition is proved.
\end{proof} 

\begin{prop}[Harish-Chandra homomorphism]\label{prop:HCzero}
Let $I$ be a two-sided ideal of $U(\mathfrak g)$.
Put $\mathcal{V}(I)=\{\mu\in\mathfrak a;\,D_\mathfrak a(\mu)=0\ (\forall D\in I)\}$.

\noindent
{\rm i)} Fix $\alpha\in\Psi(\mathfrak g)$.
If $\mu\in \mathcal{V}(I)$ and
\begin{equation}\label{eq:ref}
  2\frac{\ang{\mu+\rho}{\alpha}}{\ang{\alpha}{\alpha}}\notin\{1,2,3,\ldots\},
\end{equation}
then $w_\alpha.\mu\in \mathcal{V}(I)$.

\noindent
{\rm ii)} Suppose $\lambda\in\mathfrak a_\Theta^*$ and
\begin{equation}\label{eq:gapmu}
  J_\Theta(\lambda) = I + J(\lambda_\Theta).
\end{equation}
Then $w.\lambda_\Theta\notin \mathcal{V}(I)$ for $w\in W_\Theta\setminus\{e\}$.

\noindent
{\rm iii)} In addition to the assumption of\/ {\rm ii)}, suppose $\lambda_\Theta+\rho$ is dominant
and
\begin{equation}\label{eq:bigideal}
I \supset \Ann\bigl(M(\lambda_\Theta)\bigr).
\end{equation}
Then $I=\Ann\bigl(M_\Theta(\lambda)\bigr)$ and
\begin{equation}\label{eq:Annroots}
   \mathcal V(I) = 
   \{w.\lambda_\Theta;\, w\in W(\Theta)\}.
\end{equation}
\end{prop}
\begin{proof} i)
Note that $\mu\in \mathcal{V}(I)$ if and only if $I\subset\Ann\bigl(L(\mu)\bigr)$.
It is known by \cite{Jo} that 
$\Ann\bigl(L(\mu)\bigr) \subset \Ann\bigl(L(w_\alpha.\mu)\bigr)$
if \eqref{eq:ref} holds, which implies i).

ii)
Since $I\subset\Ann\bigl(M_\Theta(\lambda)\bigr)\subset\Ann\bigl(L(\lambda_\Theta)\bigr)$
we have $\lambda_\Theta\in \mathcal{V}(I)$.
Put $W'=\{w\in W_\Theta\setminus\{e\};\, w.\lambda_\Theta\in \mathcal{V}(I)\}$.
Then, by Lemma~\ref{lem:gap} with $\mathbf V=I$,
$w_\alpha\notin W'$ for any $\alpha\in\Theta$.
Suppose $W'\ne\emptyset$.
Let $w'$ be an element of $W'$ with the minimal length.
Then there exists $\alpha\in\Theta$ such that
the length of $w''=w_\alpha w'$ is smaller than that of $w'$.
Then $w''\ne e$ and
\[
 2\frac{\ang{w'.\lambda_\Theta+\rho}{\alpha}}{\ang{\alpha}{\alpha}} = 
 2\frac{\ang{w'\rho}{\alpha}}{\ang{\alpha}{\alpha}} < 0.
\]
Hence by i), we have $w''.\mu\in \mathcal{V}(I)$, which is a contradiction.

iii)
It immediately follows from Proposition~\ref{prop:2toleft}
that $I=\Ann\bigl(M_\Theta(\lambda)\bigr)$.
Since $\Ann\bigl(M(\lambda_\Theta)\bigr)
=\sum_{\Delta\in Z(\mathfrak g)}U(\mathfrak g)\bigr(\Delta-\Delta_{\mathfrak a}(\lambda_\Theta)\bigr)$,
$\mathcal V(I) \subset \{w.\lambda_\Theta;\, w\in W\}.$
Let $w=w(\Theta)w_\Theta\in W$ with $w(\Theta)\in W(\Theta)$ and $w_\Theta\in W_\Theta$.
Suppose $w(\Theta)\ne e$.  
Then there exists $\alpha\in\Psi(\mathfrak g)$
such that the length of 
$w_\alpha w(\Theta)$ is less than that of $w(\Theta)$.
For this root $\alpha$ we have 
$w_\alpha w(\Theta)\in W(\Theta)$ and
$w(\Theta)^{-1}\alpha, w_\Theta^{-1}w(\Theta)^{-1}\alpha\in\Sigma(\mathfrak g)^-\setminus\Sigma(\mathfrak g_\Theta)$.
The assumption thereby implies
\[
 2\frac{\ang{w.\lambda_\Theta+\rho}{\alpha}}{\ang{\alpha}{\alpha}}\notin\{1,2,3,\ldots\}.
\]
Hence $(w_\alpha w).\lambda_\Theta\in \mathcal{V}(I)$ provided that $w.\lambda_\Theta\in \mathcal{V}(I)$,
which assures 
\begin{equation}\label{eq:roots_disj}
\mathcal{V}(I)\cap
 \bigl\{\bigl(W(\Theta)w_\Theta \bigr).\lambda_\Theta;\,
 w_\Theta\in W_\Theta\setminus\{e\} \bigr\}=\emptyset
\end{equation}
by ii) and the induction on the length of $w(\Theta)$.
Similarly we can show that $\mathcal{V}(I)\supset\{w.\lambda_\Theta;\, w\in W(\Theta)\}$
if 
\begin{equation}\label{eq:regroot}
 2\frac{\ang{\lambda_\Theta+\rho}{\alpha}}{\ang{\alpha}{\alpha}}\notin\{1,2,3,\ldots\}
 \quad\bigl(\forall\alpha\in
 \Sigma(\mathfrak g)^+\setminus\Sigma(\mathfrak g_\Theta)\bigr).
\end{equation}
Let us remove the condition~\eqref{eq:regroot} by Proposition~\ref{prop:ann_str}.
Since $U(\mathfrak g)=\mathcal{H}\oplus \Ann\bigl(M(\lambda_\Theta)\bigr)$,
we have only to show for each finite dimensional $\mathfrak g$-module $\mathbf{V}$
\begin{equation}\label{eq:phivanish}
 \bigl(\varphi(v)\bigr)_{\mathfrak a}(w.\lambda_\Theta)=0\quad\Bigl(
 \forall\varphi\in
 \Hom_{\mathfrak g}\left(\mathbf{V}, \mathcal{H}\cap\Ann\bigl(M_\Theta(\lambda)\bigr)\right),
\forall v\in\mathbf{V}, \forall w\in W(\Theta)\Bigr).
\end{equation}
For $D\in U(\mathfrak g)$ we denote by $D^\lambda$
a unique element of $U(\bar{\mathfrak n}_\Theta)$ such that $D-D^\lambda\in J_\Theta(\lambda)$. 
Then $\varphi\in\Hom_{\mathfrak g}\left(\mathbf{V}, \mathcal{H}\right)$
belongs to $\Hom_{\mathfrak g}\left(\mathbf{V}, \mathcal{H}\cap\Ann\bigl(M_\Theta(\lambda)\bigr)\right)$
if and only if $\varphi(v)^\lambda=0$ for $v\in V$.
Let $k=\dim\mathbf{V}_0$ and take 
$\varphi_1,\ldots,\varphi_k \in \Hom_{\mathfrak g}\left(\mathbf{V}, \mathcal{H}\right)$ 
so that they constitute a basis.
Note that for $v\in\mathbf{V}$ and $i=1,\ldots,k$,
$\varphi_i(v)^\lambda$ are 
$U(\bar{\mathfrak n}_\Theta)$-valued polynomials in $\lambda$.
Let $\ell=k-\dim\mathbf{V}^{\mathfrak g_\Theta}$.
Then by Proposition~\ref{prop:ann_str}
there exist an open neighborhood $S\subset\mathfrak a_\Theta^*$ of the point in question
and complex-valued rational functions $a_{ij}(\lambda)$ on $S$ such that
\[
a_{1j}(\lambda)\varphi_1+a_{2j}(\lambda)\varphi_2+\cdots+a_{kj}(\lambda)\varphi_k
\quad(j=1,\ldots,\ell)
\]
form a basis of $\Hom_{\mathfrak g}\left(\mathbf{V}, \mathcal{H}\cap\Ann\bigl(M_\Theta(\lambda)\bigr)\right)$
for any $\lambda\in S$.
Since generic $\lambda\in S$ satisfy \eqref{eq:regroot},
\eqref{eq:phivanish} holds for any $\lambda\in S$.
\end{proof}
On the existence of a two-sided ideal $I_\Theta(\lambda)$ satisfying
\eqref{eq:gapdef}, we have 
\begin{thm}\label{thm:gapexist}
Suppose $\lambda_\Theta+\rho$ is dominant.
Then the following four conditions are equivalent.\smallskip

\noindent {\rm i)}
$J_\Theta(\lambda)=\Ann\bigl(M_\Theta(\lambda)\bigr) + J(\lambda_\Theta)$.

\noindent {\rm ii)}
If $\beta\in\Sigma(\mathfrak g)^+\setminus\Sigma(\mathfrak g_\Theta)$
satisfies $\ang{\lambda_\Theta+\rho}{\beta}=0$,
then $\ang{\beta}{\alpha}=0$ for all $\alpha\in\Theta$.

\noindent {\rm iii)}
$W(\Theta).\lambda_\Theta \cap W_\Theta.\lambda_\Theta=\{\lambda_\Theta\}$.

\noindent {\rm iv)}
If $w_\Theta\in W_\Theta$ satisfies
$\bigl(W(\Theta)w_\Theta\bigr).\lambda_\Theta \cap W(\Theta).\lambda_\Theta\ne\emptyset$,
then $w_\Theta=e$.\smallskip

\noindent
In particular, if $\lambda_\Theta+\rho$ is regular,
these conditions are satisfied.
\end{thm}
\begin{proof}
iv) $\Rightarrow$ iii) is obvious.

iii) $\Rightarrow$ ii).
Suppose there exist $\beta\in\Sigma(\mathfrak g)^+\setminus\Sigma(\mathfrak g_\Theta)$
and $\alpha\in\Theta$ such that 
$\ang{\lambda_\Theta+\rho}{\beta}=0$ and $\ang{\beta}{\alpha}\ne0$.
For $\gamma\in\Sigma(\mathfrak g_\Theta)^+$
we have
\[
2\frac{\ang{\lambda_\Theta+\rho}{w_\beta\gamma}}{\ang{w_\beta\gamma}{w_\beta\gamma}}
=2\frac{\ang{\lambda_\Theta+\rho}{\gamma}}{\ang{\gamma}{\gamma}}
=2\frac{\ang{\rho}{\gamma}}{\ang{\gamma}{\gamma}}
\in\{1,2,\ldots\},
\]
which shows $\ang{\beta}{\gamma}\leq 0$ and $w_\beta\in W(\Theta)$.
In particular $\ang{\beta}{\alpha}< 0$ and hence $w_\alpha w_\beta \in W(\Theta)$.
Now we get $(w_\alpha w_\beta).\lambda_\Theta=w_\alpha.\lambda_\Theta$,
a contradiction.

ii) $\Rightarrow$ i).
For each $\alpha\in\Theta$ we define the $\mathfrak g$-homomorphism
$M(\lambda_\Theta-\alpha)\rightarrow M(\lambda_\Theta)$
by $D\bmod J(\lambda_\Theta-\alpha)\mapsto DX_{-\alpha} \bmod J(\lambda_\Theta)$.
This is an injection and therefore we identify its image with $M(\lambda_\Theta-\alpha)$.
Note that
\[
\sum_{\alpha\in\Theta}M(\lambda_\Theta-\alpha)=
\Bigl(J(\lambda_\Theta)+\sum_{\alpha\in\Theta}U(\mathfrak g)X_{-\alpha} \Bigr)
/ J(\lambda_\Theta) =J_\Theta(\lambda) / J(\lambda_\Theta)
\]
and we have a surjection $P(\lambda_\Theta-\alpha)\rightarrow M(\lambda_\Theta-\alpha)$
by Proposition~\ref{prop:Ofacts}~ii).
Moreover it is clear that the condition \eqref{eq:proj_rep}
with $(\mu,\mu')=(\lambda_\Theta,\lambda_\Theta-\alpha)$
holds for each $\alpha\in\Theta$.
Hence by Proposition~\ref{prop:2toleft}
we have a two-sided ideal $I$ containing $\Ann\bigl(M(\lambda_\Theta)\bigr)$
such that $IM(\lambda_\Theta)=J_\Theta(\lambda) / J(\lambda_\Theta)$.
Then $I=\Ann\bigl(M_\Theta(\lambda)\bigr)$ and $J_\Theta(\lambda)=I+J(\lambda_\Theta)$.

i) $\Rightarrow$ iv) follows from \eqref{eq:Annroots} and \eqref{eq:roots_disj}.
\end{proof}
\begin{rem}
Through $I_{\pi,\Theta}$,
we will get in \S\ref{sec:exm}
many sufficient conditions for \eqref{eq:gapann},
which are effective even if $\lambda_\Theta+\rho$ is not dominant.
\end{rem}
\begin{defn}[extremal low weight]
For a simple root $\alpha\in\Psi(\mathfrak g)$, 
we call a minimal element of 
$\{\varpi\in\mathcal W(\pi);\,\ang{\varpi}{\alpha} \ne 0\}$ 
under the order $\le$ in Definition~\ref{def:shift} {\it an extremal low weight of $\pi$ 
with respect to\/} $\alpha$.
\end{defn}
Since $\pi$ is a faithful representation, $\pi(X_{-\alpha})$ is not zero
and therefore an extremal low weight $\varpi_\alpha$ with respect to $\alpha$ always exists
but it may not be unique.
The main purpose in this section is to calculate 
the function 
\begin{equation}\label{eqn:r_func}
\mathfrak a_\Theta^*\ni\lambda\mapsto
 \Bigl(q_{\pi,\Theta}(F_\pi;\lambda)_{\varpi_\alpha\varpi_\alpha}
 \Bigr)_{\mathfrak a}(\lambda_\Theta-\alpha)
\end{equation}
on $\mathfrak a_\Theta^*$.
If for any $\alpha\in\Theta$ there exists $\varpi_\alpha$ such that 
the value of the corresponding function \eqref{eqn:r_func} does not vanish, Lemma~\ref{lem:gap} 
assures \eqref{eq:gap}.
\begin{lem}\label{lem:maxroot}
Fix $\alpha\in\Psi(\mathfrak g)$ and let $\varpi_\alpha$ be an extremal low weight of $\pi$
with respect to $\alpha$.
For $\lambda=\sum_{\beta\in\Psi(\mathfrak g)}m_\beta\beta\in\mathfrak a^*$ put
$|\lambda| = \sum_{\beta\in\Psi(\mathfrak g)}m_\beta$.
Then there exists $\{\gamma_1,\ldots,\gamma_K\}\subset\Psi(\mathfrak g)$
with $\gamma_K=\alpha$
such that the following \eqref{eq:maxL1}--\eqref{eq:maxL8} hold by denoting 
\begin{equation}\label{eq:maxL0}
\varpi_i =\varpi_\alpha - \sum_{i\le\nu< K}\gamma_\nu.
\end{equation}
\begin{align}
 K = |\varpi_\alpha &- \bar\pi| +1 \text{ and }\varpi_1 = \bar\pi,\label{eq:maxL1}\\
 \ang{\varpi_i}{\gamma_i} &< 0\text{ for }i=1,\ldots,K,\label{eq:maxL2}\\
 \ang{\varpi_i}{\gamma_j} &= 0\text{ if } 1\le i<j\le K,\label{eq:maxL3}\\
 \ang{\gamma_i}{\gamma_j} &\ne 0\text{ if and only if }|i-j|\le 1,
  \label{eq:maxL4}\\
 \{\varpi_1,\ldots,\varpi_{K-1}\}
   &=\{\varpi'\in\mathcal W(\pi);\,\varpi'<\varpi_\alpha\},\label{eq:maxL5}\\
 \varpi_i\text{ is an extrem}&\text{al low weight of }\pi
 \text{ with respect to }\gamma_i
 \text{ for }i=1,\ldots,K,\label{eq:maxL6}\\
 \text{the multiplicity }&\text{of the weight space of the weight $\varpi_i$
  equals $1$}.\label{eq:maxL8}
\end{align}
The sequence $\gamma_1,\ldots,\gamma_K$ is unique 
by the condition $\varpi_1,\ldots,\varpi_K\in\mathcal W(\pi)$.
The part of the partially ordered set of the weights of $\pi$ 
which are smaller or equal to $\varpi_\alpha$ is as follows:
\begin{equation}
 \varpi_1=\bar\pi\xrightarrow{\gamma_1}\varpi_2\xrightarrow{\gamma_2}\varpi_3
 \xrightarrow{\gamma_3}\cdots\cdots\xrightarrow{\gamma_{K-1}}
 \varpi_K=\varpi_\alpha\xrightarrow{\gamma_K=\alpha}
\end{equation}
\end{lem}
\begin{proof}
Let $\gamma_1,\ldots,\gamma_K$ be a sequence of $\Psi(\mathfrak g)$
satisfying \eqref{eq:maxL1},
$\gamma_K=\alpha$, and $\varpi_1,\ldots,\varpi_K\in\mathcal W(\pi)$
under the notation \eqref{eq:maxL0}.
The existence of such a sequence is clear.
We shall prove by the induction on $K$
that such a sequence is unique and 
that it satisfies \eqref{eq:maxL2}--\eqref{eq:maxL6}.

By the minimality of $\varpi_\alpha$ we have $\ang{\varpi_i}{\alpha}=0$
for $i = 1, \ldots, K-1$.
Hence $\ang{\gamma_i}{\alpha}=\ang{\varpi_{i+1}-\varpi_i}{\alpha}=0$
for $i = 1, \ldots, K-2$
and
$\ang{\gamma_{K-1}}{\alpha}=
\ang{\varpi_\alpha-\varpi_{K-1}}{\alpha}=\ang{\varpi_\alpha}{\alpha}<0$.
Thus we get $\gamma_i \ne \alpha$
for $i=1,\ldots,K-1$.
Moreover $\varpi_\alpha - \gamma_i \notin \mathcal W(\pi)$ 
for $i=1,\ldots,K-2$
because $\ang{\varpi_\alpha - \gamma_i}{\alpha}=\ang{\varpi_\alpha}{\alpha}\ne0$
and $\varpi_\alpha$ is minimal.
This means
$
 \{\varpi'\in \mathcal W(\pi);\,\varpi'<\varpi_\alpha\}
 = \{\varpi_{K-1}\}\cup
 \{\varpi'\in \mathcal W(\pi);\,\varpi'<\varpi_{K-1}\}
$.

Suppose $\ang{\varpi_{K-1}}{\gamma_{K-1}}\ge0$.
Then $\varpi_{K-1}-\gamma_{K-1}\in\mathcal W(\pi)$ because
$\varpi_{K-1}+\gamma_{K-1} = \varpi_\alpha\in\mathcal W(\pi)$.
Hence $\ang{\varpi_{K-1}-\gamma_{K-1}}{\alpha} = -\ang{\gamma_{K-1}}\alpha > 0$,
which contradicts with the minimality of $\varpi_\alpha$.
Thus we get $\ang{\varpi_{K-1}}{\gamma_{K-1}}<0$.

Suppose $\varpi_{K-1}$ is not an extremal low weight with respect to
$\gamma_{K-1}$. 
Then there exists an extremal low weight $\varpi'$ with respect to 
$\gamma_{K-1}$ such that $\varpi'<\varpi_{K-1}$.
Then $\mathcal W(\pi)\ni\varpi' + \gamma_{K-1} < \varpi_\alpha$ and 
$\ang{\varpi'+\gamma_{K-1}}{\alpha} = \ang{\gamma_{K-1}}{\alpha}= 0$
by the minimality of $\varpi_\alpha$.
It is a contradiction.
Hence $\varpi_{K-1}$ is an extremal low weight with respect
to $\gamma_{K-1}$.

Now by the induction hypothesis we obtain the uniqueness and \eqref{eq:maxL2}--\eqref{eq:maxL6}.
Note that \eqref{eq:maxL8} follows from the uniqueness and the following lemma
because $V=U(\mathfrak n)v_{\bar\pi}$ with a lowest weight vector $v_{\bar\pi}$
of $\pi$.
\end{proof}
\begin{lem}
$U(\mathfrak n)$ is generated by $\{X_\gamma;\,\gamma\in\Psi(\mathfrak g)\}$
as a subalgebra of $U(\mathfrak g)$.
\end{lem}
\begin{proof}
Let $U$ denote the algebra generated by
$\{X_\gamma;\,\gamma\in\Psi(\mathfrak g)\}$.
It is sufficient to show that $X_\beta\in U$ for $\beta\in\Sigma(\mathfrak g)^+$,
which is proved by the induction on $|\beta|$ as follows. 
If $|\beta|>1$, there exists $\gamma\in\Psi(\mathfrak g)$ such that 
$\beta'=\beta-\gamma\in\Sigma(\mathfrak g)^+$.
Then $X_\beta=C(X_\gamma X_{\beta'} - X_{\beta'} X_\gamma)$ with a constant
$C\in\mathbb C$.  Hence the condition $X_\gamma,\,X_{\beta'}\in U$ implies 
$X_\beta\in U$.
\end{proof}
\begin{rem}\label{rem:lowweight}
By virtue of \eqref{eq:maxL4} the Dynkin diagram of 
the system $\{\gamma_1,\ldots,\gamma_K\}$ in Lemma~\ref{lem:maxroot}
is of type $A_K$ or $B_K$ or $C_K$ or $F_4$ or $G_2$
where $\gamma_1$ and $\gamma_K$ correspond to the end points of 
the diagram.  
Note that 
\begin{equation}\label{eq:simseq}
\ang{\bar\pi}{\gamma_1}<0 \text{ and } 
\ang{\bar\pi}{\gamma_i}=0 \text{ for } i=2,\ldots,K.
\end{equation}
Conversely if a subsystem $\{\gamma_1,\ldots,\gamma_K\}\subset\Psi(\mathfrak g)$
satisfies \eqref{eq:maxL4} and \eqref{eq:simseq} then
$\bar\pi+\gamma_1+\cdots+\gamma_{K-1}$ is an extremal low weight with respect to $\gamma_K$.
Hence we have at most three different extremal 
low weights of $\pi$ with respect to a fixed $\alpha\in\Psi(\mathfrak g)$.
\end{rem}

The next lemma is studied in \cite[Lemma~3.5]{O-Cl}.
It gives the solutions for the recursive equations
which play key roles in the calculation of \eqref{eqn:r_func}.

\begin{lem}\label{lem:recursion}
For $k=0, 1, \ldots$ and $\ell=1,2,\ldots$,
define the polynomial $f(k,\ell)$ in the variables
$s_1, \ldots, s_{\ell-1}, \mu_1, \mu_2, \ldots$
recursively by 
\begin{equation}\label{eq:zenka_f}
 f(k,\ell) = \begin{cases}
     1
     &\, \text{if \ }k=0,\\
    f(k-1,\ell)(\mu_\ell - \mu_k)
    + \sum\limits_{\nu=1}^{\ell - 1} s_\nu f(k-1,\nu)
    &\, \text{if \ }k\ge1.\\
    \end{cases}
\end{equation}
Moreover for $k=1, 2, \ldots$ and $\ell=1,2,\ldots$,
define the polynomial $g(k,\ell)$ in the variables
$t,s_1, \ldots, s_{\ell-1}, \mu_1, \mu_2, \ldots$
recursively by 
\begin{equation}\label{eq:zenka_g}
 g(k,\ell) = \begin{cases}
     1&\quad\quad\quad\text{\ \ if \ }k=1,\\
     g(k-1,\ell)(t - \mu_k)
     + f(k-1,\ell)
     &\quad\quad\quad\text{\ \ if \ }k>1.\\
    \end{cases}
\end{equation}
Then the following \eqref{eq:zenka_f1}--\eqref{eq:zenka_g1} hold.
\begin{align}
f(k,\ell)&=0\text{\ \ for \ }k\ge\ell,\label{eq:zenka_f1}\\
f(\ell-1,\ell)&=\prod_{\nu=1}^{\ell-1}(\mu_\ell-\mu_\nu+s_\nu),\label{eq:zenka_f2}\\
g(k,\ell)&=\prod_{\nu=1}^{\ell-1}(t-\mu_\nu+s_\nu)
\prod_{\nu=\ell+1}^k(t-\mu_\nu)\text{\ \ for \ }k\ge\ell.\label{eq:zenka_g1}
\end{align}
\end{lem}

Now recall \eqref{eq:zenka} with $\Theta=\emptyset$.
Let $F_{ii}^k\in U(\mathfrak a)$ be the element in \eqref{eq:zenka} 
corresponding to the weight $\varpi_i$ for $i=1,\ldots,K$ under the notation in 
Lemma~\ref{lem:maxroot}. 
Then Lemma~\ref{lem:maxroot} and Lemma~\ref{lem:casimir} iv) with $\ell=1$, 
$\beta=\varpi_i-\varpi_\nu\in\Sigma(\mathfrak g)^+$ $(1\le\nu<i)$ and 
$\varpi= \varpi_\nu$ show that \eqref{eq:zenka} is reduced to
\begin{multline}\label{eq:Fi_zenka}
  F_{ii}^k - F_{ii}^{k-1}\bigl(\varpi_i - \mu_k 
   + D_\pi(\varpi_i)\bigr)\\
  \equiv
  \sum_{1\le\nu<i}\ang{\varpi_\nu}{\varpi_i-\varpi_\nu}F_{\nu\nu}^{k-1}
  \mod U(\mathfrak g)\mathfrak n.
\end{multline}
Since $\ang{\varpi_i}{\lambda_\Theta}=\ang{\varpi_i}{\lambda_\Theta-\alpha}$
for $i=1,\ldots,K-1$, 
\eqref{eq:Fi_zenka} inductively implies
\begin{equation}\label{eq:lowF}
\bigl(F_{ii}^k\bigr)_{\mathfrak a}(\lambda_\Theta)=
\bigl(F_{ii}^k\bigr)_{\mathfrak a}(\lambda_\Theta-\alpha)
\quad\text{ for }i=1,\ldots,K-1\text{ and }k=0,1,\ldots.
\end{equation}
From \eqref{eq:maxL3} we have
\[
\ang{\varpi_\nu}{\varpi_i-\varpi_\nu}
 = \ang{\varpi_\nu}{\gamma_\nu+\cdots+\gamma_{i-1}}
 = \ang{\varpi_\nu}{\gamma_\nu}
\]
and hence
\begin{align*}
F_{i+1i+1}^k - F_{ii}^k&\equiv
F_{i+1i+1}^{k-1}\bigl(\varpi_{i+1} - \mu_k 
   + D_\pi(\varpi_{i+1})\bigr)\\
&\quad +F_{ii}^{k-1}\ang{\varpi_{i}}{\varpi_{i+1}-\varpi_{i}}
-F_{ii}^{k-1}\bigl(\varpi_i - \mu_k 
   + D_\pi(\varpi_i)\bigr)
\mod U(\mathfrak g)\mathfrak n\\
&=(F_{i+1i+1}^{k-1}-F_{ii}^{k-1})\bigl(\varpi_{i+1} - \mu_k +
  D_\pi(\varpi_{i+1})
  \bigr) + F_{ii}^{k-1}\gamma_i.
\end{align*}
The last equality above follows from Lemma~\ref{lem:order} i)
with $\varpi=\varpi_i$ and $\varpi'=\varpi_{i+1}$ 
because $\gamma_i = \varpi_{i+1}-\varpi_{i}\in\Psi(\mathfrak g)$.
Hence by the induction on $k$ we have
\[
 F_{i+1i+1}^k\equiv F_{ii}^k\mod U(\mathfrak g)\mathfrak n
        + U(\mathfrak g)\gamma_i.
\]

Now consider general $\Theta\subset\Psi(\mathfrak g)$.
Define integers $n_0$, $n_1,\ldots,n_L$ with 
$n_0=0<n_1<\cdots<n_L=K$
such that
\[
 \{n_1,\ldots,n_{L-1}\} = \{\nu\in\{1,\ldots,K-1\};
 \gamma_\nu\notin\Theta\}.
\]

If $n_{\ell-1}<\nu<n_\ell$,
then $\gamma_\nu\in\Theta$,
which implies  
$\ang{\gamma_{\nu}}{\lambda_\Theta} = 0$ and hence
\[
\bigl(F_{\nu+1\nu+1}^k\bigr)_{\mathfrak a}(\lambda_\Theta)
=\bigl(F_{\nu\nu}^k\bigr)_{\mathfrak a}(\lambda_\Theta).
\]
We note that
$\varpi_{n_0+1}|_{\mathfrak a_\Theta}<_\Theta
\varpi_{n_1+1}|_{\mathfrak a_\Theta}<_\Theta\cdots<_\Theta
\varpi_{n_{L-1}+1}|_{\mathfrak a_\Theta}$ and
\[
 \{\varpi_{n_0+1},\ldots,\varpi_{n_{L-1}+1}\}
 = \{\varpi'\in\overline{\mathcal W}_\Theta(\pi);\,
   \varpi'\le\varpi_\alpha\}.
\]

Put $\mu_\ell=\ang{\varpi_{n_{\ell-1}+1}}{\lambda_\Theta}+D_\pi(\varpi_{n_{\ell-1}+1})$
for $\ell=1,\ldots,L$.
Since $\prod_{\ell=1}^L(x-\mu_\ell)$ is a divisor of $q_{\pi,\Theta}(x;\lambda)$,
we can take $\mu_\ell$ for $\ell=L+1, L+2, \ldots, L'=\deg_x q_{\pi,\Theta}(x;\lambda)$
so that $q_{\pi,\Theta}(x;\lambda)=\prod_{\ell=1}^{L'}(x-\mu_\ell).$

For $k=0, 1, \ldots,L'$ and $\ell=1, 2, \ldots, L$ we define 
\[
f(k,\ell)=\bigl(F_{n_{\ell-1}+1,n_{\ell-1}+1}^k\bigr)_{\mathfrak a}(\lambda_\Theta)
=\cdots=\bigl(F_{n_\ell,n_\ell}^k\bigr)_{\mathfrak a}(\lambda_\Theta).
\]
Then putting
\[s_\ell = \sum_{n_{\ell-1}<\nu\le n_\ell}\ang{\varpi_\nu}{\gamma_{\nu}},\]
we have from \eqref{eq:Fi_zenka} with $i=n_{\ell-1}+1$
\[
f(k,\ell) = f(k-1,\ell)(\mu_\ell - \mu_k)+\sum_{j=1}^{\ell-1}s_j f(k-1,j).
\]
From
\eqref{eq:lowF} and \eqref{eq:Fi_zenka} with $i=n_L=K$
we also have 
\begin{multline*}
\bigl(F_{KK}^k\bigr)_{\mathfrak a}(\lambda_\Theta-\alpha) = 
\bigl(F_{KK}^{k-1}\bigr)_{\mathfrak a}(\lambda_\Theta-\alpha)
\bigl(\ang{\varpi_\alpha}{\lambda_\Theta-\alpha}
+D_\pi(\varpi_\alpha)-\mu_k\bigr)\\
+\sum_{j=1}^{L-1}s_j f(k-1,j)
+\left(\sum_{\nu=n_{L-1}+1}^{K-1}\ang{\varpi_\nu}{\gamma_\nu}\right) f(k-1,L).
\end{multline*}
Hence by Lemma~\ref{lem:order} i)
\begin{multline*}
\frac{f(k,L)-\bigl(F_{KK}^k\bigr)_{\mathfrak a}(\lambda_\Theta-\alpha)}{\ang{\varpi_\alpha}{\alpha}}\\
=\frac{f(k-1,L)-\bigl(F_{KK}^{k-1}\bigr)_{\mathfrak a}(\lambda_\Theta-\alpha)}
{\ang{\varpi_\alpha}{\alpha}}
\bigl(\ang{\varpi_\alpha}{\lambda_\Theta-\alpha}
+D_\pi(\varpi_\alpha)-\mu_k\bigr)\\
+f(k-1,L).
\end{multline*}
Now applying Lemma~\ref{lem:recursion}
to
\[g(k,L)=\frac{f(k,L)-\bigl(F_{KK}^k\bigr)_{\mathfrak a}(\lambda_\Theta-\alpha)}{\ang{\varpi_\alpha}{\alpha}}\]
with $t=\ang{\varpi_\alpha}{\lambda_\Theta-\alpha}
+D_\pi(\varpi_\alpha)$,
we obtain 
\[
\begin{aligned}
  \Bigr(q_{\pi,\Theta}(F_\pi;\lambda)_{\varpi_\alpha\varpi_\alpha}\Bigr)_{\mathfrak a}
  &(\lambda_\Theta-\alpha) \\
  &= 
  \left(F_{KK}^{L'}\right)_{\mathfrak a}
  (\lambda_\Theta-\alpha)\\
  &=-\ang{\varpi_\alpha}{\alpha}
  \prod_{\ell=1}^{L-1} \Bigl(\ang{\varpi_\alpha}{\lambda_\Theta-\alpha}
+D_\pi(\varpi_\alpha)-\mu_\ell+s_\ell\Bigr)\\
&\qquad \cdot
\prod_{\ell=L+1}^{L'}
\Bigl(\ang{\varpi_\alpha}{\lambda_\Theta-\alpha}
+D_\pi(\varpi_\alpha)-\mu_\ell\Bigr) \\
&=-\ang{\varpi_\alpha}{\alpha}
  \prod_{\ell=1}^{L-1}
  \Bigl(
    \ang{\varpi'_\alpha-\varpi_{n_\ell}}{\lambda_\Theta}
    +D_\pi(\varpi'_\alpha) -D_\pi(\varpi_{n_\ell+1})
  \Bigr)\\
&\qquad \cdot \!\!\!\!\!\! \prod_{
     (\mu,C)\in\Omega_{\pi,\Theta}\setminus
     \Omega_{\pi,\Theta}^{\varpi_\alpha}}
     \Bigl(
    \ang{\varpi_\alpha'-\mu}{\lambda_\Theta}+D_\pi(\varpi_\alpha')-C
    \Bigr).
\end{aligned}\]
Here we put $\varpi'_\alpha=\varpi_\alpha+\alpha\in\mathcal{W}(\pi)$ and
\begin{equation}\label{eq:partialOmega}
 \Omega_{\pi,\Theta}^{\varpi_0} = 
  \{(\varpi|_{\mathfrak a_\Theta},D_\pi(\varpi));\,
    \varpi\in\overline{\mathcal W}_\Theta(\pi),\ \varpi\le\varpi_0\}
\end{equation}
for $\varpi_0\in\mathcal{W}(\pi)$.
To deduce the last equality,
we have used
\[
\mu_\ell - s_\ell = \ang{\varpi_{n_\ell}}{\lambda_\Theta}+
  D_\pi(\varpi_{n_\ell+1}) \qquad
\text{ if } 1\le \ell \le L-1.\]

\begin{defn}\label{def:gap}
Suppose $\alpha\in\Theta$ and
$\varpi_\alpha$ is an extremal low weight of $\pi$ with 
respect to $\alpha$. 
Put $\varpi_\alpha'=\varpi_\alpha+\alpha\in\mathcal W(\pi)$ and
\[
 \{\varpi_1,\ldots,\varpi_K\} = \{\varpi\in\overline{\mathcal W}(\pi);
  \,\varpi\le\varpi_\alpha\}
\]
with $\varpi_1<\varpi_2<\cdots<\varpi_K$ and define $n_0=0<n_1<\cdots<n_L<K$ so that
\[
 \{\varpi_{n_0+1},\ldots,\varpi_{n_L+1}\} = 
 \{\varpi\in\overline{\mathcal W}_\Theta(\pi);\,\varpi \le \varpi_\alpha\}.
\]
Under the notation in Definition~\ref{def:shift} and \eqref{eq:partialOmega}, define
\begin{multline}\label{eq:gapcond}
 r_{\alpha,\varpi_\alpha}(\lambda) = 
   \prod_{(\mu,C)\in\Omega_{\pi,\Theta}\setminus\Omega_{\pi,\Theta}^{\varpi_\alpha}}
    \Bigl(
    \ang{\lambda_\Theta}{\varpi_\alpha'-\mu}+D_\pi(\varpi_\alpha')-C
    \Bigr)\\
  \cdot\prod_{i=1}^{L}
   \Bigl(
    \ang{\lambda_\Theta}{\varpi_\alpha-\varpi_{n_i}}
    -\ang{\alpha}{\varpi_\alpha}
    +D_\pi(\varpi_\alpha) -D_\pi(\varpi_{n_i+1})
    \Bigr).
\end{multline}
If there is no extremal low weights with respect to $\alpha$
other than $\varpi_\alpha$, 
we use the simple symbol $r_\alpha(\lambda)$ for $r_{\alpha,\varpi_\alpha}(\lambda)$.
\end{defn}
\begin{rem}\label{rm:Agap}
In the above definition we have the following.

\noindent i)
If the lowest weight $\bar\pi$ is an extremal low weight of $\pi$ 
with respect to $\alpha$, then $L=0$.

\noindent ii)
The second factor
\[
  \prod_{i=1}^{L}
   \Bigl(
    \ang{\lambda_\Theta}{\varpi_\alpha-\varpi_{n_i}}
    -\ang{\alpha}{\varpi_\alpha}
    +D_\pi(\varpi_\alpha) -D_\pi(\varpi_{n_i+1})
   \Bigr)
\]
is not identically zero because 
$\varpi_{n_i}|_{\mathfrak a_\Theta} <_\Theta \varpi_{n_i+1}|_{\mathfrak a_\Theta}
\leq_\Theta \varpi_\alpha|_{\mathfrak a_\Theta}$.

\noindent iii)
For $\varpi$ and $\varpi'\in\mathcal W(\pi)$
\begin{equation}
 \ang{\lambda_\Theta}{\varpi-\varpi'}+D_\pi(\varpi) - D_\pi(\varpi')
 =\ang{\lambda_\Theta+\rho}{\varpi-\varpi'}
 +\frac{\ang{\varpi'}{\varpi'}-\ang{\varpi}{\varpi}}2.
\end{equation}

\noindent iv) 
Put $\gamma_\nu=\varpi_{\nu+1}-\varpi_\nu$ for $\nu=1,\ldots,K-1$
and $\gamma_K=\alpha$.
If
\begin{equation}\label{eq:eqabs}
-2\frac{\ang{\varpi_\nu}{\gamma_\nu}}{\ang{\gamma_\nu}{\gamma_\nu}}\ 
\Bigl(= -2\frac{\ang{\gamma_{\nu-1}}{\gamma_\nu}}{\ang{\gamma_\nu}{\gamma_\nu}}
 \text{ if }\nu>1\Bigr) = 1,
\end{equation}
then $\ang{\varpi_\nu}{\varpi_\nu} = \ang{\varpi_{\nu+1}}{\varpi_{\nu+1}}$.

\noindent v)
Suppose 
$
 \frac{2\ang{\bar\pi}{\gamma_1}}{\ang{\gamma_1}{\gamma_1}}=-1
$
and the Dynkin diagram of the system $\{\gamma_1,\ldots,\gamma_K\}$
is of type $A_K$ or of type $B_K$ with short root $\gamma_K$
or of type $G_2$ with short root $\gamma_2$.
Then it follows from Lemma~\ref{lem:maxroot} and Lemma~\ref{lem:order}~i) that
\begin{multline}\label{eq:pregap}
 \ang{\lambda_\Theta}{\varpi_\alpha-\varpi_{n_i}}-
   \ang{\alpha}{\varpi_\alpha}+D_\pi(\varpi_\alpha)-D_\pi(\varpi_{n_i+1})\\
 =\ang{\lambda_\Theta}{\varpi_\alpha-\varpi_{n_i}}  + D_\pi(\varpi_\alpha)
  - D_\pi(\varpi_{n_i}) \\
 =\ang{\lambda_\Theta+\rho}{\varpi_\alpha-\varpi_{n_i}}
 =\ang{\lambda_\Theta+\rho}{\gamma_{n_i}+\cdots+\gamma_{K-1}} 
\end{multline}
for $i=1,\ldots,L$.

\end{rem}
\begin{thm}[gap]\label{thm:gap}
Let $\varpi_\alpha$ be an extremal low weight with 
respect to $\alpha\in\Theta$.
Then 
\[
  X_{-\alpha}\in I_{\pi,\Theta}(\lambda)
   + J(\lambda_\Theta)
  \quad\text{if }r_{\alpha,\varpi_\alpha}(\lambda)\ne0.
\]
If for all $\alpha\in\Theta$ there exists an extremal low weight $\varpi_\alpha$
with respect to $\alpha$ such that $r_{\alpha,\varpi_\alpha}(\lambda)\ne0$,
then 
\[ J_\Theta(\lambda) = 
  I_{\pi,\Theta}(\lambda)
  + J(\lambda_\Theta).
\]
\end{thm}
By Proposition~\ref{prop:HCzero}~iii)
we have the following corollary.
\begin{cor}[annihilator]\label{cor:ann}
If $\lambda_\Theta+\rho$ is dominant
and if for all $\alpha\in\Theta$ there exists an extremal low weight $\varpi_\alpha$
with respect to $\alpha$ such that $r_{\alpha,\varpi_\alpha}(\lambda)\ne0$, then
$I_{\pi,\Theta}(\lambda)=\Ann\bigl(M_\Theta(\lambda)\bigr)$.
\end{cor}

\begin{rem}\label{rem:badcase}
It does not always hold that for each $\alpha\in\Theta$
there exists an extremal low weight $\varpi_\alpha$
with respect to $\alpha$ such that the function $r_{\alpha,\varpi_\alpha}(\lambda)$
is not identically zero.
In fact we construct counter examples in Appendix~\ref{app:badcase}. %
However this condition is valid for many $\pi$
as we see below.
\end{rem}

Recall the notation in Proposition~\ref{prop:free}.
\begin{lem}\label{lem:free_nonzero}
Suppose $\varpi_\alpha$ is an extremal low weight with 
respect to $\alpha\in\Theta$.
The function $r_{\alpha,\varpi_\alpha}(\lambda)$ is not identically zero
if the space
\[ V(\varpi_\alpha|_{\mathfrak a_\Theta})=
   \sum_{\varpi\in\mathcal{W}(\pi);\,\varpi|_{\mathfrak a_\Theta}=\varpi_\alpha|_{\mathfrak a_\Theta}}
   \!\!\! V_\varpi
\]
is irreducible as a $\mathfrak g_\Theta$-module.
\end{lem}
\begin{proof}
In this case we have $\mu|_{\mathfrak a_\Theta} \ne \varpi_\alpha'|_{\mathfrak a_\Theta}$
for $(\mu,C)\in\Omega_{\pi,\Theta}\setminus\Omega_{\pi,\Theta}^{\varpi_\alpha}$
and the first factor of \eqref{eq:gapcond} is not identically zero.
\end{proof}
\begin{prop}\label{prop:goodweight}
Use the notation in\/ {\rm Lemma~\ref{lem:maxroot}}
and suppose $\gamma_K=\alpha\in\Theta$.
The function $r_{\alpha,\varpi_\alpha}(\lambda)$ is not identically zero
if either one of the following conditions is satisfied.

\noindent
{\rm i)}
$\{\gamma_1,\ldots,\gamma_K\} \subset \Theta$.

\noindent
{\rm ii)}
The connected component of the Dynkin diagram of\/ $\Theta$
containing $\alpha$ is orthogonal to $\bar{\pi}$.
$\Theta\setminus\{\gamma_1,\ldots,\gamma_K\}$
is orthogonal to $\{\gamma_1,\ldots,\gamma_{K-1}\}$.
Moreover
the Dynkin diagram of 
the system $\{\gamma_1,\ldots,\gamma_{K-1}\}$
is of type $A_{K-1}$.
\end{prop}
\begin{proof}
i)
Since $\varpi_\alpha|_{\mathfrak a_\Theta}=\bar{\pi}|_{\mathfrak a_\Theta}$
and 
$V(\bar{\pi}|_{\mathfrak a_\Theta})$
is an irreducible $\mathfrak g_\Theta$-module,
the claim follows from Lemma~\ref{lem:free_nonzero}.

ii)
Suppose $\varpi\in\overline{\mathcal W}_\Theta(\pi)$ satisfies
$\varpi|_{\mathfrak a_\Theta}=\varpi_\alpha|_{\mathfrak a_\Theta}$.
Then we can write
\[
\varpi=\bar{\pi}\ +\ \sum_{i=1}^K m_i \gamma_i\ + \!\!\!\!
 \sum_{\beta\in\Theta\setminus\{\gamma_1,\ldots,\gamma_K\}}
  n_\beta \beta
\]
with non-negative integers $m_i$ and $n_\beta$.
Put
\begin{align*}
\Theta'&=\{\gamma_i;\,m_i>0\},\\
\Theta''&=\{\beta;\,n_\beta>0\},
\end{align*}
and define
\[
V'=\sum\{V_{\varpi'};\,
\varpi'\in \bar\pi\ + \!\!\!\!
\sum_{\beta\in\Theta'\cup\Theta''}\mathbb{Z}\,\beta 
\}.
\]
Since $V'$ is an irreducible $\mathfrak{g}_{\Theta'\cup\Theta''}$-module
with lowest weight $\bar\pi$
and $\{0\} \subsetneq V_\varpi \subset V'$,
each connected component of 
the Dynkin diagram of the system $\Theta'\cup\Theta''$
is not orthogonal to $\bar{\pi}$.

Suppose $\gamma_K\in\Theta'$.
Then the condition ii) implies $\Theta'=\{\gamma_1, \ldots, \gamma_K\}$
and therefore $\varpi'_\alpha=\varpi_\alpha+\alpha\leq\varpi$.
However it is clear $\dim V_{\varpi'_\alpha}=1$
and $\varpi'_\alpha\notin\overline{\mathcal W}_\Theta(\pi)$.
Thus we have $\varpi'_\alpha<\varpi$.
In this case, by Lemma~\ref{lem:order}~ii),
we have $D(\varpi'_\alpha)<D(\varpi)$.

Suppose $\gamma_K\notin\Theta'$.
Then $\Theta'$ is orthogonal to $\Theta''$ and 
hence we have the direct sum decomposition
\[
\mathfrak{g}_{\Theta'\cup\Theta''}
=\mathfrak{a}_{\Theta'\cup\Theta''} \oplus \mathfrak{m}_{\Theta'}
 \oplus \mathfrak{m}_{\Theta''}.
\]
Since 
$\varpi$ is the lowest weight of a $\mathfrak{m}_{\Theta''}$-submodule
of $V'$, which is an irreducible $\mathfrak{m}_{\Theta'}
 \oplus \mathfrak{m}_{\Theta''}$-module,
$\Theta''$ must be empty.
On the other hand, we see
$\Theta'=\{\gamma_1, \ldots, \gamma_{K'}\}$
with $K'<K$.
Now we can find each weight $\varpi'$ of 
the $\mathfrak{g}_{\Theta'}$-module $V'$ 
is in the form
\[
\varpi'=\bar\pi\ +\ \sum_{i=1}^{K'} m'_i \gamma_i
\qquad\text{with }
-2\frac{\ang{\bar\pi}{\gamma_1}}{\ang{\gamma_1}{\gamma_1}}\geq
m'_1\geq m'_2\geq \cdots \geq m'_{K'}\geq 0
\]
and its multiplicity is one (cf.~Example~\ref{ex:A}~ii)).
Fix $v\in V_\varpi\setminus\{0\}$.
Take $i=1,\ldots,K'$ so that $m_i>m_{i+1}$.
Then $X_{-\gamma_i} v\ne 0$ and therefore
$\gamma_i\notin\Theta$. 
Since $\varpi|_{\mathfrak a_\Theta}=\varpi_\alpha|_{\mathfrak a_\Theta}$,
we conclude $i=K'$ and $m_{K'}=1$.
It shows
\[ \varpi = \gamma_1 + \cdots + \gamma_{K'} \leq \varpi_\alpha. \]

Thus we have proved the function \eqref{eq:gapcond} is not identically zero.
\end{proof}
\begin{rem}
The condition i) of the proposition is satisfied if the lowest weight 
$\bar\pi$ (or equivalently, the highest weight $\pi$) of $(\pi, V)$ is regular.
\end{rem}
\begin{prop}\label{prop:goodrep}

\noindent
{\rm i) (multiplicity free representation)}
Suppose $\dim V_\varpi=1$ for any $\varpi\in\mathcal W(\pi)$.
Then for any extremal low weight $\varpi_\alpha$
with respect to $\alpha\in\Theta$,
the function $r_{\alpha,\varpi_\alpha}(\lambda)$ is not identically zero.

\noindent
{\rm ii) (adjoint representation)}
Suppose $\mathfrak g$ is simple and 
$\pi$ is the adjoint representation of $\mathfrak g$.
Suppose $\alpha\in\Theta$.
If the Dynkin diagram of $\Psi(\mathfrak g)$ is of type $A_r$,
then we have just two extremal low weights $\varpi_\alpha$ with respect to $\alpha$.
If the diagram is not of type $A_r$,
then we have a unique $\varpi_\alpha$.
In either case, there is at least one $\varpi_\alpha$ 
such that $r_{\alpha,\varpi_\alpha}(\lambda)$ is not identically zero.

\noindent
{\rm iii) (minuscule representation)}
Suppose $(\pi, V)$ is minuscule.
Then for any $\alpha\in\Theta$
there is a unique extremal low weight $\varpi_\alpha$
with respect to $\alpha$.
Moreover the function $r_{\alpha}(\lambda)$ is not identically zero. 
\end{prop}
\begin{proof}
i) Thanks to Proposition~\ref{prop:free}~i),
$V(\varpi_\alpha|_{\mathfrak a_\Theta})$ is an irreducible $\mathfrak g_\Theta$-module.
Hence our claim follows from Lemma~\ref{lem:free_nonzero}. 

ii)
The lowest weight of the adjoint representation
is $-\alpha_{\max}$.
Hence by Remark~\ref{rem:lowweight} we can determine
the number of extremal low weights 
from the completed Dynkin diagram of each type,
which is shown in \S\ref{sec:exm}.

Note that $\mathcal{W}(\pi)=\Sigma(\mathfrak g)\cup\{0\}$.
Suppose $\varpi_\alpha\notin\Sigma(\mathfrak g_\Theta)$.
Then Proposition~\ref{prop:free}~ii) assures 
the irreducibility of
$V(\varpi_\alpha|_{\mathfrak a_\Theta})$.
Hence $r_{\alpha,\varpi_\alpha}(\lambda)$ is not identically zero.

Suppose $\varpi_\alpha\in\Sigma(\mathfrak g_\Theta)$.
Take $\{\gamma_1,\ldots,\gamma_K\}\subset\Psi(\mathfrak g)$ as in
Lemma~\ref{lem:maxroot} and put
\[
\varpi_i=-\alpha_{\max} + \gamma_1 + \cdots + \gamma_{i-1}
\qquad\text{for }i=1,\ldots,K.
\]
Let $\Theta_1$ denote the connected component of the Dynkin diagram of
$\Theta$ containing $\gamma_K=\alpha$.
Then we can find an integer $K'\in\{1,\ldots,K-1\}$
such that $\{\gamma_1,\ldots,\gamma_{K'}\}\subset \Psi(\mathfrak g)\setminus\Theta_1$
and $\{\gamma_{K'+1},\ldots,\gamma_K\}\subset \Theta_1$.
Then it follows from Lemma~\ref{lem:maxroot}
that the root vectors $X_{\varpi_i}$ for $i=1,\ldots,K'$
are lowest weight vectors of $\pi|_{\mathfrak m_{\Theta_1}}$.
These lowest weight vectors generate 
the irreducible ${\mathfrak m_{\Theta_1}}$-submodules
belonging to the same equivalence class
because $\{\gamma_1,\ldots,\gamma_{K'-1}\}$
is orthogonal to $\Theta_1$.
On the other hand, we have
$\varpi_{K'+1}\in\overline{\mathcal W}_\Theta(\pi)$.
Then it follows from Proposition~\ref{prop:free}~ii)
that $\varpi_{K'+1}\in\Sigma(\mathfrak{g}_{\Theta_1})^-$.
Since $\varpi_{K'}-\varpi_{K'+1}=-\gamma_{K'}\in\Sigma(\mathfrak g)$,
$[X_{-\varpi_{K'+1}}, X_{\varpi_{K'}}]\ne0$.
It shows the equivalence class above is not the class of the trivial representation.
Hence $\Theta_1$ is not orthogonal to $\bar\pi=\varpi_1$.
Now we can take another extremal low weight $\varpi'_\alpha$
with respect to $\alpha$ which satisfies
the condition i) of Proposition~\ref{prop:goodweight}.

iii)
Since a minuscule representation is of multiplicity free,
we have only to show the uniqueness of $\varpi_\alpha$.
Let $[\mathfrak g,\mathfrak g]=\mathfrak{g}_1\oplus\cdots\oplus\mathfrak{g}_m$
be the decomposition into simple Lie algebras.
Then $\pi|_{[\mathfrak g,\mathfrak g]}$ is a tensor product
of faithful minuscule representations of $\mathfrak{g}_i$
for $i=1,\ldots,m$.
Hence, from Proposition~\ref{prop:minuscule}~v),
each connected component of the Dynkin diagram of $\Psi(\mathfrak g)$,
which corresponds to some $\Psi(\mathfrak g_i)$,
has just one root $\gamma$ which is not orthogonal to $\bar\pi$.
Now the uniqueness follows from Remark~\ref{rem:lowweight}.
\end{proof}

We conclude this section with a discussion of the commutative case.
Consider $F_\pi = \Bigl(F_{ij}\Bigr)_{
  \substack{1\le i \le N\\1\le j \le N}
 }$ as an element of $M(N, S(\mathfrak g))$. 
Then we have
\begin{thm}[coadjoint orbit]\label{thm:orbit}
Put 
\begin{equation}
\begin{aligned}
 \bar \Omega_{\pi,\Theta}&=\{\varpi|_{\mathfrak a_\Theta};\,
   \varpi\in\overline{\mathcal W}_\Theta(\pi)\}\\
 \bar q_{\pi,\Theta}(x;\lambda) &= \prod_{\mu\in\bar\Omega_{\pi,\Theta}}
   \bigl(x-\mu(\lambda)\bigr),\\
 \bar r_\Theta(\lambda) &= 
  \prod_{\mu,\,\mu'\in\bar\Omega_{\pi,\Theta},\ \mu\ne\mu'}
  \bigl(\mu(\lambda) - \mu'(\lambda)\bigr).
\end{aligned}
\end{equation}
Then if $\bar r_\Theta(\lambda)\ne0$,
\[
 \sum_{i,j}S(\mathfrak g)
  \bar q_{\pi,\Theta}(F_\pi;\lambda)_{ij} + 
  \sum_{f\in I(\mathfrak g)}S(\mathfrak g)\bigl(f-f(\lambda_\Theta)\bigr)
  = \{f\in S(\mathfrak g); f|_{\Ad(G)\lambda_\Theta}=0\}.
\]
Here $I(\mathfrak g)$ is the space of the $\ad(\mathfrak g)$-invariant 
elements in the symmetric algebra $S(\mathfrak g)$ of $\mathfrak g$
and $G$ a connected complex Lie group with Lie algebra $\mathfrak g$.
\end{thm}
\begin{proof}
Let $\{v_i;\,i=1,\ldots,N\}$ be a base of $V$
such that each $v_i$ is a weight vector with weight $\varpi_i$.
Then
\[
  d\bar q_{\pi,\Theta}(F_\pi;\lambda)_{ij}|_{\lambda_\Theta} =
  \begin{cases}
  0&\text{if }\ang{\varpi_i-\varpi_j}{\lambda_\Theta}\ne0,\\
   \prod_{\mu\in\bar\Omega_{\pi,\Theta}\setminus
       \{\varpi_i|_{\mathfrak a_\Theta}\}}
   \bigl(\ang{\varpi_i}{\lambda_\Theta}-\mu(\lambda)\bigr)dF_{ij}
   &\text{if }\ang{\varpi_i-\varpi_j}{\lambda_\Theta}=0.
  \end{cases}
\]

For $\alpha\in\Sigma(\mathfrak g)\setminus\Sigma(\mathfrak g_\Theta)$
there exists a pair of weights of $\pi$ whose difference equals $\alpha$ and
therefore $\bar r_\Theta(\lambda)\ne0$ implies $\ang{\alpha}{\lambda_\Theta}\ne0$,
which assures that the centralizer of $\lambda_\Theta$ in $\mathfrak g$ equals 
$\mathfrak g_\Theta$.
Since 
\[\mathfrak g_\Theta = \sum_{
   i=j\text{ or }\varpi_i-\varpi_j\text{ is a root of }
   \mathfrak g_\Theta}\mathbb CF_{ij}
\]
and $[H,F_{ij}]=(\varpi_i-\varpi_j)(H)F_{ij}$ for $H\in\mathfrak a$,
we can prove the theorem as in the same way as in the proof of 
\cite[Theorem~4.11]{O-Cl}.
\end{proof}

\begin{rem}
There is a natural projection 
$\bar p_{\pi,\Theta}:\Omega_{\pi,\Theta}\to\bar\Omega_{\pi,\Theta}$.
We say that $\mu\in\bar\Omega_{\pi,\Theta}$ is {\it ramified in the quantization\/}
of $\bar q_{\pi,\Theta}$ to $q_{\pi,\Theta}$ if 
$\bar p_{\pi,\Theta}^{-1}(\mu)$ is not a single element.

If $\pi$ is of multiplicity free, 
then there is no ramified element in $\bar\Omega_{\pi,\Theta}$ (cf.~Proposition~\ref{prop:free}~i)).
In this case, consider $\mathfrak g$ as an abelian Lie algebra acting on $S(\mathfrak g)$
by the multiplication
and define the $\mathfrak g$-module
$M^0_\Theta(\lambda)=S(\mathfrak g)/\sum_{X\in\mathfrak p_\Theta}S(\mathfrak g)(X-\lambda_\Theta(X))$.
Then taking a ``classical limit" as in \cite{O-Cl},
we can prove $\bar q_{\pi,\Theta}(F_\pi;\lambda)M^0_\Theta(\lambda)=0$.
Moreover if $\bar r_\Theta(\lambda)\ne0$,
the polynomial $\bar q_{\pi,\Theta}(x;\lambda)$ is minimal in the obvious sense.
\end{rem}
\section{Examples}\label{sec:exm}
In this section we give the explicit form of the characteristic polynomials 
of some small dimensional representations $\pi$ of classical and exceptional Lie algebras $\mathfrak g$.
(As in the previous sections,
we always assume that $\mathfrak g$ and $\pi$ satisfy \eqref{eq:setting}.)
In some special cases we also calculate the global minimal polynomials.
Note that if $q_\pi(x)=\prod_{1\le i\le m}(x-\varpi_i-C_i)$
with suitable $\varpi_i\in\mathfrak a^*$ and $C_i\in\mathbb C$ is the characteristic 
polynomial, then the global minimal polynomial $q_{\pi,\Theta}(x,\lambda)$ for a given $\Theta$ equals 
$\prod_{i\in I}(x-\ang{\varpi_i}{\lambda_\Theta+\rho} - C_i)$ with 
a certain subset $I$ of $\{1,\ldots,m\}$.

It is clear that if the dimension of $\pi$ is small,
then the degree of $q_{\pi,\Theta}(x,\lambda)$ is small,
which means the corresponding ideal $I_{\pi,\Theta}(\lambda)$ is generated by
elements with small degrees.
In such a case, for an extremal low weight $\varpi_\alpha$ of $\pi$ with respect to $\alpha\in\Theta$,
the degree of the polynomial $r_{\alpha,\varpi_{\alpha}}(\lambda)$ defined by \eqref{eq:gapcond}
is also small and hence
the assumptions on $\lambda$ of
Theorem~\ref{thm:gap} and Corollary~\ref{cor:ann}
become very weak.

\def\lin{\hskip-.8pt\text{---}\hskip-.8pt}
\def\linlin{\lin\lin}
\def\longlin{\hskip-.8pt\text{---}\cdots\text{---}\hskip-.8pt}
\def\rarw{\hskip-3.7pt\Longrightarrow\hskip-.1pt}
\def\larw{\hskip-2.8pt\Longleftarrow\hskip-1pt}
\def\vlin{\hskip2.4pt\hbox{\vrule height15pt width.2pt}}
\def\pup#1{\phantom{\lower#1pt\hbox{|}}}
\def\pdw#1{\phantom{\raise#1pt\hbox{|}}}
\def\dcirc{\hskip-.4pt\raise2pt\hbox{${}_\circledcirc$}\hskip-.7pt}

\def\mL#1#2{\!\displaystyle\mathop{\mbox{$#1$}}_{\mbox{$#2$}}\!}
\def\lU#1{\makebox[0pt]{\quad\quad\quad\quad$\swarrow{#1}$}}
\def\lD#1{\makebox[0pt]{\quad\quad\quad\quad$\searrow{#1}$}}
\def\RA{\rangle}

\begin{lem}[bilinear form]
Let $(\ ,\ )$ be a symmetric bilinear form on $\mathfrak a^*$
and let $\mathfrak a^*=\mathfrak a^*_1\oplus \mathfrak a^*_2$ be a direct
sum of linear subspaces with $(\mathfrak a^*_1,\mathfrak a^*_2)
=\ang{\mathfrak a^*_1}{\mathfrak a^*_2}=0$.
If there exists $C\in\mathbb C\setminus\{0\}$ such that
\[
  (\mu,\mu')=C\ang\mu{\mu'}\quad(\forall\mu,\,\mu'
   \in\mathfrak a^*_1),
\]
then
\[
  C = \sum_{\varpi\in\mathcal W(\pi)}
      m_\pi(\varpi)\frac{(\alpha,\varpi)^2}{(\alpha,\alpha)}
     \quad\text{ for }\alpha\in\mathfrak a_1^*
     \text{ such that }(\alpha,\alpha)\ne0.
\]
Here $m_\pi(\varpi)$ denotes the multiplicity of the weight $\varpi\in\mathcal{W}(\pi)$.
\end{lem}
\begin{proof}
Let $H_\alpha\in\mathfrak a$ correspond to $\alpha$
by the bilinear form $\ang{\ }{\ }$.
Then we have
\begin{align*}
 C(\alpha,\alpha)&= C^2\ang{\alpha}{\alpha}
                  = C^2\trace\pi(H_\alpha)^2\\
                 &= C^2\sum_{\varpi\in\mathcal W(\pi)}m_\pi(\varpi)
                    \left(\ang{\alpha}{\varpi}\right)^2
                  =  \sum_{\varpi\in\mathcal W(\pi)}m_\pi(\varpi)
                    (\alpha,\varpi)^2.\qedhere
\end{align*}
\end{proof}

In the following examples $\varepsilon_1$, $\varepsilon_2,\ldots$ constitute
a base of a vector space with symmetric bilinear form $(\ ,\ )$
defined by $(\varepsilon_i,\varepsilon_j)=\delta_{ij}$.
We consider $\mathfrak a^*$ a subspace of this space
where $\varepsilon_1-\varepsilon_2$ etc.~ are suitable elements in 
$\Psi(\mathfrak g)$ (cf.~\cite{Bo}).

$C_\pi$ equals the constant $C$ in the above lemma for 
$\mathfrak a_1 = \mathfrak a\cap[\mathfrak g,\mathfrak g]$.
$C'_\pi$ is the similar constant in the case when $\mathfrak a_1$ is the
center of $\mathfrak g$.
Then we can calculate $\ang{\ }{\ }$ under the base 
$\{\varepsilon_1,\varepsilon_2,\ldots\}$ by the above lemma.

\begin{exmp}[$A_{n-1}$]\label{ex:A}
\[
\offinterlineskip
\vbox{
\halign{$ # $&$ # $&$ # $&$ # $\cr
 \alpha_1 & \alpha_2 &\!\!\alpha_{n-2}  & \alpha_{n-1}\pup7\cr
 \circ\linlin&\circ\longlin&\circ\linlin&\circ\cr
  \pdw{12}\cr
    }
}
\quad\quad
\vbox{
\halign{$ # $&$ # $&$ # $&$ # $\cr
 1 & 1 & 1 & 1\pup7\cr
 \circ\linlin&\circ\longlin&\circ\linlin&\circ\cr
 \hskip1pt\backslash & & &\!/\cr
 \ \linlin\!&\linlin\!\bullet\!\linlin&\linlin\cr
    }
}
\]

$\Psi=\{\alpha_1 = \ve_1-\ve_2,\ldots,\alpha_{n-1}=\ve_{n-1}-\ve_n\}$

$\rho=\sum_{\nu=1}^n(\tfrac{n-1}2 - (\nu-1))\ve_\nu
            =\sum_{\nu=1}^{n-1}\tfrac{\nu(n-\nu)}2\alpha_\nu$
\medskip

\noindent
{\rm i)}
$\mathfrak g=\mathfrak{gl}_n$

$\pi=\varpi_k:=\ve_1+\cdots+\ve_k=\bigwedge^k\varpi_1$ %
(minuscule, $k=1,\ldots,n-1$)

$\dim\varpi_k = \binom nk$

$%
(\varpi_k,\rho) = \frac{k(n-k)}2
$

$\mathcal{W}(\varpi_k)=\{\ve_{\nu_1}+\cdots +\ve_{\nu_k};\,
 1\le\nu_1<\cdots<\nu_k\le n\}$

$
  C_{\varpi_k}
   = \tfrac12\sum_{1\le\nu_1<\cdots<\nu_k\le n}
        (\ve_{\nu_1}+\cdots +\ve_{\nu_k}, %
		\ve_1-\ve_2)^2
    = \binom{n-2}{k-1}
$

$ 
 C'_{\varpi_k} =
 \tfrac1n\sum_{1\le\nu_1<\cdots<\nu_k\le n}
        (\ve_{\nu_1}+\cdots +\ve_{\nu_k}, 
		\ve_1+\cdots+\ve_n)^2
   = k\binom{n-1}{k-1}
$

$
 \ang{\varepsilon_i}{\varepsilon_j} = \frac{(n-k)!(k-1)!}{n!}
 \bigl(\frac{n-1}{n-k}(n\delta_{ij}-1)+\frac1k\bigr)
$

$q_{\varpi_k}(x) = 
\prod_{1\le i_1<\cdots<i_k\le n}\big(x - (\ve_{i_1} + \cdots + \ve_{i_k}) - 
\frac{k!(n-k)!}{2(n-2)!}\big)$
\medskip

\noindent
{\rm ii)}
$\mathfrak g=\mathfrak{gl}_n$

$V=V_m:=\{\text{homogeneous polynomials of $(x_1,\ldots,x_n)$ with degree $m$}\}$

$\pi= m\varepsilon_1$ (multiplicity free,\ $m=1,2,\cdots$)

$\mathcal W(m\ve_1) = \{m_1\ve_1 + \cdots + m_n\ve_n;
 \, m_1+\cdots+m_n = m,\ m_j\in\mathbb Z_{\ge 0}\}$

$\dim m\ve_1 = {}_nH_m = \binom{n+m-1}{m} = \frac{(n+m-1)!}{m!(n-1)!}$

$
C_{m\ve_1} = 
\frac12\sum_{m_1+\cdots+m_n=m}(m_1\ve_1 + \cdots + m_n\ve_n,\ve_1-\ve_2)^2
$

$
 = \frac12\sum_{k=0}^m\sum_{m_1=0}^k(k-2m_1)^2{}_{n-2}H_{m-k}.
$

$
 = \frac1{3!}\sum_{k=0}^mk(k+1)(k+2)\frac{(m+n-k-3)!}{(n-3)!(m-k)!}
$

$
 = \frac1{3!(n-3)!}\sum_{k=0}^mk(k+1)(k+2)(m+n-(k+3))\cdots(m+n-(k+n-1))
$

$
 = \cdots = \frac{(m+n)!}{(n+1)!(m-1)!}
$

$C'_{m\ve_1}
 = \frac1n\sum_{m_1+\cdots+m_n=m}(m_1\ve_1 + \cdots + m_n\ve_n,\ve_1+\cdots+\ve_n)^2
 = \frac{m^2}n{}_nH_m
$ 

$ = \frac{(n+m-1)!}{(m-1)!n!}m
 = \frac{m(m+1)\cdots(m+n-1)}{n!}m
$

$q_{m\ve_1}(x) =
 \prod_{\substack{m_1+\cdots + m_n = m\\ m_i\in\mathbb Z_{\ge0}}}
 \bigl(
   x - \sum_{i=1}^nm_i\ve_i -
  \frac{m(m + n - 1) - \sum_{i=1}^nm_i^2}{2C_{m\ve_1}}
 \bigr)
$

\noindent
{\rm iii)}
$\mathfrak g=\mathfrak{sl}_n$

$\pi=\varpi_1 + \varpi_{n-1} = \ve_1 - \ve_n$ (adjoint)

$\dim(\varpi_1+\varpi_n) = n^2 - 1$

$C_{\varpi_1+\varpi_{n-1}}= 2n$

$(\varpi_1+\varpi_{n-1},\rho)=n-1$

$q_{\varpi_1+\varpi_{n-1}}(x)=(x-\frac12)
\prod_{1\le i<j\le n}\bigl((x-\frac{n-1}{2n})^2
 - (\ve_i-\ve_j)^2\bigr)$
\medskip

In \cite{O-Cl} we choose 
$\Psi'=\{\alpha'_1=\varepsilon_2-\varepsilon_1,\ldots,
\alpha'_{n-1}=\varepsilon_n-\varepsilon_{n-1}\}$
as a fundamental system of $\mathfrak{gl}_n$ and then $\bar\pi=\varpi_1$ is
the lowest weight of the natural representation $\pi$ of $\mathfrak{gl}_n$.
For a strictly increasing sequence 
\begin{equation}\label{eq:seq}
 n_0=0<n_1<\cdots<n_L=n
\end{equation}
we put $n'_j=n_j - n_{j-1}$ and
$\Theta=\bigcup_{k=1}^L\bigcup_{n_{k-1}<\nu< n_k}\{\alpha'_\nu\}$ and study
the minimal polynomial $q_{\pi,\Theta}(x;\lambda)$ in \cite{O-Cl} for
$\lambda=\bigl(\lambda_k\bigr)\in\mathbb C^L\simeq\mathfrak a_\Theta^*$.
Define $\rho'=-\rho$ and put
\begin{equation}\label{eq:barlambda}
 \bar\lambda_1\varepsilon_1+\cdots+\bar\lambda_n\varepsilon_n
 =\rho'+ \sum_{k=1}^L\lambda_k\Bigl(\sum_{n_{k-1}<\nu\le n_k}\varepsilon_\nu\Bigr).
\end{equation}
The partially ordered set of the weights of $\pi$ is as follows
\[
 \varepsilon_1\xrightarrow{\ \alpha'_1\ }\varepsilon_2\xrightarrow{\ \alpha'_2\ }\cdots
 \cdots\xrightarrow{\alpha'_{n_k-1}}\varepsilon_{n_k}\xrightarrow{\ \alpha'_{n_k}\ }
 \varepsilon_{n_k+1}\xrightarrow{\alpha'_{n_k+1}}\cdots\cdots
 \xrightarrow{\ \alpha'_{n-1}\ }\varepsilon_n.
\]
Then 
$\overline{\mathcal W}_\Theta(\pi)=
\{\varepsilon_{n_0+1},\ldots,\varepsilon_{n_{L-1}+1}\}$ and 
Theorem~\ref{thm:min} says
\begin{align*}
 q_{\pi,\Theta}(x,\lambda) 
 &=\prod_{k=1}^L\bigl(x - \lambda_k - \frac12
  (\varepsilon_1 - \varepsilon_{n_{k-1}+1},
   \varepsilon_1 + \varepsilon_{n_{k-1}+1} - 2\rho') \bigr)\\
 &= \prod_{k=1}^L\bigl(x - \lambda_k - n_{k-1}\bigr)
\end{align*}
and it follows from Remark~\ref{rm:Agap} that
\[
 r_{\alpha'_i}(\lambda)
  = \prod_{\nu=k+1}^L\bigl(\bar\lambda_{i+1} - \bar\lambda_{n_{\nu-1}+1} \bigr)
     \prod_{\nu=1}^{k-1}\bigl(\bar\lambda_{i} - 
       \bar\lambda_{n_\nu} \bigr)
\]
 in Definition~\ref{def:gap} if $n_{k-1}<i<n_k$.  
This result coincides with \cite[Theorem~4.4]{O-Cl}.
Note that if $\lambda$ satisfies the condition:
\begin{equation}\label{eq:dexist}
\ang{\lambda+\rho'}{\beta}=0\text{ with }\beta\in\Sigma(\mathfrak g)\ 
\Rightarrow\ \forall\alpha'\in\Theta\ \ang{\beta}{\alpha'}=0,
\end{equation}
then $r_{\alpha'}(\lambda)\ne0$ for each $\alpha'\in\Theta$.

Let $\pi_{\varpi_k}$ be the minuscule representation $\varpi_k$ in i)
and we here adopt the fundamental system $\Psi'$ as above.
The decomposition 
\begin{equation}\label{eq:glrest}
 \pi_{\varpi_k}|_{\mathfrak g_\Theta} = \bigoplus_{
  \substack{
   k_1+\cdots+k_L=k\\
   0\le k_j\le n'_j\ (j=1,\ldots,L)
  }
 }
 \pi_{k_1,\ldots,k_L}
\end{equation}
is a direct consequence of Proposition~\ref{prop:free}~i).
Here $\pi_{k_1,\ldots,k_L}$ denotes the irreducible representation of 
$\mathfrak g_\Theta$ with lowest weight 
$\sum_{j=1}^L(\ve_{n_{j-1}+1}+\cdots+\ve_{n_{j-1}+k_j})$.
Then by Proposition~\ref{prop:minfree}~i) we have
\begin{multline*}
 q_{\pi_{\varpi_k},\Theta}(x;\lambda) = 
 \prod_{
  \substack{
   k_1+\cdots+k_L=k\\
   0\le k_j\le n'_j\ (j=1,\ldots,L)
  }
 }
 \Bigl(
  x - \sum_{i=1}^n\sum_{j=1}^L\sum_{\nu=1}^{k_j}\bar\lambda_i
     \ang{\varepsilon_i}{\varepsilon_{n_{j-1}+\nu}}
    - \frac{k!(n-k)!}{2(n-2)!}
 \Bigr)\\
 =
 \prod_{
  \substack{
   k_1+\cdots+k_L=k\\
   0\le k_j\le n'_j\ (j=1,\ldots,L)
  }
 }
 \Bigl(
  x - C''_{\varpi_k}(n-1)\sum_{j=1}^Lk_j\bigl(\lambda_j+n_{j-1}+\frac{k_j-k}2\bigr) \\
   +C''_{\varpi_k}(k-1)\sum_{j=1}^Ln'_j\bigl(\lambda_j+n_{j-1}+\frac{n'_j-n}2\bigr)
 \Bigr)
\end{multline*}
with $C''_{\varpi_k}=\frac{(n-k-1)!(k-1)!}{(n-1)!}$.
To deduce the final form we have used the relation
$\sum_{j=1}^Ln'_jn_{j-1}=\frac{n^2-\sum_{j=1}^Ln'_j}2$.
\begin{rem}
Put $\mathfrak g'_\Theta = [\mathfrak g_\Theta, \mathfrak g_\Theta]$.
Then the irreducible decomposition of $\pi_{\varpi_k}|_{\mathfrak g'_\Theta}$
is not of multiplicity free if and only if there exist an integer $K$ and subsets 
$I$ and $J$ of $\{1,\ldots,L\}$ such that
\[
  K = \sum_{i\in I}n'_i = \sum_{j\in J}n'_j\le k,\ K\le n-k\text{ and }I\ne J.
\]
This is clear from \eqref{eq:glrest} because 
$\pi_{k_1,\ldots,k_L}|_{\mathfrak g'_\Theta}=
\pi_{k'_1,\ldots,k'_L}|_{\mathfrak g'_\Theta}$ if and only if
$k_i=k'_i$ or $(k_i,k'_i)=(0,n'_i)$ or $(n'_i,0)$ for $i=1,\ldots,L$.
\end{rem}
\end{exmp}
\begin{exmp}[$B_n$]
$\mathfrak g=\mathfrak{o}_{2n+1}$
\[
\offinterlineskip
\vbox{
\halign{$ # $&$ # $&$ # $&$ # $\cr
 \alpha_1 & \alpha_2 &\!\! \alpha_{n-1}  & \alpha_n\pup7\cr
 \circ\linlin&\circ\longlin&\circ\rarw&\circ\cr
  \pdw{15}\cr
    }
}
\qquad\qquad
\vbox{
\halign{$ # $&$ # $&$ # $&$ # $\cr
  1 & 2 & 2 & 2\pdw4\cr
  \circ\linlin&\circ\longlin&\circ\rarw&\circ\quad \raise1pt\hbox{$n\ge3$}\cr
  &\vlin\cr
  &\bullet\cr
    }
}
\]

$\Psi=\{\alpha_1 = \ve_1-\ve_2,\ldots,\alpha_{n-1}=\ve_{n-1}-\ve_n,\,
\alpha_n=\ve_n\}$

$\rho=\sum_{\nu=1}^n(n - \nu + \frac12)\ve_\nu
     =\sum_{\nu=1}^n\frac{\nu(2n-\nu)}2\alpha_\nu$
\medskip

\noindent
{\rm i)}
$\pi=\varpi_1:=\ve_1$ (multiplicity free)

$\dim\varpi_1 = 2n+1$

$(\varpi_1,\rho) = n - \frac 12$

$C_{\varpi_1}=\sum(\pm\ve_\nu,\ve_1)^2+(0,\ve_1)^2 = 2$

$q_{\varpi_1}(x) = (x - \frac{n}2)
\prod_{i=1}^n\bigl((x - \frac{2n-1}4)^2
 - \ve_i^2\bigr)$
\medskip

\noindent
{\rm ii)}
$\pi=\varpi_n:=\frac12(\ve_1+\cdots+\ve_n)$ (minuscule)

$\dim\varpi_n = 2^n$

$%
(\varpi_n,\rho) = \frac{(2n-1)+(2n-3)+\cdots+1}4=\frac{n^2}4
$

$C_{\varpi_n}=\sum(\pm\ve_1\pm\cdots\pm\ve_n,\ve_1)^2
    = 2^n$

$q_{\varpi_n}(x) = \prod_{c_1=\pm1,\cdots,c_n=\pm1}
\bigl(x - \frac12(c_1\ve_1+\cdots+c_n\ve_n) - 
\frac{n^2}{2^{n+2}}\bigl)$
\medskip

\noindent
{\rm iii)}
$\pi=\varpi_2 := \ve_1+\ve_2$ (adjoint)
$\cdots\ \varpi_2$ is not a fundamental weight if $n=2$.

$\dim\varpi_2 = n(2n+1)$

$C_{\varpi_2} = 4n-2$

$(\varpi_2,\rho) = 2n - 2$

$\ve_1=\varpi_2 - \alpha_2 - \cdots - \alpha_n$

$q_{\varpi_2}(x)=(x-\frac12)
\prod_{1\le i<j\le n}\bigl((x-\frac{n-1}{2n-1})^2
 - (\ve_i-\ve_j)^2\bigr)\bigl((x-\frac{n-1}{2n-1})^2
 - (\ve_i+\ve_j)^2\bigr)
\prod_{i=1}^n\bigl((x-\frac{4n-3}{8n-4})^2 - \ve_i^2\bigr)
$
\medskip

Choose $\Psi'=\{\alpha_1'=\varepsilon_2-\varepsilon_1,\ldots,
\alpha_{n-1}'=\varepsilon_n-\varepsilon_{n-1},\,\alpha_n'=-\varepsilon_n\}$
as a fundamental system.
Then the partially ordered set of the weights of the natural representation $\pi$
of $\mathfrak o_{2n+1}$ is shown by
\begin{align*}
 &\varepsilon_1\xrightarrow{\alpha'_1}\varepsilon_2\xrightarrow{\alpha'_2}\cdots
 \cdots\xrightarrow{\alpha'_{n_k-1}}\varepsilon_{n_k}\xrightarrow{\alpha'_{n_k}}
 \varepsilon_{n_k+1}\xrightarrow{\alpha'_{n_k+1}}\cdots\cdots
 \xrightarrow{\alpha'_{n-1}}\varepsilon_n
 \xrightarrow{\alpha'_n}0\\
 &\xrightarrow{\alpha'_n}-\varepsilon_n\xrightarrow{\alpha'_{n-1}}
 \cdots\cdots
 \xrightarrow{\alpha'_{n_k+1}}-\varepsilon_{n_k+1}
 \xrightarrow{\alpha'_{n_k}}-\varepsilon_{n_k}
 \xrightarrow{\alpha'_{n_k-1}}%
 \cdots\cdots
 \xrightarrow{\alpha'_1}-\varepsilon_1.
\end{align*}
Here we use the same notation as in \eqref{eq:seq} and \eqref{eq:barlambda}.
Put
$\Theta=\bigcup_{k=1}^L\bigcup_{n_{k-1}<\nu< n_k}\{\alpha'_\nu\}$ and 
$\bar\Theta=\Theta\cup\{\alpha'_n\}$.
Then
\begin{align*}
 \overline{\mathcal W}_{\bar\Theta}(\pi)
   &=\{\varepsilon_{n_0+1},\ldots,\varepsilon_{n_{L-1}+1},
      -\varepsilon_{n_{L-1}},\ldots,-\varepsilon_{n_1}\},\\
 \overline{\mathcal W}_\Theta(\pi)
  &=\overline{\mathcal W}_{\bar\Theta}(\pi)\cup\{0,\,-\varepsilon_n\}.
\end{align*}
Hence by Theorem~\ref{thm:min}
\begin{align*}
 q_{\pi,\Theta}(x;\lambda)&=
 \Bigl(x-\frac14(\varepsilon_1,\varepsilon_1-2\rho')\Bigr)\\
 &\quad\cdot \prod_{j=1}^L\Bigl(x-\frac12\lambda_j-
   \frac14(\varepsilon_1-\varepsilon_{n_{j-1}+1},
   \varepsilon_1+\varepsilon_{n_{j-1}+1}-2\rho')\Bigr)\\
 &\quad \cdot\prod_{j=1}^L\Bigl(x+\frac12\lambda_j-
   \frac14(\varepsilon_1+\varepsilon_{n_j},
    \varepsilon_1-\varepsilon_{n_j}-2\rho')\Bigr)\\
 &=\Bigl(x-\frac n2\Bigr)
  \prod_{j=1}^L\Bigl(x-\frac{\lambda_j}2-\frac{n_{j-1}}2\Bigr)
   \Bigl(x+\frac{\lambda_j}2-\frac{2n-n_j}2\Bigr),
\\
 q_{\pi,\bar\Theta}(x;\lambda)&=
  \Bigl(x-\frac14(\varepsilon_1-\varepsilon_{n_{L-1}+1},
   \varepsilon_1+\varepsilon_{n_{L-1}+1}-2\rho')\Bigr)\\
 &\quad \cdot \prod_{j=1}^{L-1}\Bigl(x-\frac12\lambda_j-
   \frac14(\varepsilon_1-\varepsilon_{n_{j-1}+1},
   \varepsilon_1+\varepsilon_{n_{j-1}+1}-2\rho')\Bigr)\\
 &\quad\cdot\prod_{j=1}^{L-1}
    \Bigl(x+\frac12\lambda_j-\frac14(\varepsilon_1+\varepsilon_{n_j},
    \varepsilon_1-\varepsilon_{n_j}-2\rho')\Bigr)\\
 &=\Bigl(x-\frac{n_{L-1}}2 \Bigr)
   \prod_{j=1}^{L-1}\Bigl(x-\frac{\lambda_j}2-\frac{n_{j-1}}2\Bigr)
   \Bigl(x+\frac{\lambda_j}2-\frac{2n-n_j}2\Bigr).
\end{align*}
Moreover if $n_{k-1}<i<n_k$, 
\begin{align*}
 2^{2L}r_{\alpha'_i,\Theta}(\lambda)
 &= \prod_{\nu=1}^{k-1}
    \Bigl(\bar\lambda_i - \bar\lambda_{n_\nu}\Bigr)
    \prod_{\nu=k+1}^L
     \Bigl(\bar\lambda_{i+1}-\bar\lambda_{n_{\nu-1}+1}\Bigr)\\
 &\quad\cdot\Bigl(\bar\lambda_{i+1}-\frac12\Bigr)
     \prod_{\nu=1}^L\bigl(\bar\lambda_{i+1}+\bar\lambda_{n_{\nu}}\Bigr)\\
&=\frac12\prod_{\nu=1}^{k-1}
    \Bigl(\bar\lambda_i - \bar\lambda_{n_\nu}\Bigr)
    \prod_{\nu=k+1}^L
     \Bigl(\bar\lambda_{i+1}-\bar\lambda_{n_{\nu-1}+1}\Bigr)\\
 &\quad\cdot\Bigl(\bar\lambda_i+\bar\lambda_{i+1}\Bigr)
     \prod_{\nu=1}^L\bigl(\bar\lambda_{i+1}+\bar\lambda_{n_{\nu}}\Bigr),\\
 2^{2L-2}r_{\alpha'_i,\bar\Theta}(\lambda)
 &= \prod_{\nu=1}^{k-1}
      \Bigl(\bar\lambda_i - \bar\lambda_{n_\nu}\Bigr)
    \prod_{\nu=k+1}^{L}
     \Bigl(\bar\lambda_{i+1}-\bar\lambda_{n_{\nu-1}+1}\Bigr)
     \prod_{\nu=1}^{L-1}\bigl(\bar\lambda_{i+1}+\bar\lambda_{n_{\nu}}\Bigr),\\
 2^{2L-2}r_{\alpha'_n,\bar\Theta}(\lambda)
 &= \prod_{\nu=1}^{L-1} \Bigl(\bar\lambda_n - \bar\lambda_{n_\nu} \Bigr)
    \prod_{\nu=1}^{L-1} \Bigl(\frac12 + \bar\lambda_{n_\nu}\Bigr)
 = (-1)^{L-1}\prod_{\nu=1}^{L-1}\Bigl(\bar\lambda_n - \bar\lambda_{n_\nu}\Bigr)^2.
\end{align*}
Here we denote $r_\alpha(\lambda)$ corresponding to $\Theta$ and $\bar\Theta$ by 
$r_{\alpha,\Theta}(\lambda)$ and $r_{\alpha,\bar\Theta}(\lambda)$, respectively.
Note that $r_{\alpha',\bar\Theta}(\lambda)\ne0$ for $\alpha'\in\bar\Theta$
under the condition \eqref{eq:dexist} for $\bar\Theta$.
Moreover suppose $\lambda+\rho'$ is dominant.
Then $\bar\lambda_i+\bar\lambda_{i+1}=2\bar\lambda_{i+1}-1=
-2\frac{\ang{\lambda+\rho'}{-\varepsilon_{i+1}}}{\ang{-\varepsilon_{i+1}}{-\varepsilon_{i+1}}}-1\ne0$
and hence $r_{\alpha',\Theta}(\lambda)\ne0$ for $\alpha'\in\Theta$
under the condition \eqref{eq:dexist}.
\end{exmp}
\begin{exmp}[$C_n$]
$\mathfrak g=\mathfrak{sp}_n$
\[
\offinterlineskip
\vbox{
\halign{$ # $&$ # $&$ # $&$ # $\cr
 \alpha_1 & \alpha_2 & \!\!\alpha_{n-1}  & \alpha_n\pup8\cr
 \circ\linlin&\circ\longlin&\circ\larw&\circ\cr
    }
}
\qquad\qquad
\vbox{
\halign{$ # $&$ # $&$ # $&$ # $&$ # $\cr
  & 2 & 2 & 2 & 1\pup5\cr
  \bullet\rarw&\circ\linlin&\circ\longlin&\circ\larw&\circ\quad n\ge2\cr
    }
}\]

$\Psi=\{\alpha_1 = \ve_1-\ve_2,\ldots,\alpha_{n-1}=\ve_{n-1}-\ve_n,\,
\alpha_n=2\ve_n\}$

$\rho=\sum_{\nu=1}^n(n - \nu + 1)\ve_\nu
     =\sum_{\nu=1}^{n-1}\frac{\nu(2n-\nu +1)}2\alpha_\nu
     +\frac{n(n+1)}4\alpha_n$
\medskip

\noindent
i)
$\pi=\varpi_1:=\ve_1$  (minuscule)

$\dim\varpi_1 = 2n$

$C_{\varpi_1} = \sum(\pm\ve_\nu,\ve_1)^2 = 2$

$%
 (\varpi_1,\rho) = n$

$q_{\varpi_1}(x) = \prod_{i=1}^n\bigl((x - \frac{n}2)^2
 - \ve_i^2)$
\medskip

\noindent
\rm{ii)}
$\pi=2\varpi_1=2\ve_1$ (adjoint)

$\dim2\varpi_1 = n(2n+1)$

$C_{2\varpi_1} = 4(n+1)$

$(2\varpi_1,\rho) = 2n$

$q_{2\varpi_1}(x)=(x-\frac12)
\prod_{i=1}^n\bigl((x-\frac{n}{2n+2})^2
 - 2\ve_i^2\bigr)
\prod_{1\le i<j\le n}\bigl((x-\frac{2n+1}{4n+4})^2
 - (\ve_i-\ve_j)^2\bigr)\bigl((x-\frac{2n+1}{4n+4})^2
 - (\ve_i+\ve_j)^2\bigr)
$
\medskip

Choose $\Psi'=\{\alpha_1'=\varepsilon_2-\varepsilon_1,\ldots,
\alpha_{n-1}'=\varepsilon_n-\varepsilon_{n-1},\,\alpha_n'=-2\varepsilon_n\}$
as a fundamental system.
The partially ordered set of the weights of the natural representation $\pi$
of $\mathfrak{sp}_n$ is shown by
\begin{align*}
 &\varepsilon_1\xrightarrow{\alpha'_1}\varepsilon_2\xrightarrow{\alpha'_2}\cdots
 \cdots\xrightarrow{\alpha'_{n_k-1}}\varepsilon_{n_k}\xrightarrow{\alpha'_{n_k}}
 \varepsilon_{n_k+1}\xrightarrow{\alpha'_{n_k+1}}\cdots\cdots
 \xrightarrow{\alpha'_{n-1}}\varepsilon_n\\
 &\xrightarrow{\alpha'_n}-\varepsilon_n\xrightarrow{\alpha'_{n-1}}
 \cdots\cdots
 \xrightarrow{\alpha'_{n_k+1}}-\varepsilon_{n_k+1}
 \xrightarrow{\alpha'_{n_k}}-\varepsilon_{n_k}
 \xrightarrow{\alpha'_{n_k-1}}%
 \cdots\cdots
 \xrightarrow{\alpha'_1}-\varepsilon_1.
\end{align*}
Under the same notation as in the previous example, we have
\begin{align*}
 \overline{\mathcal W}_{\bar\Theta}(\pi) &=
 \{\varepsilon_{n_0+1},\ldots,\varepsilon_{n_{L-1}+1},
  -\varepsilon_{n_{L-1}},\ldots,-\varepsilon_{n_1}\},\\
 \overline{\mathcal W}_{\Theta}(\pi) &=
 \overline{\mathcal W}_{\bar\Theta}(\pi)\cup\{-\varepsilon_n\}.
\end{align*}
If $n_{k-1}<i<n_k$, it follows from Theorem~\ref{thm:min},
and Remark~\ref{rm:Agap} that
\begin{align*}
 q_{\pi,\bar\Theta}(x;\lambda) &=
 \prod_{j=1}^{L}\Bigl(
   x - \frac{\lambda_{j}}2 - \frac{n_{j-1}}2
  \Bigr)
 \prod_{j=1}^{L-1}\Bigl(
   x + \frac{\lambda_{j}}2 - \frac{2n-n_{j}+1}2
  \Bigr),\\
 q_{\pi,\Theta}(x;\lambda) &=
 \prod_{j=1}^L\Bigl(
   x - \frac{\lambda_{j}}2 - \frac{n_{j-1}}2
  \Bigr)\Bigl(
   x + \frac{\lambda_{j}}2 - \frac{2n-n_{j}+1}2
  \Bigr),\\
 2^{2L-1}r_{\alpha'_i,\Theta}(\lambda) &=
  \prod_{\nu=1}^{k-1}
   \Bigl(\bar\lambda_i - \bar\lambda_{n_\nu}\Bigr)
  \prod_{\nu=k+1}^L
   \Bigl(\bar\lambda_{i+1} - \bar\lambda_{n_{\nu-1}+1}\Bigr)
  \prod_{\nu=1}^L
   \Bigl(\bar\lambda_{i+1} + \bar\lambda_{n_\nu}\Bigr),
   \\
 2^{2L-2}r_{\alpha'_i,\bar\Theta}(\lambda) &=
  \prod_{\nu=1}^{k-1}
   \Bigl(\bar\lambda_i - \bar\lambda_{n_\nu}\Bigr)
  \prod_{\nu=k+1}^L
   \Bigl(\bar\lambda_{i+1} - \bar\lambda_{n_{\nu-1}+1}\Bigr)
  \prod_{\nu=1}^{L-1}
   \Bigl(\bar\lambda_{i+1} + \bar\lambda_{n_\nu}\Bigr),\\
 2^{2L-2}r_{\alpha'_n,\bar\Theta}(\lambda) &=
   \prod_{\nu=1}^{L-1} \bar\lambda_{n_\nu}
   \prod_{\nu=1}^{L-1}\Bigl(\bar\lambda_n - \bar\lambda_{n_\nu}\Bigr).
\end{align*}
If the condition \eqref{eq:dexist} holds,
then we have $r_{\alpha',\Theta}(\lambda)\ne0$
and $r_{\alpha',\bar\Theta}(\lambda)\ne0$
for $\alpha'\in\Theta$.
Moreover suppose $\ang{\lambda}{\alpha'_n}=0$ and $\lambda+\rho'$ is dominant.
In this case $\bar\lambda_n=-1$ and 
$\bar\lambda_{n_\nu}=
-2\frac{\ang{\lambda+\rho'}{\varepsilon_n-\varepsilon_{n_\nu}}}{\ang{\varepsilon_n-\varepsilon_{n_\nu}}{\varepsilon_n-\varepsilon_{n_\nu}}}-1\ne0$.
Hence $r_{\alpha'_n,\bar\Theta}(\lambda)\ne0$
under the condition \eqref{eq:dexist} for $\bar\Theta$.
\end{exmp}
\begin{exmp}[$D_n$]\label{ex:Dn}
$\mathfrak g=\mathfrak{o}_{2n}$
\[
\offinterlineskip
\vbox{
\halign{$ # $&$ # $&$ # $&$ # $\cr
 \alpha_1    & \alpha_2    &\!\!\alpha_{n-2}& \alpha_{n-1}\pup7\cr
 \circ\linlin&\circ\longlin&\circ\linlin&\circ\cr
             &             &\vlin\cr
             &             &\circ\cr
             &             &\alpha_n\pdw7\cr
   }
}
\quad\quad
\vbox{
\halign{$ # $&$ # $&$ # $&$ # $\cr
 1 & 2 & 2 & 1\pup4\cr
 \circ\linlin&\circ\longlin&\circ\linlin&\circ\quad \raise1pt\hbox{$n\ge4$}\cr
             &\vlin        &\vlin\cr
             &\bullet      &\circ\cr
             &             &1\pdw7\cr
   }
}
\]

$\Psi=\{\alpha_1 = \ve_1-\ve_2,\ldots,\,\alpha_{n-1}=\ve_{n-1}-\ve_n,\,
\alpha_n=\ve_{n-1}+\ve_n\}$

$\rho=\sum_{\nu=1}^n(n - \nu)\ve_\nu
   =\sum_{\nu=1}^{n-2}\frac{\nu(2n-\nu-1)}2\alpha_\nu
   +\frac{n(n-1)}4(\alpha_{n-1}+\alpha_n)$
\medskip

\noindent
{\rm i)}
$\pi=\varpi_1:=\ve_1$ (minuscule)

$\dim\varpi_1=2n$

$C_{\varpi_1}=\sum(\pm\ve_\nu,\ve_1)^2 = 2$

$%
(\varpi_1,\rho)=n - 1$

$q_{\varpi_1}(x) = 
\prod_{i=1}^n\bigl((x - \frac{n-1}2)^2
 - \ve_i^2\bigr)$
\medskip

\noindent
{\rm ii)}
$\pi=\begin{cases}
\varpi_{n-1}:=\frac12(\ve_1+\cdots+\ve_{n-1} - \ve_n) & \text{ (minuscule)} \\
\varpi_{n}:=\frac12(\ve_1+\cdots+\ve_{n-1} + \ve_n) & \text{ (minuscule)} \\
\end{cases}$

$\dim\varpi_{n-1} = \dim\varpi_n = 2^{n-1}$ 

$C_{\varpi_{n-1}}=C_{\varpi_n}=
\sum(\pm\ve_1\pm\cdots\pm\ve_n,\ve_1)^2=2^{n-1}$

$%
(\varpi_{n-1},\rho) = (\varpi_n,\rho) = \frac{n(n-1)}4$.

$q_{\varpi_{n-1}}(x) = 
\prod_{\substack{c_1=\pm1,\ldots,c_n=\pm1\\ c_1\cdots c_n=-1}}
\bigl(x - \frac12(c_1\ve_1+\cdots+ c_n\ve_n) - \frac{n(n-1)}{2^{n+1}}
\bigr)$

$q_{\varpi_n}(x) = 
\prod_{\substack{c_1=\pm1,\ldots,c_n=\pm1\\ c_1\cdots c_n=1}}
\bigl(x - \frac12(c_1\ve_1+\cdots+ c_n\ve_n) - \frac{n(n-1)}{2^{n+1}}
\bigr)$
\medskip

\noindent
{\rm iii)}
$\pi=\varpi_2:=\ve_1 + \ve_2$ (adjoint)

$\dim\varpi_2 = n(2n-1)$

$C_{\ve_1 + \ve_2} = 4(n-1)$

$(\varpi_2,\rho) = 2n-3$

$q_{\varpi_2}(x) = (x-\frac12)\prod_{1\le i<j\le n}
\bigl((x - \frac{2n-3}{4n-4})^2 - (\ve_i - \ve_j)^2\bigr)
\bigl((x - \frac{2n-3}{4n-4})^2 - (\ve_i + \ve_j)^2\bigr)
$
\medskip

\noindent
Note that the coefficient of
$\varepsilon_1\varepsilon_2\cdots\varepsilon_n$ in the polynomial 
$\sum_{\substack{c_1=\pm1,\ldots,c_n=\pm1\\ c_1\cdots c_n=1}}(c_1\ve_1+\cdots+ c_n\ve_n)^n$ of $(\varepsilon_1,\ldots,\varepsilon_n)$
does not vanish.
Hence 
\begin{equation}
Z(\mathfrak g) = \mathbb C[\trace F_{\varpi_1}^2,\trace F_{\varpi_1}^4,
\ldots,\trace F_{\varpi_1}^{2(n-1)},
\trace F_{\varpi_n}^n].
\end{equation}
Choose $\Psi'=\{\alpha_1'=\varepsilon_2-\varepsilon_1,\ldots,
\alpha_{n-1}'=\varepsilon_n-\varepsilon_{n-1},\,
\alpha_n'=-\varepsilon_n-\varepsilon_{n-1}\}$
as a fundamental system.
Then the partially ordered set of the weights of the natural representation $\pi$
of $\mathfrak o_{2n}$ is shown by
\begin{align*}
  \varepsilon_1\xrightarrow{\alpha_1'}\varepsilon_2\xrightarrow{\alpha_2'}
  \cdots\cdots
  \xrightarrow{\alpha_{n-2}'}&\,\varepsilon_{n-1}\xrightarrow{\alpha_{n-1}'}\varepsilon_n
  \\[-6pt]
  &\downarrow{{}^{\alpha_n'}}
   \phantom{\xrightarrow{\alpha_{n-1}'}\ }
   \downarrow{{}^{\alpha_n'}}\\[-6pt]
  &-\varepsilon_n\xrightarrow{\alpha_{n-1}'}
    -\varepsilon_{n-1}\xrightarrow{\alpha_{n-2}'}\cdots\cdots
    \xrightarrow{\alpha_2'}-\varepsilon_2\xrightarrow{\alpha_1'}-\varepsilon_1.
\end{align*}
Use the notation as in \eqref{eq:seq} and \eqref{eq:barlambda}.
Put
$\Theta=\bigcup_{k=1}^L\bigcup_{n_{k-1}<\nu< n_k}\{\alpha'_\nu\}$.
If $\alpha'_{n-1}\in\Theta$, we also put 
$\bar\Theta=\Theta\cup\{\alpha'_n\}$.

Then
\[
 \pi|_{\mathfrak g_\Theta} = 
  \bigoplus_{j=0}^{L-1}\pi_{\varepsilon_{n_j+1}}
  \oplus\bigoplus_{j=1}^L\pi_{-\varepsilon_{n_j}},\qquad
 \pi|_{\mathfrak g_{\bar\Theta}} =
  \bigoplus_{j=0}^{L-1}\pi_{\varepsilon_{n_j+1}}
  \oplus\bigoplus_{j=1}^{L-1}\pi_{-\varepsilon_{n_j}}.
\]
Here $\pi_\ve$ denotes the irreducible representation of $\mathfrak{g}_\Theta$
or $\mathfrak{g}_{\bar\Theta}$ with lowest weight $\ve$.
Hence if $n_{k-1}<i<n_k$, 
\begin{align*}
 q_{\pi,\Theta}(x;\lambda) &=
 \prod_{j=1}^L\left(
  x - \frac{\lambda_j}2 - \frac{n_{j-1}}2
 \right)\left(
  x + \frac{\lambda_j}2 - \frac{2n-n_j-1}2
 \right),
 \\
 q_{\pi,\bar\Theta}(x;\lambda) &=
 \prod_{j=1}^L\left(
  x - \frac{\lambda_j}2 - \frac{n_{j-1}}2
 \right)
 \prod_{j=1}^{L-1}\left(
  x + \frac{\lambda_j}2 - \frac{2n-n_j-1}2
 \right),\\
 2^{2L-1}r_{\alpha'_i,\Theta}(\lambda) &=
  \prod_{\nu=1}^{k-1}
   \Bigl(\bar\lambda_i - \bar\lambda_{n_\nu}\Bigr)
  \prod_{\nu=k+1}^L
   \Bigl(\bar\lambda_{i+1} - \bar\lambda_{n_{\nu-1}+1}\Bigr)
  \prod_{\nu=1}^{L}
   \Bigl(\bar\lambda_{i+1} + \bar\lambda_{n_\nu}\Bigr),
 \\
 2^{2L-2}r_{\alpha'_i,\bar\Theta}(\lambda) &=
  \prod_{\nu=1}^{k-1}
   \Bigl(\bar\lambda_i - \bar\lambda_{n_\nu}\Bigr)
  \prod_{\nu=k+1}^L
   \Bigl(\bar\lambda_{i+1} - \bar\lambda_{n_{\nu-1}+1}\Bigr)
  \prod_{\nu=1}^{L-1}
   \Bigl(\bar\lambda_{i+1} + \bar\lambda_{n_\nu}\Bigr),
 \\
 2^{2L-2}r_{\alpha'_n,\bar\Theta}(\lambda) &=
  (-1)^{L-1}
  \prod_{\nu=1}^{L-1}
   \Bigl(\bar\lambda_n - \bar\lambda_{n_\nu}\Bigr)
   \Bigl(\bar\lambda_{n-1} - \bar\lambda_{n_\nu}\Bigr).
\end{align*}
If $\ang{\lambda}{\alpha'_n}=0$ and $i+1=n_k<n$,
then $\bar\lambda_{i+1}+\bar\lambda_{n_k}=2(\bar\lambda_{i+1}+\bar\lambda_n)$.
Hence $r_{\alpha',\bar\Theta}(\lambda)\ne0$ for $\alpha'\in\bar\Theta$
under the condition \eqref{eq:dexist} for $\bar\Theta$.

Now suppose $\alpha'_{n-1}\notin\Theta$.
Then $n_{L-1}=n-1$.
If $\lambda_L=0$, then $q'_{\pi,\Theta}(F_\pi;\lambda)M_\Theta(\lambda)=0$
by Corollary~\ref{cor:involution} with
\[
  q'_{\pi,\Theta}(x;\lambda) =
  \left(
   x - \frac{\lambda_L}2 - \frac{n-1}2
  \right)
  \prod_{j=1}^{L-1}\left(
   x - \frac{\lambda_j}2 - \frac{n_{j-1}}2
  \right)\left(
   x + \frac{\lambda_j}2 - \frac{2n-n_j-1}2
  \right).
\]
The analogue of $r_{\alpha'_i,\Theta}(\lambda)$ in this case is
\[
  r'_{\alpha'_i,\Theta}(\lambda) =
  2^{2-2L}\bar\lambda_{i+1}\prod_{\nu=1}^{k-1}
   \Bigl(\bar\lambda_i - \bar\lambda_{n_\nu}\Bigr)
  \prod_{\nu=k+1}^{L-1}
   \Bigl(\bar\lambda_{i+1} - \bar\lambda_{n_{\nu-1}+1}\Bigr)
  \prod_{\nu=1}^{L-1}
   \Bigl(\bar\lambda_{i+1} + \bar\lambda_{n_\nu}\Bigr).
\]
If $i+1=n_k$ then $\bar\lambda_{i+1}+\bar\lambda_{n_k}=2(\bar\lambda_{i+1}+\bar\lambda_n)$.
Hence $r'_{\alpha',\bar\Theta}(\lambda)\ne0$ for $\alpha'\in\Theta$
under the condition \eqref{eq:dexist}.

Let $\pi_{\varpi_{n-1}}$ be the half spin representation $\varpi_{n-1}$ in ii)
and we here use the fundamental system $\Psi'$ defined above.
\[
 \pi_{\varpi_{n-1}}|_{\mathfrak g_\Theta} = 
 \bigoplus_{(k_1,\ldots,k_L)\in\mathbf K_\Theta}
 \pi_{k_1,\ldots,k_L},\quad
 \pi_{\varpi_{n-1}}|_{\mathfrak g_{\bar\Theta}} = 
 \bigoplus_{(k_1,\ldots,k_L)\in\mathbf K_{\bar\Theta}}
 \pi_{k_1,\ldots,k_L},
\]
where
\begin{align*}
 \mathbf K_\Theta &= \{(k_1,\ldots,k_L)\in\mathbb Z^L;\,
  0 \le k_j\le n'_j\ (j=1,\ldots,L),\\
  &\qquad\qquad\qquad\qquad
  n-k_1-\cdots-k_L\equiv1\mod 2\},\\
 \mathbf K_{\bar\Theta} &= 
  \{(k_1,\ldots,k_L)\in \mathbf K_\Theta;\,
   k_L \ge n'_L-1\}\qquad\text{(Note $\alpha'_{n-1}\in\Theta$ and $n'_L>1$)}
\end{align*}
and $\pi_{k_1,\ldots,k_L}$ is the irreducible representation of $\mathfrak g_\Theta$
or $\mathfrak g_{\bar\Theta}$ with lowest weight
\[
  \sum_{j=1}^L\frac12(\varepsilon_{n_{j-1}+1}+\cdots + \varepsilon_{n_{j-1}+k_j}
  - \varepsilon_{n_{j-1}+k_j+1}-\cdots - \varepsilon_{n_j}).
\]
Then for $\Theta'=\Theta$ or $\bar\Theta$
\begin{multline*}
 q_{\pi_{\varpi_{n-1}},\Theta'}(x;\lambda)
 = \prod_{(k_1,\ldots,k_L)\in\mathbf K_{\Theta'}}
 \Bigl(
    x - \frac{n(n-1)}{2^{n+1}}\\
   {} - \frac1{2^{n}}\sum_{j=1}^L(\bar\lambda_{n_{j-1}+1}+\cdots +
   \bar\lambda_{n_{j-1}+k_j}
   - \bar\lambda_{n_{j-1}+k_j+1}-\cdots - \bar\lambda_{n_j})
 \Bigr).
\end{multline*}
If $n'_L>1$, then
\begin{align*}
 r_{\alpha'_{n-1},\Theta'}(\lambda)
 &= \prod_{
   \substack{(k_1,\ldots,k_L)\in\mathbf K_{\Theta'}\\
    (k_1,\ldots,k_L)\ne(n'_1,\ldots,n'_{L-1}, n'_L-1)}
  }2^{1-n}\\
  &\quad\cdot\Bigl(
   \sum_{j=1}^{L}(\bar\lambda_{n_{j-1}+k_j+1}+\cdots+
    \bar\lambda_{n_j-1} + \bar\lambda_{n_j})
    -\bar\lambda_{n-1}\Bigr).
\end{align*}
\end{exmp}
\begin{exmp}[$E_6$]
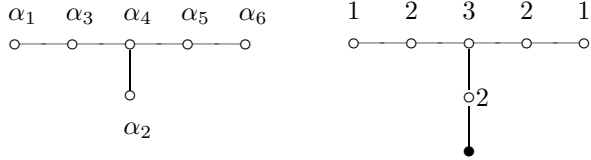
\label{ex:E6}
\[
\offinterlineskip
\vbox{
\halign{$ # $&$ # $&$ # $&$ # $&$ # $\cr
 \alpha_1    & \alpha_3   &\alpha_4 &\alpha_5 &\alpha_6\pup7\cr
 \circ\linlin&\circ\linlin&\circ\linlin&\circ\linlin&\circ\cr
			 &            &\vlin\cr
             &            &\circ\cr
             &            &\alpha_2\pdw8\cr
 \pdw{3}\cr
   }
}
\quad\quad
\vbox{
\halign{$ # $&$ # $&$ # $&$ # $&$ # $\cr
 1 & 2 & 3 & 2 & 1\pup7\cr
 \circ\linlin&\circ\linlin&\circ\linlin&\circ\linlin&\circ\cr
			 &            &\vlin\cr
             &            &\circ\lower1pt\hbox{2}\cr
             &            &\vlin\cr
		     &            &\bullet\cr
   }

}
\]

$\Psi=\{\alpha_1 = \frac12(\ve_1+\ve_8)-
\frac12(\ve_2+\ve_3+\ve_4+\ve_5+\ve_6+\ve_7),\,
\alpha_2 = \ve_1 + \ve_2, \alpha_3 = \ve_2 - \ve_1,\,
\alpha_4 = \ve_3 - \ve_2,\,
\alpha_5 = \ve_4 - \ve_3,\,
\alpha_6 = \ve_5 - \ve_4\}$

$\rho=\ve_2 + 2\ve_3 + 3\ve_4+ 4\ve_5 + 4(\ve_8-\ve_7-\ve_6)
     = 8\alpha_1 + 11\alpha_2 + 15\alpha_3 + 21\alpha_4 + 15\alpha_5
       + 8\alpha_6$
\medskip

\noindent
{\rm i)}
$\pi=\begin{cases}
\varpi_1 := \frac23(\ve_8 - \ve_7 -\ve_6) & \text{ (minuscule)}\\
\varpi_6 := \frac13(\ve_8 - \ve_7 - \ve_6) + \ve_5 & \text{ (minuscule)}
\end{cases}$

$\dim\varpi_1 = \dim\varpi_6 = 27$

$C_{\varpi_1} = C_{\varpi_6} = 6
\qquad\text{(see below)}$

$(\varpi_1, \rho) = (\varpi_6,\rho) = 8$

$q_{\varpi_i}(x) = 
\prod_{\varpi\in W_{E_6}\varpi_i}\bigl(x - \varpi - \frac43\bigr)$
 for $i=1$ and $6$.
\medskip

\noindent
{\rm ii)}
$\pi=\varpi_2 := \frac12(\ve_1+\ve_2+\ve_3+\ve_4+\ve_5-\ve_6-\ve_7+\ve_8)$
(adjoint)

$\dim\varpi_2 = 78$

$C_{\varpi_2} = 24$

$(\varpi_2,\rho) = 11$

$q_{\varpi_2} = (x-\frac12)
\prod_{\alpha\in\Sigma(E_6)}\bigl(x - \alpha - \frac{11}{24}
\bigr)$
\medskip

Expressing a weight by the linear combination of the fundamental weights $\varpi_j$,
we indicate the weight by the symbol arranging the coefficients in the corresponding 
position of the Dynkin diagram. 
For example, $\varpi=\sum_{j=1}^6m_j\varpi_j$
is indicated by the symbol $m_1m_3\mL{m_4}{m_2}m_5m_6$.
Moreover for a positive integer $m$ we will sometimes write $\bar m$ in place of $-m$.

Let $\pi$ be the minuscule representation $\varpi_1$ in i).
Then the partially ordered set of the weights of $\pi$ is
shown by the following.
Here the number $j$ beside an arrow represents $-\alpha_j$.
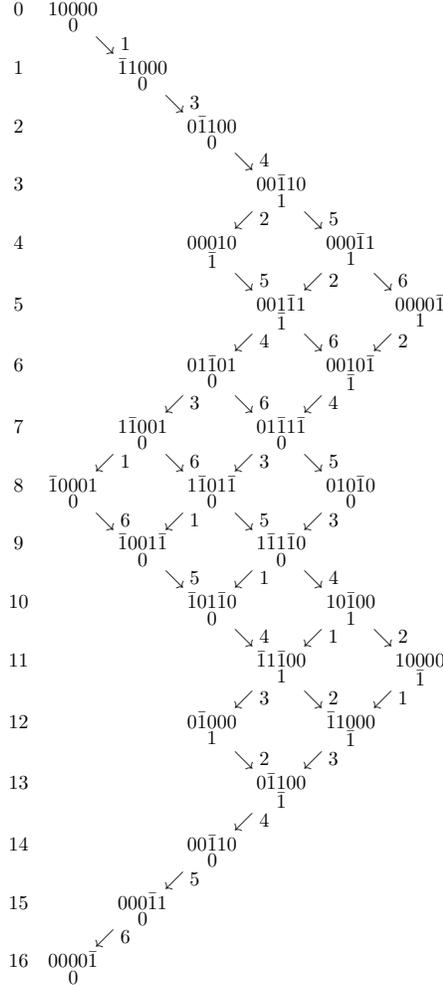
\[
\scalebox{0.75}{$
 \begin{matrix}
0&  10\mL0000 &           &           &           &           &           \\
 &  \lD1      &           &           &           &           &           \\
1&            & {\bar1}1\mL0000 &     &           &           &           \\
 &            & \lD3      &           &           &           &           \\
2&            &           & 0{\bar1}\mL1000 &     &           &           \\
 &            &           & \lD4      &           &           &           \\
3&            &           &           & 00\mL{\bar1}110 &     &           \\
 &            &           & \lU2      & \lD5      &           &           \\
4&            &           & 00\mL0{\bar1}10 &     & 00\mL01{\bar1}1 &     \\
 &            &           & \lD5      & \lU2      & \lD6      &           \\
5&            &           &           & 00\mL1{\bar1}{\bar1}1 & & 00\mL010{\bar1} \\
 &            &           & \lU4      & \lD6      & \lU2      &           \\
6&            &           & 01\mL{\bar1}001 &     & 00\mL1{\bar1}0{\bar1} &\\
 &            &  \lU3     & \lD6      & \lU4      &           &           \\
7&            & 1{\bar1}\mL0001 &     & 01\mL{\bar1}01{\bar1} &&          \\
 &  \lU1      & \lD6      & \lU3      & \lD5      &           &           \\
8&  {\bar1}0\mL0001 &     & 1{\bar1}\mL001{\bar1} &     & 01\mL00{\bar1}0 &\\
 &  \lD6      & \lU1      & \lD5      & \lU3      &           &           \\
9&      & {\bar1}0\mL001{\bar1} &     & 1{\bar1}\mL10{\bar1}0 &&          \\
 &            &  \lD5     & \lU1      &  \lD4     &           &           \\
10&           &           & {\bar1}0\mL10{\bar1}0 && 10\mL{\bar1}100 &    \\
 &            &           & \lD4      &  \lU1     & \lD2      &           \\
11&           &           &           & {\bar1}1\mL{\bar1}100 && 10\mL0{\bar1}00 \\
 &            &           & \lU3      & \lD2      & \lU1      &           \\
12&           &           & 0{\bar1}\mL0100 && {\bar1}1\mL0{\bar1}00 &    \\
 &            &           & \lD2      & \lU3      &           &           \\
13&           &           &           & 0{\bar1}\mL1{\bar1}00 &&          \\
 &            &           & \lU4      &           &           &           \\
14&           &           & 00\mL{\bar1}010 &     &           &           \\
 &            & \lU5      &           &           &           &           \\
15&           & 00\mL00{\bar1}1 &     &           &           &           \\
 &  \lU6      &           &           &           &           &           \\
16& 00\mL000{\bar1} &     &           &           &           &           
 \end{matrix}
$}
\]

The type $A_5$ corresponding to $\{\alpha_1,\alpha_3,\ldots,\alpha_6\}$
is contained in type $E_6$.
The highest weights of the restriction $(E_6,\pi)|_{A_5}$ are
$\varpi_1=10\mL0000$, $\varpi_5 - \varpi_2
=w_2w_4w_3w_1\varpi_1=00\mL0{\bar1}10$ and $\varpi_1-\varpi_2
=w_2w_4w_5w_6w_3w_4w_5(\varpi_5 - \varpi_2)=10\mL0{\bar1}00$.
Here we put $w_j=w_{\alpha_j}$.
Hence
$(E_6,\pi)|_{A_5} = 2(A_5,\varpi_1) + (A_5,\varpi_4)$
and
$C_{\varpi_1} = C_{\varpi_6} =
2\binom{5-1}{1-1} + \binom{5-1}{2-1} = 6$.

Now use the fundamental system
$\Psi'=\{\alpha'_1=-\alpha_1,\ldots,\alpha'_6=-\alpha_6\}$.
Then the lowest weight $\bar\pi$ of $\pi$ equals $\varpi_1$. 
Putting $\Theta_i = \Psi'\setminus\{\alpha'_i\}$,
we have
\begin{align*}
 \overline{\mathcal W}_{\Theta_1}(\pi) &=
  \bigl\{10\mL0000,\,\bar11\mL0000,\,\bar10\mL0001\bigr\},\\
 \overline{\mathcal W}_{\Theta_2}(\pi)  &=
  \bigl\{10\mL0000,\, 00\mL0{\bar1}10,\, 10\mL0{\bar1}00\bigr\},\\
 \overline{\mathcal W}_{\Theta_3}(\pi)  &=
  \bigl\{10\mL0000,\, 0\bar1\mL1000,\, 1\bar1\mL0001,\, 0\bar1\mL0100\bigr\},\\
 \overline{\mathcal W}_{\Theta_4}(\pi)  &=
  \bigl\{10\mL0000,\, 00\mL{\bar1}110,\, 01\mL{\bar1}001,\, 10\mL{\bar1}100,\,
  00\mL{\bar1}010\bigr\},\\
 \overline{\mathcal W}_{\Theta_5}(\pi)  &=
  \bigl\{10\mL0000,\, 00\mL01{\bar1}1,\, 01\mL00{\bar1}0,\, 00\mL00{\bar1}1\bigr\},\\
 \overline{\mathcal W}_{\Theta_6}(\pi)  &=
  \bigl\{10\mL0000,\, 00\mL010{\bar1},\, 00\mL000{\bar1}\bigr\}.
\end{align*}
If we identify $\mathfrak{a}_{\Theta_i}^*$ with $\mathbb C$
by $\lambda_{\Theta_i}=\lambda\varpi_i$
and put $\bar\pi - \Lambda = \sum_{j}m^j_\Lambda\alpha_j$ 
for $\Lambda\in\mathcal{W}(\pi)$,
then Proposition~\ref{prop:minfree}~i) implies
\begin{equation}\label{eq:minmin}
 q_{\pi,\Theta_i}(x;\lambda) = 
 \prod_{\Lambda\in\overline{\mathcal W}_{\Theta_i}(\pi)}
 \Bigl(x - (\ang{\bar\pi}{\varpi_i} - m^i_\Lambda\ang{\alpha_i}{\varpi_i})\lambda
  - \sum_{j} m^j_\Lambda\ang{\alpha_j}{\rho}\Bigr).
\end{equation}
Since $\ang{\alpha_j}{\varpi_j}=\ang{\alpha_j}{\rho}=
\frac12\ang{\alpha_j}{\alpha_j}=\frac16$ and
\begin{align*}
\ang{\varpi_1}{\varpi_1}&=\frac29,
&\ang{\varpi_1}{\varpi_2}&=\frac16,
&\ang{\varpi_1}{\varpi_3}&=\frac5{18},\\
\ang{\varpi_1}{\varpi_4}&=\frac13,
&\ang{\varpi_1}{\varpi_5}&=\frac29,
&\ang{\varpi_1}{\varpi_6}&=\frac19,
\end{align*}
we get
\begin{align*}
q_{\pi,\Theta_1}(x;\lambda) &= 
 \Bigl(x - \frac29\lambda\Bigr)
 \Bigl(x - \frac1{18}\lambda - \frac16\Bigr)
 \Bigl(x + \frac19\lambda - \frac43\Bigr),\\
q_{\pi,\Theta_2}(x;\lambda) &= 
 \Bigl(x - \frac16\lambda\Bigr)
 \Bigl(x - \frac23\Bigr)
 \Bigl(x + \frac16\lambda - \frac{11}6\Bigr),\\
q_{\pi,\Theta_3}(x;\lambda) &= 
 \Bigl(x - \frac{18}5\lambda\Bigr)
 \Bigl(x - \frac19\lambda - \frac13\Bigr)
 \Bigl(x + \frac1{18}\lambda - \frac76\Bigr)
 \Bigl(x + \frac29\lambda - 2\Bigr),\\
q_{\pi,\Theta_4}(x;\lambda) &= 
 \Bigl(x - \frac13\lambda\Bigr)
 \Bigl(x - \frac16\lambda - \frac12\Bigr)
 \Bigl(x - 1\Bigr)
 \Bigl(x + \frac16\lambda - \frac53\Bigr)
 \Bigl(x + \frac13\lambda - \frac73\Bigr),\\
q_{\pi,\Theta_5}(x;\lambda) &= 
 \Bigl(x - \frac29\lambda\Bigr)
 \Bigl(x - \frac1{18}\lambda - \frac23\Bigr)
 \Bigl(x + \frac19\lambda - \frac43\Bigr)
 \Bigl(x + \frac5{18}\lambda - \frac52\Bigr),\\
q_{\pi,\Theta_6}(x;\lambda) &= 
 \Bigl(x - \frac19\lambda\Bigr)
 \Bigl(x + \frac1{18}\lambda - \frac56\Bigr)
 \Bigl(x + \frac29\lambda - \frac83\Bigr). 
\end{align*}
\end{exmp}
\begin{exmp}[$E_7$]\label{ex:E7}
\[
\offinterlineskip
\vbox{
\halign{$ # $&$ # $&$ # $&$ # $&$ # $&$ # $\cr
 \alpha_1    & \alpha_3   &\alpha_4 &\alpha_5& \alpha_6&\alpha_7\pup7\cr
 \circ\linlin&\circ\linlin&\circ\linlin&\circ\linlin&\circ\linlin&\circ\cr
             &            &\vlin\cr
             &            &\circ\cr
             &            &\alpha_2\pdw8\cr
   }
}
\quad\quad
\vbox{
\halign{$ # $&$ # $&$ # $&$ # $&$ # $&$ # $&$ # $\cr
 & 2 & 3 & 4 & 3 & 2 & 1\pup7\cr
 \bullet\linlin&\circ\linlin&\circ\linlin&\circ\linlin&\circ\linlin&\circ\linlin&\circ\cr
             &            &            &\vlin\cr
             &            &            &\circ\cr
             &            &            &2\pdw8\cr
   }
}
\]

$\Psi=\{\alpha_1 = \frac12(\ve_1+\ve_8)-
\frac12(\ve_2+\ve_3+\ve_4+\ve_5+\ve_6+\ve_7),\,
\alpha_2 = \ve_1 + \ve_2, \alpha_3 = \ve_2 - \ve_1,\, 
\alpha_4 = \ve_3 - \ve_2,\,
\alpha_5 = \ve_4 - \ve_3,\,
\alpha_6 = \ve_5 - \ve_4,\,
\alpha_7 = \ve_6 - \ve_5
\}$

$\rho= \ve_2 + 2\ve_3 + 3\ve_4+ 4\ve_5 + 5\ve_6 - 
       \frac{17}2\ve_7 + \frac{17}2\ve_8
     = 17\alpha_1 + \frac{49}2\alpha_2 + 33\alpha_3 + 48\alpha_4
       + \frac{75}2\alpha_5 + 26\alpha_6 + \frac{27}2\alpha_7$
\medskip

\noindent
{\rm i)}
$\pi=\varpi_7 := \ve_6 + \frac12(\ve_8 - \ve_7)$ (minuscule)

$\dim\varpi_7=56$.

$C_{\varpi_7}=12\qquad\text{(see below)}$

$(\varpi_7,\rho) = \frac{27}2$

$q_{\varpi_7}(x) = 
\prod_{\varpi\in W_{E_7}\varpi_7}\bigl(x - \varpi - \frac{9}{8}
\bigr)$
\medskip

\noindent
{\rm ii)}
$\pi=\varpi_1 := \ve_8 - \ve_2$
(adjoint)

$\dim\varpi_1 = 133$

$C_{\varpi_1} = 36$

$(\varpi_1,\rho) = 17$

$q_{\varpi_1}(x) = (x-\frac12)
\prod_{\alpha\in\Sigma(E_7)}
\bigl(x - \alpha - \frac{17}{36}\bigr)$
\medskip

Let $\pi$ be the minuscule representation $\varpi_7$ in i).
Then the diagram of the partially ordered set of the weights of $\pi$ is as 
follows.
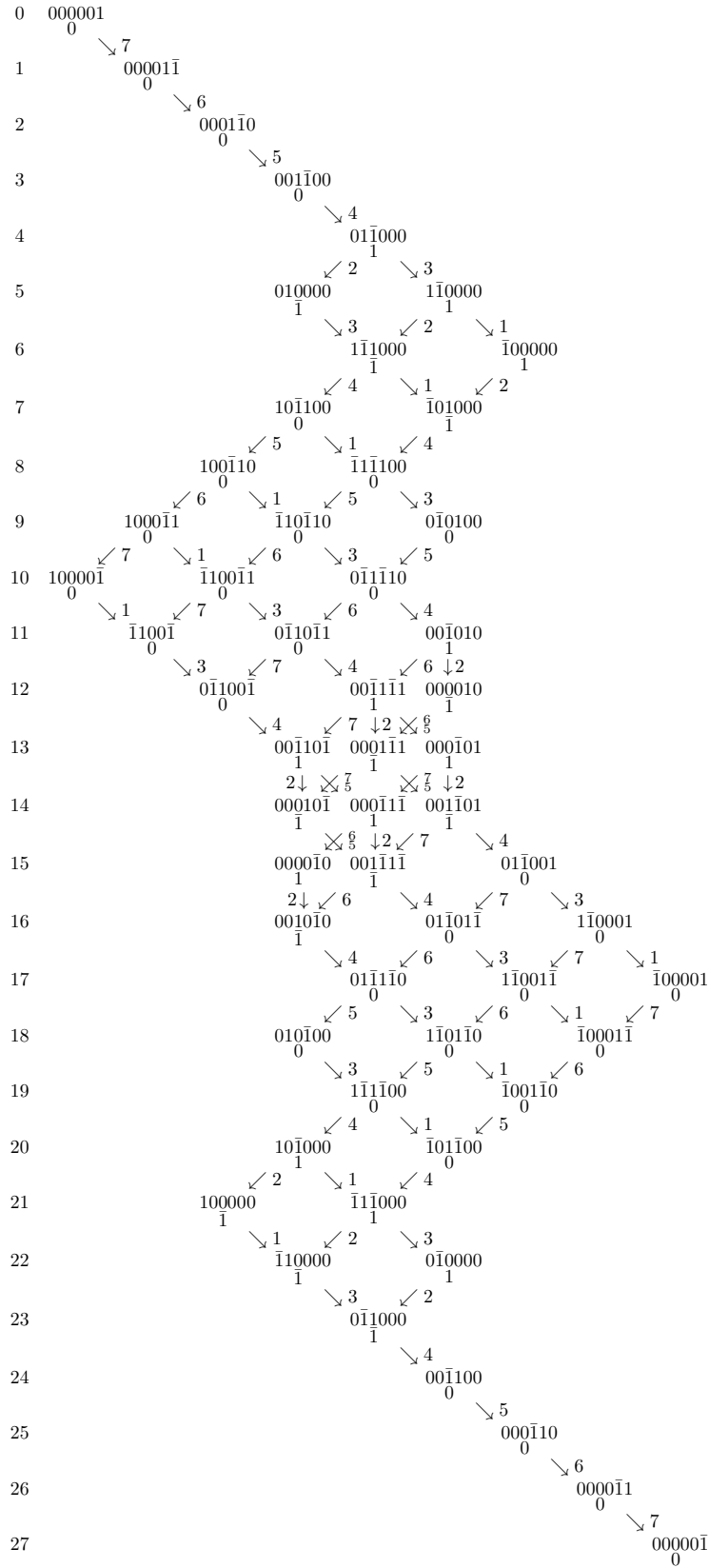
\[
\scalebox{.75}{$
\begin{matrix}
0&  00\mL00001& & & & & & & &\\
&  \lD7&&&&&&&&\\
1&  &00\mL0001{\bar1}&&&&&&&\\
&  &\lD6&&&&&&&\\
2&  &&00\mL001{\bar1}0&&&&&&\\
&  &&\lD5&&&&&&\\
3&  &&&00\mL10{\bar1}00\\
&  &&&\lD4&&&&&\\
4&  &&&&01\mL{\bar1}1000&&&&\\
&  &&&\lU2&\lD3&&&&\\
5&  &&&01\mL0{\bar1}000&&1{\bar1}\mL01000\\
&  &&&\lD3&\lU2&\lD1&&&\\
6&  &&&&1{\bar1}\mL1{\bar1}000&&{\bar1}0\mL01000\\
&  &&&\lU4&\lD1&\lU2\\
7&  &&&10\mL{\bar1}0100&&{\bar1}0\mL1{\bar1}000\\
&  &&\lU5&\lD1&\lU4\\
8&  &&10\mL00{\bar1}10&&{\bar1}1\mL{\bar1}0100\\
&  &\lU6&\lD1&\lU5&\lD3\\
9&  &10\mL000{\bar1}1&&{\bar1}1\mL00{\bar1}10&&0{\bar1}\mL00100\\
&  \lU7&\lD1&\lU6&\lD3&\lU5\\
10&  10\mL0000{\bar1}&&{\bar1}1\mL000{\bar1}1&&0{\bar1}\mL10{\bar1}10\\
&  \lD1&\lU7&\lD3&\lU6&\lD4\\
11&  &{\bar1}1\mL000{\bar1}&&0{\bar1}\mL100{\bar1}1&&00\mL{\bar1}1010\\
&  &\lD3&\lU7&\lD4&\lU6&\downarrow\!2\\
12&  &&0{\bar1}\mL1000{\bar1}&&00\mL{\bar1}11{\bar1}1&00\mL0{\bar1}010\\
&  &&\lD4&\lU7&\makebox[0pt]{\quad\ \lD{{}_5}}\makebox[0pt]{\quad\ \lU{{}^6}}
  \downarrow\!2\\
13&  &&&00\mL{\bar1}110{\bar1}&00\mL0{\bar1}1{\bar1}1&00\mL01{\bar1}01\\
&  &&&\lU{{}^7}\lD{{}_5}\!\!\!\!2\!\downarrow&\lU{{}^7}\lD{{}_5}&\downarrow\!2\\
14&  &&&00\mL0{\bar1}10{\bar1}&00\mL01{\bar1}1{\bar1}&00\mL1{\bar1}{\bar1}01\\
&  &&&\lU{{}^6}\lD{{}_5}&\makebox[0pt]{\quad\lU7}\downarrow\!2&\lD4\\
15&  &&&00\mL010{\bar1}0&00\mL1{\bar1}{\bar1}1{\bar1}&&01\mL{\bar1}0001\\
&  &&&\lU6\!\!\!2\!\downarrow&\lD4&\lU7&\lD3\\
16&  &&&00\mL1{\bar1}0{\bar1}0&&01\mL{\bar1}001{\bar1}&&1{\bar1}\mL00001\\
&  &&&\lD4&\lU6&\lD3&\lU7&\lD1\\
17&  &&&&01\mL{\bar1}01{\bar1}0&&1{\bar1}\mL0001{\bar1}&&{\bar1}0\mL00001\\
&  &&&\lU5&\lD3&\lU6&\lD1&\lU7\\
18&  &&&01\mL00{\bar1}00&&1{\bar1}\mL001{\bar1}0&&{\bar1}0\mL0001{\bar1}\\
&  &&&\lD3&\lU5&\lD1&\lU6\\
19&  &&&&1{\bar1}\mL10{\bar1}00&&{\bar1}0\mL001{\bar1}0\\
&  &&&\lU4&\lD1&\lU5\\
20&  &&&10\mL{\bar1}1000&&{\bar1}0\mL10{\bar1}00\\
&  &&\lU2&\lD1&\lU4\\
21&  &&10\mL0{\bar1}000&&{\bar1}1\mL{\bar1}1000\\
&  &&\lD1&\lU2&\lD3\\
22&  &&&{\bar1}1\mL0{\bar1}000&&0{\bar1}\mL01000\\
&  &&&\lD3&\lU2\\
23&  &&&&0{\bar1}\mL1{\bar1}000\\
&  &&&&\lD4\\
24&  &&&&&00\mL{\bar1}0100\\
&  &&&&&\lD5\\
25&  &&&&&&00\mL00{\bar1}10\\
&  &&&&&&\lD6\\
26&  &&&&&&&00\mL000{\bar1}1\\
&  &&&&&&&\lD7\\
27&  &&&&&&&&00\mL0000{\bar1}
\end{matrix}$}
\]
Here we use the similar notation as in Example~\ref{ex:E6}.

The type $A_6$ corresponding to $\{\alpha_1,\alpha_3,\ldots,\alpha_7\}$
is contained in type $E_7$.
The highest weights of the restriction $(E_7,\pi)|_{A_6}$ are
$\varpi_7=00\mL00001$, $\varpi_3 - \varpi_2=w_2w_4w_5w_6w_7\varpi_7=01\mL0{\bar1}000$, 
$\varpi_6 - \varpi_2=w_2w_4w_3w_1w_5w_4w_3(\varpi_3 - \varpi_2)=00\mL0{\bar1}010$ and
$\varpi_1-\varpi_2=w_2w_4w_5w_6w_7w_3w_4w_5w_6(\varpi_6 - \varpi_2)=10\mL0{\bar1}000$. 
Therefore
$(E_7,\pi)|_{A_6} = (A_6,\varpi_6) + (A_6,\varpi_2) +
(A_6,\varpi_5) + (A_6,\varpi_1)$ and
$C_{\varpi_7}=
\binom{6-1}{6-1}+\binom{6-1}{2-1}+\binom{6-1}{5-1}+\binom{6-1}{1-1}=12$.

Now use $\Psi'=-\Psi$ and put $\Theta_i = \Psi'\setminus\{\alpha'_i\}$.
Then
\begin{align*}
 \overline{\mathcal W}_{\Theta_1}(\pi) &=
  \bigl\{00\mL00001,\,\bar10\mL0100,\,\bar10\mL00001\bigr\},\\
 \overline{\mathcal W}_{\Theta_2}(\pi)  &=
  \bigl\{00\mL00001,\, 01\mL0{\bar1}000,\, 00\mL0{\bar1}010,\,
  10\mL0{\bar1}000\bigr\},\\
 \overline{\mathcal W}_{\Theta_3}(\pi)  &=
  \bigl\{00\mL00001,\, 0\bar1\mL01000,\, 0\bar1\mL00100,\, 1\bar1\mL00001,\,
  0\bar1\mL01000\bigr\},\\
 \overline{\mathcal W}_{\Theta_4}(\pi)  &=
  \bigl\{00\mL00001,\, 01\mL{\bar1}1000,\, 10\mL{\bar1}0100,\, 00\mL{\bar1}1010,\,
  01\mL{\bar1}0001,\,10\mL{\bar1}1000,\,00\mL{\bar1}0100\bigr\},\\
 \overline{\mathcal W}_{\Theta_5}(\pi)  &=
  \bigl\{00\mL00001,\, 00\mL10{\bar1}00,\, 10\mL00{\bar1}10,\,
  00\mL01{\bar1}01,\,01\mL00\bar100,\,00\mL00\bar110\bigr\},\\
 \overline{\mathcal W}_{\Theta_6}(\pi)  &=
  \bigl\{00\mL00001,\, 00\mL001{\bar1}0,\, 10\mL000{\bar1}1,\, 
   00\mL010{\bar1}0,\,00\mL000\bar11\bigr\},\\
 \overline{\mathcal W}_{\Theta_7}(\pi)  &=
  \bigl\{00\mL00001,\,00\mL0001\bar1,\,\bar10\mL0000\bar1,\,00\mL00001\bigr\}.
\end{align*}
From \eqref{eq:minmin} with 
$\ang{\alpha_i}{\varpi_i}=\ang{\alpha_i}{\rho}=\frac12\ang{\alpha_i}{\alpha_i}
=\frac1{12}$ and
\begin{align*}
\ang{\varpi_7}{\varpi_1}&=\frac1{12},
&\ang{\varpi_7}{\varpi_2}&=\frac18,
&\ang{\varpi_7}{\varpi_3}&=\frac16,
&\ang{\varpi_7}{\varpi_4}&=\frac14,\\
\ang{\varpi_7}{\varpi_5}&=\frac5{24},
&\ang{\varpi_7}{\varpi_6}&=\frac16,
&\ang{\varpi_7}{\varpi_7}&=\frac18,
\end{align*}
we have
\begin{align*}
q_{\pi,\Theta_1}(x;\lambda) &= 
 \Bigl(x - \frac1{12}\lambda\Bigr)
 \Bigl(x - \frac12\Bigr)
 \Bigl(x + \frac1{12}\lambda - \frac{17}{12}\Bigr),\\
q_{\pi,\Theta_2}(x;\lambda) &= 
 \Bigl(x - \frac18\lambda\Bigr)
 \Bigl(x - \frac1{24}\lambda - \frac5{12}\Bigr)
 \Bigl(x + \frac1{24}\lambda - 1\Bigr)
 \Bigl(x + \frac18\lambda - \frac74\Bigr),\\
q_{\pi,\Theta_3}(x;\lambda) &= 
 \Bigl(x - \frac16\lambda\Bigr)
 \Bigl(x - \frac1{12}\lambda - \frac5{12}\Bigr)
 \Bigl(x - \frac34\Bigr)
 \Bigl(x + \frac1{12}\lambda - \frac43\Bigr)
 \Bigl(x + \frac16\lambda - \frac{11}6\Bigr),\\
q_{\pi,\Theta_4}(x;\lambda) &= 
 \Bigl(x - \frac14\lambda\Bigr)
 \Bigl(x - \frac16\lambda - \frac13\Bigr)
 \Bigl(x - \frac1{12}\lambda - \frac7{12}\Bigr)
 \Bigl(x - \frac{11}{12}\Bigr) \\
 &\qquad\cdot
 \Bigl(x + \frac1{12}\lambda - \frac54\Bigr)
 \Bigl(x + \frac16\lambda - \frac53\Bigr)
 \Bigl(x + \frac14\lambda - 2\Bigr),\\
q_{\pi,\Theta_5}(x;\lambda) &= 
 \Bigl(x - \frac5{24}\lambda\Bigr)
 \Bigl(x - \frac18\lambda - \frac14\Bigr)
 \Bigl(x - \frac1{24}\lambda - \frac23\Bigr)
 \Bigl(x + \frac1{24}\lambda - \frac{13}{12}\Bigr) \\
 &\qquad\cdot
 \Bigl(x + \frac18\lambda - \frac32\Bigr)
 \Bigl(x + \frac5{24}\lambda - \frac{25}{12}\Bigr),\\
q_{\pi,\Theta_6}(x;\lambda) &= 
 \Bigl(x - \frac16\lambda\Bigr)
 \Bigl(x - \frac1{12}\lambda - \frac16\Bigr)
 \Bigl(x - \frac34\Bigr)
 \Bigl(x + \frac1{12}\lambda - \frac54\Bigr)
 \Bigl(x + \frac16\lambda - \frac{13}6\Bigr),\\
q_{\pi,\Theta_7}(x;\lambda) &= 
 \Bigl(x - \frac18\lambda\Bigr)
 \Bigl(x - \frac1{24}\lambda - \frac1{12}\Bigr)
 \Bigl(x + \frac1{24}\lambda - \frac56\Bigr)
 \Bigl(x + \frac18\lambda - \frac94\Bigr). 
\end{align*}
\end{exmp}
\newpage %

\begin{exmp}[$E_8$]\label{ex:E8}
\[
\offinterlineskip\vbox{
\halign{$ # $&$ # $&$ # $&$ # $&$ # $&$ # $&$ # $\cr
 \alpha_1    & \alpha_3   &\alpha_4 &\alpha_5& \alpha_6&\alpha_7&\alpha_8\pup7\cr
 \circ\linlin&\circ\linlin&\circ\linlin&\circ\linlin&\circ\linlin&\circ\linlin&\circ\cr
             &            &\vlin\cr
             &            &\circ\cr
             &            &\alpha_2\pdw8\cr
   }
}
\quad\quad
\vbox{
\halign{$ # $&$ # $&$ # $&$ # $&$ # $&$ # $&$ # $&$ # $\cr
 2 & 4 & 6 & 5 & 4 & 3 & 2\pup7\cr
 \circ\linlin&\circ\linlin&\circ\linlin&\circ\linlin&\circ\linlin&\circ\linlin&\circ\linlin&\bullet\cr
             &            &\vlin\cr
             &            &\circ\cr
             &            &3\pdw8\cr
   }
}\]

$\Psi=\{\alpha_1 = \frac12(\ve_1+\ve_8)-
\frac12(\ve_2+\ve_3+\ve_4+\ve_5+\ve_6+\ve_7),\,
\alpha_2 = \ve_1 + \ve_2, \alpha_3 = \ve_2 - \ve_1,\, 
\alpha_4 = \ve_3 - \ve_2,\,
\alpha_5 = \ve_4 - \ve_3,\,
\alpha_6 = \ve_5 - \ve_4,\,
\alpha_7 = \ve_6 - \ve_5,\,
\alpha_8 = \ve_7 - \ve_6
\}$

$\rho= \ve_2 + 2\ve_3 + 3\ve_4+ 4\ve_5 + 5\ve_6 + 
       6\ve_7 + 23\ve_8
     = 46\alpha_1 + 68\alpha_2 + 91\alpha_3
       + 135\alpha_4 + 110\alpha_5 + 84\alpha_6
       + 57\alpha_7 + 29\alpha_8$
\medskip

\noindent
{\rm i)}
$\pi=\alpha_{\max} := \ve_7 + \ve_8$ (adjoint)

$\dim\alpha_{\max} = 248$ ($m_{\alpha_{\max}}(0)=8$)

$C_{\alpha_{\max}}=60$

$(\alpha_{\max},\rho) = 29$

$q_{\alpha_{\max}}(x) = (x-\frac12)
\prod_{\alpha\in\Sigma(E_8)}\bigl(x - \alpha - \frac{29}{60}
 \bigr)$\medskip
 
Let $\pi$ be the adjoint representation $\alpha_{\max}$
and $\alpha_{\max}=\sum_{i=1}^8n_i\alpha_i$,
that is, $n_1=2$, $n_2=3,\ldots.$
Put $\Theta_i=\Psi\setminus\{\alpha_i\}$
for $i=1,\ldots,8$.
The irreducible decomposition of $\mathfrak g$ as a ${\mathfrak g}_{\Theta_i}$-module
is given by Proposition~\ref{prop:free}~ii).
In this case $L_{\Theta_i}$ in the proposition equals
$\{-n_i, -n_i+1, \ldots, n_i\}$.
Suppose $\mathbf m\in L_{\Theta_i}\setminus\{0\}$.
Then $V(\mathbf m)$ is a minuscule representation 
since $E_8$ is simply-laced.
Let $\varpi_i$ ($j=1,\ldots,8$) be the fundamental weights.
If we write the lowest weight and the highest weight of $V(\mathbf m)$
by $\alpha_{\mathbf m}=\sum_{j=1}^8c_j\varpi_j$ 
and $\alpha'_{\mathbf m}=\sum_{j=1}^8c'_j\varpi_j$ respectively, 
we clearly have
\[
c_i=\begin{cases}
1 & \text{if }\mathbf m \ne 1,-n_i,\\
2 & \text{if }\mathbf m = 1,
\end{cases}
\qquad
c'_i=\begin{cases}
-1 & \text{if }\mathbf m \ne -1,n_i,\\
-2 & \text{if }\mathbf m =- 1,
\end{cases}
\]
and $\alpha_{\mathbf m}=-\alpha'_{-\mathbf m}$.
Since we know the highest weights and the lowest weights of 
minuscule representations of ${\mathfrak g}_{\Theta_i}$
by the previous examples,
starting with  $\alpha_{\max}=\varpi_8=00\mL000001$,
we can determine $\alpha_{\mathbf m}$ and $\alpha'_{\mathbf m}$
for $\mathbf m\in L_{\Theta_i}\setminus\{0\}$ step by step.
For example, suppose $i=4$. Then $L_{\Theta_4}=\{-6, -5, \ldots,6\}$
and we have
\begin{align*}
V(6): &
\begin{cases}
00\mL000001 & \text{h.w.} \\
00\mL10\bar1000 & \text{l.w.}
\end{cases}
&\rightarrow\quad
&00\mL10\bar1000-\alpha_4=01\mL{\bar1}10000 \text{ is a weight of }V(5)\\
V(5): &
\begin{cases}
01\mL{\bar1}10000 & \text{h.w.} \\
\bar10\mL1{\bar1}0000 & \text{l.w.}
\end{cases}
&\rightarrow\quad
&\bar10\mL1{\bar1}0000-\alpha_4=\bar11\mL{\bar1}01000 \text{ is a weight of }V(4)\\
V(4): &
\begin{cases}
10\mL{\bar1}01000 & \text{h.w.} \\
\bar11\mL{\bar1}01000 & \\
0\bar1\mL10000\bar1 & \text{l.w.}
\end{cases}
&\rightarrow\quad
&0\bar1\mL10000\bar1-\alpha_4=00\mL{\bar1}1100\bar1 \text{ is a weight of }V(3)\\
V(3): &
\begin{cases}
00\mL{\bar1}10100 & \text{h.w.} \\
00\mL{\bar1}1100\bar1 & \\
00\mL1{\bar1}00\bar10 & \text{l.w.}
\end{cases}
&\rightarrow\quad
&00\mL1{\bar1}00\bar10-\alpha_4=01\mL{\bar1}010\bar10 \text{ is a weight of }V(2)\\
V(2): &
\begin{cases}
01\mL{\bar1}00010 & \text{h.w.} \\
01\mL{\bar1}010\bar10 & \\
\bar10\mL100\bar100 & \text{l.w.}
\end{cases}
&\rightarrow\quad
&\bar10\mL100\bar100-\alpha_4=\bar11\mL{\bar1}11\bar100 \text{ is a weight of }V(1)\\
V(1): &
\begin{cases}
10\mL{\bar1}10001 & \text{h.w.} \\
\bar11\mL{\bar1}11\bar100 & \\
0\bar1\mL2{\bar1}{\bar1}000 & \text{l.w.}
\end{cases}
&&
\end{align*}

On the other hand,
the non-trivial irreducible subrepresentations of $V(0)$ correspond to 
the connected parts of Dynkin diagram of $\Theta_i$.
If $\sum_{j=1}^8c_j\varpi_j$ is a lowest weight of such subrepresentations,
then $c_i=1$.
Hence, if $i=4$, the lowest weights of the non-trivial irreducible subrepresentations of $V(0)$
are
\[
\bar1\bar1\mL100000,
\qquad
00\mL1{\bar2}0000,
\qquad
00\mL10\bar100\bar1.
\]
Thus we get
\[\begin{aligned}
\overline{\mathcal W}_{\Theta_4}(\pi)&=
\{
00\mL10\bar1000,
\bar10\mL1{\bar1}0000,
0\bar1\mL10000\bar1,
00\mL1{\bar1}00\bar10,
\bar10\mL100\bar100,
0\bar1\mL2{\bar1}{\bar1}000
\} \\
&\qquad\cup\{0\}\cup
\{
\bar1\bar1\mL100000,
00\mL1{\bar2}0000,
00\mL10\bar100\bar1
\}\\
&\qquad\cup
\{
-10\mL{\bar1}10001,
-01\mL{\bar1}00010,
-00\mL{\bar1}10100,
-10\mL{\bar1}01000,
-01\mL{\bar1}10000,
-00\mL000001
\}.
\end{aligned}\]
Put $\lambda_{\Theta_i}=\lambda\varpi_i$.
Then, by \eqref{eq:admin}, we have
\[\begin{aligned}
q_{\pi, \Theta_4}(x; \lambda)&=\Bigl(x-\frac12\Bigr)\Bigl(x-\frac{9}{20}\Bigr)\Bigl(x-\frac{7}{15}\Bigr)\Bigl(x-\frac{5}{12}\Bigr)\Bigl(x-\frac{1}{10}\lambda-\frac{9}{10}\Bigr)\\
&\cdot\Bigl(x-\frac{1}{12}\lambda-\frac{5}{6}\Bigr)\Bigl(x-\frac{1}{15}\lambda-\frac{11}{15}\Bigr)
      \Bigl(x-\frac{1}{20}\lambda-\frac{13}{20}\Bigr)\Bigl(x-\frac{1}{30}\lambda-\frac{17}{30}\Bigr)\\
&\cdot\Bigl(x-\frac{1}{60}\lambda-\frac{1}{2}\Bigr)\Bigl(x+\frac{1}{60}\lambda-\frac{7}{20}\Bigr)
      \Bigl(x+\frac{1}{30}\lambda-\frac{4}{15}\Bigr)\Bigl(x+\frac{1}{20}\lambda-\frac{1}{5}\Bigr)\\
&\cdot\Bigl(x+\frac{1}{15}\lambda-\frac{2}{15}\Bigr)\Bigl(x+\frac{1}{12}\lambda-\frac{1}{12}\Bigr)\Bigl(x+\frac{1}{10}\lambda\Bigr).
\end{aligned}\]
Similarly we get
\begin{align*}
\overline{\mathcal W}_{\Theta_1}(\pi)&=\{10\mL00000{\bar1}, 2{\bar1}\mL000000\}
\cup\{0\}\cup\{10\mL0000{\bar1}0\}
\cup\{-{\bar1}0\mL010000, -00\mL000001\},\\
\overline{\mathcal W}_{\Theta_2}(\pi)&=\{{\bar1}0\mL010000, 00\mL0100{\bar1}0, 00\mL{\bar1}20000\}
\cup\{0\}\cup\{{\bar1}0\mL01000{\bar1}\}\\
&\qquad\cup\{-00\mL0{\bar1}0100, -01\mL0{\bar1}0000, -00\mL000001\},\\
\overline{\mathcal W}_{\Theta_3}(\pi)&=\{01\mL0{\bar1}0000, {\bar1}1\mL00000{\bar1}, 01\mL000{\bar1}00, {\bar1}2\mL{\bar1}00000\}
\cup\{0\}\cup\{{\bar2}1\mL000000, 01\mL0{\bar1}000{\bar1}\}\\
&\qquad\cup\{-1{\bar1}\mL000010, -0{\bar1}\mL001000, -1{\bar1}\mL010000, -00\mL000001\},\\
\overline{\mathcal W}_{\Theta_5}(\pi)&=\{00\mL001{\bar1}00, 0{\bar1}\mL001000, 00\mL0{\bar1}100{\bar1}, {\bar1}0\mL0010{\bar1}0, 00\mL{\bar1}02{\bar1}00\}\\
&\qquad\cup\{0\}\cup\{{\bar1}0\mL0{\bar1}1000, 00\mL001{\bar1}0{\bar1}\}\\
&\qquad\cup\{-01\mL00{\bar1}001, -00\mL01{\bar1}010, -10\mL00{\bar1}100, -00\mL10{\bar1}000, -00\mL000001\},\\
\overline{\mathcal W}_{\Theta_6}(\pi)&=\{00\mL0001{\bar1}0, 00\mL0{\bar1}0100, {\bar1}0\mL00010{\bar1}, 00\mL00{\bar1}2{\bar1}0\}
\cup\{0\}\cup\{0{\bar1}\mL000100, 00\mL0001{\bar1}{\bar1}\}\\
&\qquad\cup\{-00\mL010{\bar1}01, -10\mL000{\bar1}10, -00\mL001{\bar1}00, -00\mL000001\},\\
\overline{\mathcal W}_{\Theta_7}(\pi)&=\{00\mL00001{\bar1}, {\bar1}0\mL000010, 00\mL000{\bar1}2{\bar1}\}
\cup\{0\}\cup\{00\mL0{\bar1}0010, 00\mL00001{\bar2}\}\\
&\qquad\cup\{-10\mL0000{\bar1}1, -00\mL0001{\bar1}0, -00\mL000001\},\\
\overline{\mathcal W}_{\Theta_8}(\pi)&=\{00\mL000001, 00\mL0000{\bar1}2\}
\cup\{0\}\cup\{{\bar1}0\mL000001\}
\cup\{-00\mL00001{\bar1}, -00\mL000001\},
\end{align*}
and 
\begin{align*}
q_{\pi, \Theta_1}(x; \lambda)&=\Bigl(x-\frac12\Bigr)\Bigl(x-\frac{3}{10}\Bigr)
\Bigl(x-\frac{1}{30}\lambda-\frac{23}{30}\Bigr)\Bigl(x-\frac{1}{60}\lambda-\frac{1}{2}\Bigr)\\
&\cdot\Bigl(x+\frac{1}{60}\lambda-\frac{7}{60}\Bigr)\Bigl(x+\frac{1}{30}\lambda\Bigr),\\
q_{\pi, \Theta_2}(x; \lambda)&=\Bigl(x-\frac12\Bigr)\Bigl(x-\frac{11}{30}\Bigr)
\Bigl(x-\frac{1}{20}\lambda-\frac{17}{20}\Bigr)\Bigl(x-\frac{1}{30}\lambda-\frac{2}{3}\Bigr)\Bigl(x-\frac{1}{60}\lambda-\frac{1}{2}\Bigr)\\
&\cdot\Bigl(x+\frac{1}{60}\lambda-\frac{13}{60}\Bigr)\Bigl(x+\frac{1}{30}\lambda-\frac{1}{10}\Bigr)\Bigl(x+\frac{1}{20}\lambda\Bigr),\\
q_{\pi, \Theta_3}(x; \lambda)&=\Bigl(x-\frac12\Bigr)\Bigl(x-\frac{7}{15}\Bigr)\Bigl(x-\frac{23}{60}\Bigr)\\
&\cdot\Bigl(x-\frac{1}{15}\lambda-\frac{13}{15}\Bigr)\Bigl(x-\frac{1}{20}\lambda-\frac{3}{4}\Bigr)\Bigl(x-\frac{1}{30}\lambda-\frac{3}{5}\Bigr)\Bigl(x-\frac{1}{60}\lambda-\frac{1}{2}\Bigr)\\
&\cdot\Bigl(x+\frac{1}{60}\lambda-\frac{17}{60}\Bigr)\Bigl(x+\frac{1}{30}\lambda-\frac{1}{6}\Bigr)\Bigl(x+\frac{1}{20}\lambda-\frac{1}{10}\Bigr)\Bigl(x+\frac{1}{15}\lambda\Bigr),\\
q_{\pi, \Theta_5}(x; \lambda)&=\Bigl(x-\frac12\Bigr)\Bigl(x-\frac{5}{12}\Bigr)\Bigl(x-\frac{13}{30}\Bigr)
\Bigl(x-\frac{1}{12}\lambda-\frac{11}{12}\Bigr)\Bigl(x-\frac{1}{15}\lambda-\frac{4}{5}\Bigr)\\
&\cdot\Bigl(x-\frac{1}{20}\lambda-\frac{7}{10}\Bigr)\Bigl(x-\frac{1}{30}\lambda-\frac{3}{5}\Bigr)
      \Bigl(x-\frac{1}{60}\lambda-\frac{1}{2}\Bigr)\Bigl(x+\frac{1}{60}\lambda-\frac{19}{60}\Bigr)\\
&\cdot\Bigl(x+\frac{1}{30}\lambda-\frac{7}{30}\Bigr)\Bigl(x+\frac{1}{20}\lambda-\frac{3}{20}\Bigr)\Bigl(x+\frac{1}{15}\lambda-\frac{1}{15}\Bigr)\Bigl(x+\frac{1}{12}\lambda\Bigr),\\
q_{\pi, \Theta_6}(x; \lambda)&=\Bigl(x-\frac12\Bigr)\Bigl(x-\frac{11}{30}\Bigr)\Bigl(x-\frac{9}{20}\Bigr)\\
&\cdot\Bigl(x-\frac{1}{15}\lambda-\frac{14}{15}\Bigr)\Bigl(x-\frac{1}{20}\lambda-\frac{3}{4}\Bigr)\Bigl(x-\frac{1}{30}\lambda-\frac{19}{30}\Bigr)\Bigl(x-\frac{1}{60}\lambda-\frac{1}{2}\Bigr)\\
&\cdot\Bigl(x+\frac{1}{60}\lambda-\frac{4}{15}\Bigr)\Bigl(x+\frac{1}{30}\lambda-\frac{1}{6}\Bigr)\Bigl(x+\frac{1}{20}\lambda-\frac{1}{20}\Bigr)\Bigl(x+\frac{1}{15}\lambda\Bigr),\\
q_{\pi, \Theta_7}(x; \lambda)&=\Bigl(x-\frac12\Bigr)\Bigl(x-\frac{3}{10}\Bigr)\Bigl(x-\frac{7}{15}\Bigr)
\Bigl(x-\frac{1}{20}\lambda-\frac{19}{20}\Bigr)\Bigl(x-\frac{1}{30}\lambda-\frac{2}{3}\Bigr)\\
&\cdot\Bigl(x-\frac{1}{60}\lambda-\frac{1}{2}\Bigr)\Bigl(x+\frac{1}{60}\lambda-\frac{11}{60}\Bigr)\Bigl(x+\frac{1}{30}\lambda-\frac{1}{30}\Bigr)\Bigl(x+\frac{1}{20}\lambda\Bigr),\\
q_{\pi, \Theta_8}(x; \lambda)&=\Bigl(x-\frac12\Bigr)\Bigl(x-\frac{1}{5}\Bigr)
\Bigl(x-\frac{1}{30}\lambda-\frac{29}{30}\Bigr)\Bigl(x-\frac{1}{60}\lambda-\frac{1}{2}\Bigr)\\
&\cdot\Bigl(x+\frac{1}{60}\lambda-\frac{1}{60}\Bigr)\Bigl(x+\frac{1}{30}\lambda\Bigr).
\end{align*}

\end{exmp}
\begin{exmp}[$F_4$]\label{ex:F4}
\[
\offinterlineskip
\vbox{
\halign{$ # $&$ # $&$ # $&$ # $\cr
 \alpha_1 & \alpha_2 & \alpha_3  & \alpha_4\pup7\cr
 \circ\linlin&\circ\rarw&\circ\linlin&\circ\cr
    }
}
\quad\quad
\vbox{
\halign{$ # $&$ # $&$ # $&$ # $&$ # $\cr
 & 2 & 3 & 4 & 2\pup7\cr
 \bullet\linlin&\circ\linlin&\circ\rarw&\circ\linlin&\circ\cr
    }
}
\]

$\Psi=\{
\alpha_1 = \ve_2 - \ve_3,\,
\alpha_2 = \ve_3 - \ve_4,\,
\alpha_3 = \ve_4,\,
\alpha_4 = \frac12(\ve_1 - \ve_2 - \ve_3 - \ve_4)
\}$

$\rho = \frac{11}2\ve_1 + \frac52\ve_2 + \frac32\ve_3 + \frac12\ve_4
      = 8\alpha_1 + 15\alpha_2 + 21\alpha_3 + 11\alpha_4$
\medskip

\noindent
{\rm i)}
$\pi=\varpi_4:=\ve_1=\alpha_1+2\alpha_2+3\alpha_3+2\alpha_4$ (dominant short root)

$\dim\varpi_4=26$ ($m_{\varpi_4}(0)=2$) 

$C_{\varpi_4} =\sum_{\nu=1}^4 (\pm\ve_\nu,\ve_1)^2 + 
\frac14\sum(\pm\ve_1\pm\ve_2\pm\ve_3\pm\ve_4,\ve_1)^2
=2 + \frac{16}4=6$

$(\varpi_4,\rho) = \frac{11}2$

$q_{\varpi_4}(x) = \bigl(x - 1\bigr)
  \prod_{\substack{\alpha\in\Sigma(F_4)\\|\beta|<|\alpha_{\max}|}}
    \bigr(x - \beta - \frac{11}{12}\bigr)$
\medskip

\noindent
{\rm ii)}
$\pi=\varpi_1:=\ve_1+\ve_2$ (adjoint)

$\dim\varpi_1 = 52$

$C_{\varpi_1} = 18$

$(\varpi_1,\rho) = 8$

$q_{\varpi_1}(x)=(x-\frac12)
\prod_{\substack{\alpha\in\Sigma(F_4)\\|\alpha|=|\alpha_{\max}|}}
 \bigl(x-\alpha-\frac49\bigr)
\prod_{\substack{\beta\in\Sigma(F_4)\\|\beta|<|\alpha_{\max}|}}
 \bigl(x-\beta-\frac{17}{36}\bigr)$
\medskip

Let $\pi$ be the representation $\varpi_4$ in i).
Then the diagram of the partially ordered set of the weights of $\pi$ is as 
follows.
Here the weight $00\RA00$ is the only weight with the multiplicity 2 and hence 
indicated by $[00\RA00]$.
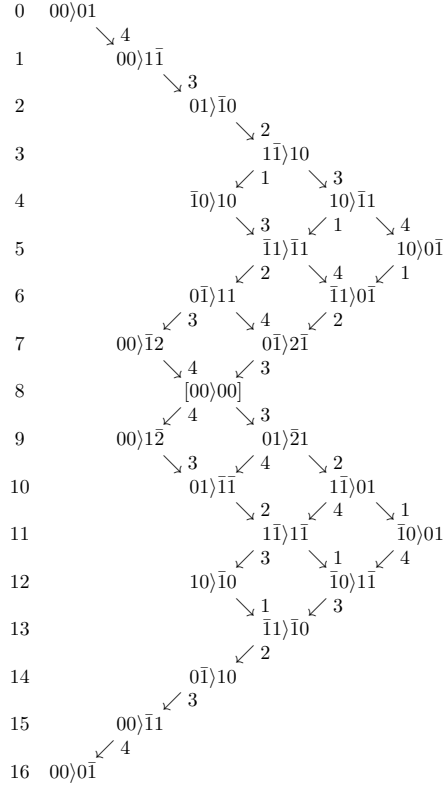
\[
\scalebox{0.75}{$
 \begin{matrix}
0&  00\RA01   &           &           &           &           &           \\
 &  \lD4      &           &           &           &           &           \\
1&            & 00\RA1{\bar1}   &     &           &           &           \\
 &            & \lD3      &           &           &           &           \\
2&            &           & 01\RA{\bar1}0 &       &           &           \\
 &            &           & \lD2      &           &           &           \\
3&            &           &           & 1{\bar1}\RA10 &       &           \\
 &            &           & \lU1      & \lD3      &           &           \\
4&            &           & {\bar1}0\RA10   &     & 10\RA{\bar1}1   &     \\
 &            &           & \lD3      & \lU1      & \lD4      &           \\
5&            &           &           & {\bar1}1\RA{\bar1}1&& 10\RA0{\bar1}\\
 &            &           & \lU2      & \lD4      & \lU1      &           \\
6&            &           & 0{\bar1}\RA11   &     & {\bar1}1\RA0{\bar1}  &\\
 &            &  \lU3     & \lD4      & \lU2      &           &           \\
7&            & 00\RA{\bar1}2   &     & 0{\bar1}\RA2{\bar1} & &           \\
 &            & \lD4      & \lU3      &           &           &           \\
8&            &           & [00\RA00] &           &           &           \\
 &            & \lU4      & \lD3      &           &           &           \\
9&            & 00\RA1{\bar2}   &     & 01\RA{\bar2}1   &     &           \\
 &            & \lD3      & \lU4      &  \lD2     &           &           \\
10&           &           & 01\RA{\bar1}{\bar1}   && 1{\bar1}\RA01&        \\
 &            &           & \lD2      &  \lU4     & \lD1      &           \\
11&           &           &           &1{\bar1}\RA1{\bar1}&& {\bar1}0\RA01\\
 &            &           & \lU3      & \lD1      & \lU4      &           \\
12&           &           & 10\RA{\bar1}0   & & {\bar1}0\RA1{\bar1}   &   \\
 &            &           & \lD1      & \lU3      &           &           \\
13&           &           &           & {\bar1}1\RA{\bar1}0 & &           \\
 &            &           & \lU2      &           &           &           \\
14&           &           & 0{\bar1}\RA10   &     &           &           \\
 &            & \lU3      &           &           &           &           \\
15&           & 00\RA{\bar1}1   &     &           &           &           \\
 & \lU4       &           &           &           &           &           \\
16& 00\RA0{\bar1}   &     &           &           &           &           
 \end{matrix}
 $}
\]
Now use $\Psi'=\{\alpha'_1=-\alpha_1,\ldots,\alpha'_4=-\alpha_4\}$
and put $\Theta_i = \Psi'\setminus\{\alpha'_i\}$.
Then we have
\begin{align*}
 \overline{\mathcal W}_{\Theta_1}(\pi) &=
  \bigl\{00\RA01,\,\bar10\RA10,\,\bar10\RA01\bigr\},\\
 \overline{\mathcal W}_{\Theta_2}(\pi)  &=
  \bigl\{00\RA01,\,1\bar1\RA10,\, 0\bar1\RA11,\,1\bar1\RA01,\,0\bar1\RA10\bigr\},\\
 \overline{\mathcal W}_{\Theta_3}(\pi)  &=
  \bigl\{00\RA01,\,01\RA\bar10,\,10\RA\bar11,\, 00\RA\bar12,\,00\RA00,\,01\RA\bar21,\,
  10\RA\bar10,\,00\RA\bar11\},\\
 \overline{\mathcal W}_{\Theta_4}(\pi)  &=
  \bigl\{00\RA0\bar1,\,00\RA1\bar1,\,10\RA0\bar1,\, 00\RA00,\,00\RA1\bar2,\,
  00\RA0\bar1\bigr\}
\end{align*}
and
\begin{align*}
q_{\pi,\Theta_1}(x;\lambda) &= 
 \Bigl(x - \frac16\lambda\Bigr)
 \Bigl(x - \frac12\Bigr)
 \Bigl(x + \frac16\lambda - \frac43\Bigr),\\
q_{\pi,\Theta_2}(x;\lambda) &= 
 \Bigl(x - \frac13\lambda\Bigr)
 \Bigl(x - \frac16\lambda - \frac13\Bigr)
 \Bigl(x - \frac34\Bigr)
 \Bigl(x + \frac16\lambda - \frac76\Bigr)
 \Bigl(x + \frac13\lambda - \frac53\Bigr),\\
q_{\pi,\Theta_3}(x;\lambda) &=  
 \Bigl(x - \frac14\lambda\Bigr)
 \Bigl(x - \frac16\lambda - \frac16\Bigr)
 \Bigl(x - \frac1{12}\lambda - \frac5{12}\Bigr)
 \Bigl(x - \frac56\Bigr) \\
 &\qquad\cdot
 \Bigl(x - 1\Bigr)
 \Bigl(x + \frac1{12}\lambda - 1\Bigr)
 \Bigl(x + \frac16\lambda - \frac43\Bigr)
 \Bigl(x + \frac14\lambda - \frac74\Bigr),\\
q_{\pi,\Theta_4}(x;\lambda) &= 
 \Bigl(x - \frac{11}6\lambda\Bigr)
 \Bigl(x - \frac1{12}\lambda - \frac1{12}\Bigr)
 \Bigl(x - \frac12\Bigr)
 \Bigl(x - 1\Bigr) \\
 &\qquad\cdot
 \Bigl(x + \frac1{12}\lambda - 1\Bigr)
 \Bigl(x + \frac16\lambda - \frac{11}6\Bigr).
\end{align*}
The extremal low weights of $\pi$ with respect to $\Psi'$ are as follows:
\begin{align*}
\varpi_{\alpha'_1}&=\varpi_4-\alpha_4-\alpha_3-\alpha_2
  =\alpha_1+\alpha_2+2\alpha_3+\alpha_4,\\
\varpi_{\alpha'_2}&=\varpi_4-\alpha_4-\alpha_3
  =\alpha_1+2\alpha_2+2\alpha_3+\alpha_4,\\
\varpi_{\alpha'_3}&=\varpi_4-\alpha_4
  =\alpha_1+2\alpha_2+3\alpha_3+\alpha_4,\\
\varpi_{\alpha'_4}&=\varpi_4
  =\alpha_1+2\alpha_2+3\alpha_3+2\alpha_4.
\end{align*}
None of them is a member of $\Sigma(\mathfrak g_\Theta)\cup\{0\}$
for any $\Theta\subsetneq\Psi'$.
Hence by Proposition~\ref{prop:free}~i)
and Lemma~\ref{lem:free_nonzero},
the functions
$r_{\alpha'_i}(\lambda)\ (i=1,2,3,4)$
are not identically zero.
\end{exmp}

\begin{exmp}[$G_2$]\label{ex:G2}
\[
\offinterlineskip
\vbox{
\halign{$ # $&$ # $\cr
 \alpha_1 & \alpha_2\pup7\cr
 \circ\!\Lleftarrow\;&\circ\cr
    }
}
\quad\quad\quad\quad
\vbox{
\halign{$ # $&$ # $&$ # $\cr
 3 & 2\pup7\cr
 \circ\!\Lleftarrow\;&\circ\linlin&\bullet\cr
    }
}\]

$\Psi=\{
\alpha_1 = \ve_1 - \ve_2,\,
\alpha_2 = -2\ve_1 + \ve_2 + \ve_3
\}$

$\rho = -\ve_1 -2\ve_2 + 3\ve_3
      = 5\alpha_1 + 3\alpha_2$
\medskip

\noindent
{\rm i)}
$\pi=\varpi_1 := -\ve_2 + \ve_3=2\alpha_1+\alpha_2$ (multiplicity free)

$\dim\varpi_1=7$

$C_{\varpi_1} = \frac12\big(2\sum_{1\le i<j\le 3}(\ve_i-\ve_j,\ve_1-\ve_2)^2
+(0,\ve_1-\ve_2)^2\big) = 6$

$%
(\varpi_1,\rho) = 5$

$q_{\varpi_1}(x) = \bigl(x-1\bigr)\prod_{1\le i<j\le 3}
           \bigl((x -\frac56)^2 
           - (\ve_i - \ve_j)^2\bigr)$
\medskip

\noindent
{\rm ii)}
$\pi=\varpi_2:=-\ve_1 - \ve_2 + 2\ve_3=3\alpha_1+2\alpha_2$ (adjoint)

$\dim\varpi_2=14$

$C_{\varpi_2}=24$

$(\varpi_2,\rho) = 9$

$q_{\varpi_2}(x)=(x-\frac12)
 \prod_{\substack{\alpha\in\Sigma(G_2)\\|\alpha|=|\alpha_{\max}|}}
 \bigl(x-\alpha-\frac38\bigr)
 \prod_{\substack{\beta\in\Sigma(G_2)\\|\beta|<|\alpha_{\max}|}}
 \bigl(x-\beta-\frac{11}{24}\bigr)$
\medskip

Consider the representation $\pi$ with the highest weight $\varpi_1$.
Then as is shown in \cite{FH}, the weights of $\pi$ are indicated by
\[
 \varepsilon_2-\varepsilon_3\xrightarrow{\alpha_1}
 \varepsilon_1-\varepsilon_3\xrightarrow{\alpha_2}
 -\varepsilon_1+\varepsilon_2\xrightarrow{\alpha_1}
 \ 0\ \xrightarrow{\alpha_1}
 \varepsilon_1-\varepsilon_2\xrightarrow{\alpha_2}
 -\varepsilon_1+\varepsilon_3\xrightarrow{\alpha_1}
 -\varepsilon_2+\varepsilon_3
\]
and therefore
\begin{align*}
\overline{\mathcal W}_{\{\alpha_1\}}(\pi)&=\{\varepsilon_2-\varepsilon_3,\,
 -\varepsilon_1+\varepsilon_2,\, -\varepsilon_1+\varepsilon_3\},\\
\overline{\mathcal W}_{\{\alpha_2\}}(\pi)&=\{\varepsilon_2-\varepsilon_3,\,
 \varepsilon_1-\varepsilon_3,\, 0,\, \varepsilon_1-\varepsilon_2,\,
 -\varepsilon_2+\varepsilon_3\}.
\end{align*}
For $\lambda\in\mathfrak{a}_\Theta^*$ we put $\lambda_\Theta=\lambda_1\varpi_1 + \lambda_2\varpi_2$.
Then $\lambda_1=0$ (resp.$~\lambda_2=0$) if $\Theta=\{\alpha_1\}$ 
(resp.~$\{\alpha_2\}$) and
\begin{align*}
 q_{\pi,\{\alpha_1\}}(x;\lambda) &=
 \Bigl(x+\frac{\lambda_2}2\Bigr)
 \Bigl(x-\frac{(\alpha_1+\alpha_2,\rho)}6\Bigr)
 \Bigl(x-\frac{\lambda_2}2 - \frac{(3\alpha_1+2\alpha_2,\rho)}6\Bigr)\\
 &= \Bigl(x+\frac{\lambda_2}2\Bigr)
    \Bigl(x-\frac23\Bigr)\Bigl(x-\frac{\lambda_2}2-\frac32\Bigr),
\\
 q_{\pi,\{\alpha_2\}}(x;\lambda) &=
 \Bigl(x+\frac{\lambda_1}3\Bigr)
 \Bigl(x+\frac{\lambda_1}6-\frac{(\alpha_1,\rho)}6\Bigr)
 \Bigl(x-1\Bigr)\\
&\qquad\cdot
 \Bigl(x-\frac{\lambda_1}6-\frac{(3\alpha_1+\alpha_2,\rho)}6\Bigr)
 \Bigl(x-\frac{\lambda_1}3-\frac{(4\alpha_1+2\alpha_2,\rho)}6\Bigr)\\
&=\Bigl(x+\frac{\lambda_1}3\Bigr)
 \Bigl(x+\frac{\lambda_1}6-\frac16\Bigr)
 \Bigl(x-1\Bigr)
 \Bigl(x-\frac{\lambda_1}6-1\Bigr)
 \Bigl(x-\frac{\lambda_1}3-\frac53\Bigr).
\end{align*}
Moreover, from Remark~\ref{rm:Agap}, we get
\begin{align*}
  r_{\alpha_1}(\lambda)
  &=\ang{\lambda_\Theta+\rho}{(-\varpi_1+\alpha_1)-(-\varpi_1+\alpha_1+\alpha_2)}\\
     &\qquad\cdot\ang{\lambda_\Theta+\rho}{(-\varpi_1+\alpha_1)-(-\varpi_1+3\alpha_1+2\alpha_2)}\\
  &=2\ang{\lambda_\Theta+\rho}{\alpha_2}
      \ang{\lambda_\Theta+\rho}{\alpha_1+\alpha_2}\\
  &=\frac16(\lambda_2+1)(3\lambda_2+4),
\\
  r_{\alpha_2}(\lambda)
  &= \Bigl( 
      \ang{\lambda_\Theta}{(-\varpi_1+\alpha_1)-(-\varpi_1)}
      -\ang{\alpha_2}{-\varpi_1+\alpha_1}  
     \Bigr)\\
     &\qquad\cdot
     \Bigl(
      \ang{\lambda_\Theta}{(-\varpi_1+\alpha_1+\alpha_2)-0}
      +\ang{-\varpi_1+\alpha_1+\alpha_2}{\alpha_1}     
     \Bigr)\\     
     &\qquad\cdot\ang{\lambda_\Theta+\rho}{(-\varpi_1+\alpha_1+\alpha_2)
       -(-\varpi_1+3\alpha_1+\alpha_2)}\\
     &\qquad\cdot\ang{\lambda_\Theta+\rho}{(-\varpi_1+\alpha_1+\alpha_2)
       -(-\varpi_1+4\alpha_1+2\alpha_2)}\\
  &=-\frac29\Bigl(\ang{\lambda_\Theta+\rho}{\alpha_1}\Bigr)
      \Bigl(\ang{\lambda_\Theta+\rho}{3\alpha_1+2\alpha_2}\Bigr)
      \Bigl(\ang{\lambda_\Theta+\rho}{3\alpha_1+\alpha_2}\Bigr)^2\\
  &=-\frac1{216}(\lambda_1+1)(\lambda_1+2)^2(\lambda_1+3).
\end{align*}
Here we have used the following relations:
\[
\left\{\begin{aligned}
&-\ang{\alpha_2}{-\varpi_1+\alpha_1}
=-\ang{\alpha_2}{\alpha_1}
=\frac{\ang{\rho}{3\alpha_1+2\alpha_2}}3,\\
&\ang{-\varpi_1+\alpha_1+\alpha_2}{\alpha_1}
=-\ang{\alpha_1}{\alpha_1}
=-\frac{\ang{\rho}{3\alpha_1+\alpha_2}}3.
\end{aligned}\right.\]
Note that $\alpha_1+\alpha_2, 3\alpha_1+2\alpha_2, 
3\alpha_1+\alpha_2\in\Sigma(\mathfrak g)$
and $r_{\alpha_i}(\lambda)\ne0$ if the condition ii) of Theorem~\ref{thm:gapexist} 
(we do not assume here that $\lambda_\Theta+\rho$ is dominant) is satisfied.

Let $S(\mathfrak a)^{(m)}$ denote the space of the elements of the symmetric algebra
over $\mathfrak a$ whose degree are at most $m$.
Note that
\begin{align*}
 (\trace F_\pi^{2m})_{\mathfrak a}&\equiv
   2(\varepsilon_1-\varepsilon_2)^{2m}+2(\varepsilon_2-\varepsilon_3)^{2m}
   +2(\varepsilon_1-\varepsilon_3)^{2m}\mod S(\mathfrak a)^{(2m-1)}\\
   &\equiv 2(\varepsilon_1-\varepsilon_2)^{2m}+2(\varepsilon_1+2\varepsilon_2)^{2m}
    +2(2\varepsilon_1+\varepsilon_2)^{2m}\\
   &\quad\quad\quad\quad\quad\quad\quad\quad\quad
   \mod S(\mathfrak a)(\varepsilon_1+\varepsilon_2+\varepsilon_3),\\
 (\trace F_\pi^2)_{\mathfrak a}&\equiv
   12(\varepsilon_1^2+\varepsilon_1\varepsilon_2+\varepsilon_2^2)
   \mod S(\mathfrak a)^{(1)}+S(\mathfrak a)
  (\varepsilon_1+\varepsilon_2+\varepsilon_3),\\
  (\trace F_\pi^4)_{\mathfrak a}&\equiv
  \frac14((\trace F_\pi^2)_{\mathfrak a})^2
  \mod S(\mathfrak a)^{(3)}+S(\mathfrak a)
  (\varepsilon_1+\varepsilon_2+\varepsilon_3).
\end{align*}
Moreover $(\trace F_\pi^6 )_{\mathfrak a}$ and 
$\bigl((\trace F_\pi^2)_{\mathfrak a}\bigr)^3$
are linearly independent in 
 \[S(\mathfrak a)/\Bigl(S(\mathfrak a)^{(5)}+S(\mathfrak a)
  (\varepsilon_1+\varepsilon_2+\varepsilon_3)\Bigr).\]
Thus we have 
\begin{equation}
Z(\mathfrak g)=\mathbb C[\trace F_\pi^2,\trace F_\pi^6].
\end{equation}
\end{exmp}

\begin{prop}\label{prop:every}
We denote by $\alpha_i$ the elements in $\Psi(\mathfrak g)$
which are specified by the Dynkin diagrams in the examples in this section.

For $\alpha\in\Psi(\mathfrak g)$
define $\Lambda_\alpha\in\mathfrak a^*$ by
\begin{equation}
 2\frac{\ang{\Lambda_\alpha}{\beta}}{\ang{\beta}{\beta}}
 =\begin{cases}
   1\quad&\text{if }\beta=\alpha,\\
   0&\text{if }\beta\in\Psi(\mathfrak g)\setminus\{\alpha\}.
  \end{cases}
\end{equation}
Let $\pi_\alpha^*$ be the irreducible representation of $\mathfrak g$
with the lowest weight $-\Lambda_\alpha$ and let $\Lambda^*_\alpha$ be the highest
weight of $\pi^*_\alpha$.

\noindent
{\rm i)}
Suppose $\mathfrak g=\mathfrak{gl}_n$, $\mathfrak{sl}_n$, $\mathfrak{sp}_n$ or 
$\mathfrak{o}_{2n+1}$ and $\pi$ is the natural representation of $\mathfrak g$.
Then \eqref{eq:gap} holds for any $\Theta$ if the infinitesimal character of the
Verma module $M(\lambda_\Theta)$ is regular, that is
\begin{equation}\label{eq:regular}
 \ang{\lambda_\Theta+\rho}{\alpha}
 \ne 0\quad(\forall \alpha\in\Sigma(\mathfrak g)).
\end{equation}
If $\lambda_\Theta+\rho$ is dominant,
then \eqref{eq:gap} is equivalent to \eqref{eq:gapann}.
Moreover in\/ {\rm Proposition~\ref{prop:geninv}} we may put
$A=\{i;\, d_i<\deg_x q_{\pi,\Theta} \}$.

\noindent
{\rm ii)}
Suppose $\mathfrak g=G_2$ and $\pi$ is the non-trivial minimal dimensional 
representation of $\mathfrak g$.  Then the same statement as above holds.

\noindent
{\rm iii)}
Suppose $\mathfrak g=\mathfrak o_{2n}$ with $n\ge 4$ and $\pi$ is the natural 
representation of $\mathfrak g$. 

Suppose $\Theta\supset\{\alpha_{n-1},\alpha_n\}$.
Then \eqref{eq:gap} holds if $\lambda_\Theta+\rho$ is regular
and \eqref{eq:gap} is equivalent to \eqref{eq:gapann}
if $\lambda_\Theta+\rho$ is dominant.

Suppose $\Theta\cap\{\alpha_{n-1},\alpha_n\}=\emptyset$
and $\ang{\lambda_\Theta}{\alpha_n-\alpha_{n-1}}=0$.
In this case we may replace $q_{\pi,\Theta}(x;\lambda)$ in the definition of 
$I_{\pi,\Theta}$
by $q'_{\pi,\Theta}(x;\lambda)$ given in\/{\rm~Example~\ref{ex:Dn}}.  
Then the same statement as the previous case holds.
Note that $\deg_x q'_{\pi,\Theta}= \deg_x q_{\pi,\Theta}-1$.

In other general cases, \eqref{eq:gap} holds if the infinitesimal
character of $M(\lambda_\Theta)$ is strongly regular, that is, 
$\lambda_\Theta+\rho$ is not fixed by any non-trivial element of the Weyl group 
of the non-connected Lie group $O(2n,\mathbb C)$.
In particular, if $\Theta\cap\{\alpha_{n-1},\alpha_n\}=\emptyset$,
then \eqref{eq:gap} holds under the conditions
\eqref{eq:regular} and 
\begin{multline}\label{eq:O2n}
  \ang{\lambda_\Theta%
  +\rho}
 {2\alpha_i+\cdots+2\alpha_{n-2}+\alpha_{n-1}+\alpha_n}
 \ne0\\
 \text{for }i=2,\ldots,n-1\text{ satisfying }\alpha_{i-1}\in\Theta
 \text{ and }\alpha_i\notin\Theta.
\end{multline}

Suppose $\Theta\cap\{\alpha_{n-1},\alpha_n\}=\{\alpha_{n-1}\}$.
Then
\begin{equation}
 J_\Theta(\lambda) = I_{\pi,\Theta}(\lambda) + I_{\pi^*_{\alpha_{n-1}},\Theta}(\lambda)
   + J(\lambda_\Theta) 
\end{equation}
if \eqref{eq:regular}, \eqref{eq:O2n} and
\begin{multline}
 \ang{\lambda_\Theta+\rho}{\varpi+\Lambda_{\alpha_{n-1}}-\alpha_{n-1}}\ne 0\\
 \text{for any \ }\varpi\in\mathcal W_\Theta(\pi^*_{\alpha_{n-1}})
 \text{ \ satisfying \ }\varpi>\alpha_{n-1}-\Lambda_{\alpha_{n-1}}
\end{multline}
hold.

In\/ {\rm Proposition~\ref{prop:geninv}} we may put $r=n$ and 
$\Delta_1,\ldots,\Delta_{n-1}$ are invariant
under the outer automorphism of $\mathfrak g$ corresponding to 
$\varepsilon_n\mapsto-\varepsilon_n$ and
$A=\{i; d_i<\deg_x q_{\pi,\Theta}\}\cup\{n\}$.

\noindent{\rm iv)}
Suppose $\mathfrak g=E_n$ with $n=6$, $7$ or $8$ {\rm(}cf.~{\rm Example~\ref{ex:E6},
\ref{ex:E7}, \ref{ex:E8})}.
For $\alpha_i\in\Psi(\mathfrak g)$ put
\begin{equation}
 \iota(\alpha_i)=
  \begin{cases}
    \alpha_1&\text{if }i=1\text{ or }3,\\
    \alpha_2&\text{if }i=2,\\
    \alpha_n&\text{if }i\ge 4,
  \end{cases}
\quad\quad\quad
 {\hat\alpha}_i =
   \begin{cases}
    \alpha_i&\text{if }i=1\text{ or }2,\\
    \alpha_1+\alpha_3&\text{if }i=3,\\
    \alpha_i+\cdots+\alpha_n&\text{if }i\ge 4.
   \end{cases}
\end{equation}
Here $\iota(\alpha_i)$ satisfies 
$\#\{\beta\in\Psi(\mathfrak g);\,\ang{\iota(\alpha_i)}{\beta}<0\}\le 1$
and $\hat\alpha$ is the smallest root with $\hat\alpha\ge\alpha$ and 
$\hat\alpha\ge\iota(\alpha)$.
Let $\lambda\in\mathfrak a_\Theta^*$.
If \eqref{eq:regular} holds and moreover 
$\lambda$ satisfies \begin{multline}\label{eq:fundrep}
 2\ang{\lambda_\Theta+\rho}{\varpi + \Lambda_{\iota(\alpha)} - \hat\alpha}
 \ne\ang{\varpi}{\varpi}-\ang{\Lambda_{\iota(\alpha)}}{\Lambda_{\iota(\alpha)}}\\
 \quad\text{for }\alpha\in\Theta\text{ and }
 \varpi\in
 \overline{\mathcal W}_\Theta(\pi_{\iota(\alpha)}^*)
 \text{ satisfying }\varpi>\hat\alpha - \Lambda_{\iota(\alpha)},
\end{multline}
then
\begin{equation}\label{eq:gapexc}
  J_\Theta(\lambda) = 
  \sum_{\alpha\in\iota(\Theta)}I_{\pi_\alpha^*, \Theta}(\lambda)
  + J(\lambda_\Theta).
\end{equation}
In particular, under the notation in\/ {\rm Definition~\ref{def:shift}} the condition
\begin{multline}\label{eq:gapexc2}
 2\frac{\ang{\lambda_\Theta+\rho}{\mu}}
 {\ang{\Lambda_{\iota(\alpha)}}{\Lambda_{\iota(\alpha)}}}
 \notin[-1, 0]\\
 \quad\text{for }\alpha\in\Theta
 \text{ and }
  \mu\in R_+ 
 \text{ with }0<\mu\le
  \Lambda_{\iota(\alpha)}+\Lambda_{\iota(\alpha)}^*-\hat\alpha
\end{multline}
assures \eqref{eq:fundrep}.
Moreover,
if $\pi=\pi^*_{\alpha_1}$ or $\pi^*_{\alpha_n}$,
we may put 
$A=\{i;\, d_i<\deg_x q_{\pi,\Theta} \}$
 in\/ {\rm Proposition~\ref{prop:geninv}}.
 
\noindent{\rm v)}
Suppose $\mathfrak g=F_4$. 
For $\alpha_i\in\Psi(\mathfrak g)$ put
\begin{equation}
 \iota(\alpha_i) =
  \begin{cases}
  \alpha_1&\text{if }i\le 2,\\
  \alpha_4&\text{if }i\ge 3,\
  \end{cases}
 \quad\quad
 {\hat\alpha}_i =
  \begin{cases}
   \alpha_i&\text{if }i = 1\text{ or }4,\\
   \alpha_1+\alpha_2&\text{if }i=2,\\
   \alpha_3+\alpha_4&\text{if }i=3.
  \end{cases}
\end{equation}
Then the same statement as\/ {\rm iv)} holds
for $\pi=\pi^*_{\alpha_4}$ {\rm(}cf.\/{\rm ~Example~\ref{ex:F4})}.
\end{prop}
\begin{proof}
The statements i) and iii) are direct consequences of
\cite[Theorem~4.4]{O-Cl} (or Theorem~\ref{thm:gap})
and Theorem~\ref{thm:gapexist}.
The statement ii) is a consequence of Example~\ref{ex:G2}.

Suppose $\mathfrak g$ is $E_6$, $E_7$, $E_8$, $F_4$ or $G_2$ and $\pi$ is a 
minimal dimensional non-trivial irreducible representation of $\mathfrak g$.
Then in Proposition~\ref{prop:geninv} it follows from \cite{Me} that the elements
$\sum_{\varpi\in\mathcal W(\pi)}m_\pi(\varpi)\varpi^{d_i}$ ($i=1,\ldots,n$)
generate the algebra of the $W$-invariants of $U(\mathfrak a)$ 
(For $G_2$ we confirm it in Example~\ref{ex:G2}) and 
hence we may put $A=\{i;\, d_i<\deg_x q_{\pi,\Theta} \}$.

Suppose $\mathfrak g$ is $E_6$, $E_7$, $E_8$ or $F_4$.
Fix $\alpha\in\Theta$. 
Then Theorem~\ref{thm:gap} assures 
$X_{-\alpha}\in I_{\pi_{\iota(\alpha)}^*,\Theta}(\lambda) + J(\lambda_\Theta)$
if $r_{\alpha,\varpi_\alpha}(\lambda)\ne 0$.
Here $r_{\alpha,\varpi_\alpha}(\lambda)$ is defined by \eqref{eq:gapcond}
with $\pi=\pi_{\iota(\alpha)}^*$ and 
$\varpi_\alpha=-\Lambda_{\iota(\alpha)}+(\hat\alpha-\alpha)$.
Then %
the assumption of Remark~\ref{rm:Agap}~v) holds
and therefore the second factor $\prod_{i=1}^{L}\bigl(\cdots\bigr)$ of 
$r_{\alpha,\varpi_\alpha}(\lambda)$ in \eqref{eq:gapcond} does not vanish under the condition 
\eqref{eq:regular}.
On the other hand, 
$\varpi\in\mathcal W(\pi_{\iota(\alpha)}^*)$ which does
not satisfy $\varpi\le-\Lambda_{\iota(\alpha)}+\hat\alpha$ always satisfies 
$\varpi>-\Lambda_{\iota(\alpha)}+\hat\alpha$ because $\{\gamma_1,\ldots,\gamma_K\}$ 
in Remark~\ref{rm:Agap} is of type $A_K$ and 
$\ang{\Lambda_{\iota(\alpha)}}{\beta}=\ang{\gamma_i}{\beta}=0$
for $i=1,\ldots,K-1$ and $\beta\in\Psi(\mathfrak g)
\setminus\{\gamma_1,\ldots,\gamma_K\}$.
Hence \eqref{eq:fundrep} assures that
the first factor of $r_{\alpha,\varpi_\alpha}(\lambda)$ does not vanish.
Thus we have 
$X_{-\alpha}\in I_{\pi_{\iota(\alpha)}^*,\Theta}(\lambda) + J(\lambda_\Theta)$.
It implies \eqref{eq:gapexc}.
It is clear that \eqref{eq:fundrep} follows from \eqref{eq:gapexc2} since 
$\ang{\Lambda_{\iota(\alpha)}}{\Lambda_{\iota(\alpha)}}\ge\ang{\varpi}{\varpi}$
for $\varpi\in\mathcal W(\pi^*_{\iota(\alpha)})$.
\end{proof}

\begin{rem}\label{rm:posW}
Suppose $\mathfrak g=\mathfrak{gl}_n$ or $\mathfrak g$ is simple.
In the preceding proposition we explicitly give a two sided ideal 
$I_\Theta(\lambda)$ of $U(\mathfrak g)$ which satisfies
$J_\Theta(\lambda) = I_\Theta(\lambda) + J(\lambda_\Theta)$ if at least
\begin{equation}
 \operatorname{Re}\ang{\lambda_\Theta+\rho}{\alpha}>0
 \quad\text{for }\alpha\in\Psi(\mathfrak g).
\end{equation}
In particular, this condition is valid when $\lambda=0$.
\end{rem}

\begin{rem}\label{rm:existence}
Suppose $\mathfrak g=\mathfrak{gl}_n$.  Then in \cite{O-Cap} the generator system
of $\Ann\bigl(M_\Theta(\lambda)\bigr)$ is constructed for any $\Theta$ and 
$\lambda$ through quantizations of elementary divisors.
It shows that the zeros of the image of the Harish-Chandra
homomorphism of $\Ann\bigl(M_\Theta(\lambda)\bigr)$ equals 
$\{w.\lambda_\Theta;\,w\in W(\Theta)\}$ and proves that \eqref{eq:gapann}
holds if and only if \eqref{eq:gapgln} is not valid for any positive numbers 
$j$ and $k$ which
are smaller or equal to $L$.
Here we note that this condition for \eqref{eq:gapann} follows from this description 
of the zeros and Lemma~\ref{lem:gap} and the following Lemma with the notation in 
Example~\ref{ex:A}.
\end{rem}

\begin{lem}
Let $n_0=0<n_1<n_2<\cdots<n_L=n$ be a strictly increasing sequence
of non-negative integers.
Let $\lambda=(\lambda_1,\ldots,\lambda_L)\in\mathbb C^L$.
Define $\bar\lambda=(\bar\lambda_1,\ldots,\bar\lambda_n)\in\mathbb C^n$
by
\[
  \bar\lambda_\nu = \lambda_k + (\nu-1) -\frac{n-1}2
  \text{ if }n_{k -1}< \nu\le n_k
\]
and put
\[
 \Lambda_k = \{\bar\lambda_{n_{k-1}+1}, \bar\lambda_{n_{k-1}+2},\ldots,
  \bar\lambda_{n_k}\}.
\]
Then there exists $\nu$ with $n_{j-1}<\nu< n_j$ 
satisfying
$(\nu,\nu+1)\bar\lambda\in W(\Theta)\bar\lambda$
if and only if there exists $k\in\{1,\ldots,L\}$ such that 
\begin{equation}\label{eq:gapgln}
 \Lambda_k\cap\Lambda_j\ne\emptyset,\ \Lambda_j\not\subset\Lambda_k
 \text{ and }\Bigl(
 \mu\in\Lambda_j\setminus\Lambda_k,\ \mu'\in\Lambda_k
 \Rightarrow (\mu'-\mu)(k-j)>0
\Bigr).
\end{equation}
Here %
$(i,j)\in\mathfrak S_n$ is the transposition of $i$ and $j$ and
\begin{align*}
 W(\Theta) &= \{\sigma\in\mathfrak S_n;\, \sigma(i) < \sigma(j)\text{ if there
 exists $k$ with }n_{k-1}<i<j\le n_k\},\\
 \sigma\mu&=(\mu_{\sigma^{-1}(1)},\ldots,\mu_{\sigma^{-1}(n)})
  \text{ for }\mu=(\mu_1,\ldots,\mu_n)\in\mathbb C^n.
\end{align*}
\end{lem}
\begin{proof}
Suppose \eqref{eq:gapgln}.
Then there exists %
$m$ such that
\begin{align*}
&\begin{cases}
  j<k,\ 1\le m< n_j-n_{j-1}\text{ and }n_{k-1}+n_j-n_{j-1}-m\le n_k,\\
  \bar\lambda_{n_{j-1}+\nu}=\bar\lambda_{n_{k-1}+\nu-m}
  \text{ for } m<\nu\le n_j-n_{j-1}
\end{cases}
\intertext{or}
&\begin{cases}
  j>k,\ 1\le m< n_j-n_{j-1}\text{ and }n_k-m+1 > n_{k-1},\\
  \bar\lambda_{n_{j-1}+\nu}=\bar\lambda_{n_k+\nu-m}
  \text{ for } 1\le\nu\le m.
\end{cases}
\end{align*}
Defining $\sigma\in W(\Theta)$ by
\begin{align*}
\sigma &= (n_{j-1}+m,n_{j-1}+m+1)
   \prod_{m<\nu\le n_j-n_{j-1}}(n_{j-1}+\nu,n_{k-1}+ \nu-m),
\intertext{ or }
\sigma &= (n_{j-1}+m,n_{j-1}+m+1)
   \prod_{1\le\nu\le m}(n_{j-1}+\nu,n_k+\nu-m),
\end{align*}
respectively, we have $(\nu,\nu+1)\bar\lambda =
\sigma\bar\lambda \in W(\Theta)\bar\lambda$
with $\nu = n_{j-1}+m$.

Conversely suppose $(\nu,\nu+1)\bar\lambda=\sigma\bar\lambda$
for suitable $\nu\in\{n_{j-1}+1,\ldots,n_j-1\}$ and $\sigma\in W(\Theta)$.
Put
\begin{align*}
 \{\ell_1,\ldots,\ell_m\} 
  &= \{\ell;\, \ell \le n_{j-1}\text{ and }\bar\lambda_\ell = 
  \bar\lambda_{n_{j-1}+1}\},\\
 \{\ell'_{m+2},\ldots,\ell'_{m+m'+1}\} &=\{\ell';\, \ell' > n_j\text{ and }
  \bar\lambda_{\ell'} =  \bar\lambda_{n_j}\}
\end{align*}
and define
\[
 \begin{cases}
  \ell'_i = \ell_i + (n_j - n_{j-1}-1)\quad&\text{if }i\le m,\\
  \ell_i = \ell_i' - (n_j-n_{j-1}-1)&\text{if } i \ge m+2,\\
  \ell_{m+1} = n_{j-1}+1,\ \ell'_{m+1} = n_j.
 \end{cases}
\]
Assume that \eqref{eq:gapgln} is not valid for any $k$.
Then for $i\in I:=\{1,\ldots,m+m'+1\}$, there exist integers $N_i$ with
$n_{N_{i-1}}<\ell_i<\ell'_i\le n_{N_i}$ and therefore
$\bar\lambda_{\ell_i} = \bar\lambda_{n_{j-1}+1}$ and
$\bar\lambda_{\ell'_i} = \bar\lambda_{n_j}$.

Note that $\#I_1\le m+1$ and $\#I_2\le m'$ by denoting
\[
  I_1 = \{i\in I;\, \sigma(\ell_i)\le n_j\}\text{ and }
  I_2 = \{i\in I;\, \sigma(\ell'_i)>n_j\}.
\]
Since $\sigma(\ell_i) < \sigma(\ell_i')$, we have
$I_1\cup I_2=I$ and therefore
$\#I_1= m+1$ and $\#I_2= m'$.
Then there exists $i_0$ with $n_{j-1}<\sigma(\ell_{i_0})\le n_j$.
Since $I_1\cap I_2=\emptyset$,
we have $\sigma(\ell_{i_0}')\le n_j$,
which implies $\sigma^{-1}(\nu')=\ell_{i_0}+\nu'-n_{j-1}-1$
for $n_{j-1}<\nu'\le n_j$.
It contradicts to the assumption $(\nu,\nu+1)\bar\lambda=\sigma\bar\lambda$.
\end{proof}

\begin{rem}
Suppose $\mathfrak g=\mathfrak{gl}_n$ and $\pi$ is its natural representation.
Then the condition $r_\alpha(\lambda)\ne0$ for any $\alpha\in\Theta$ is necessary
and sufficient for \eqref{eq:gap} (cf.~\cite[Remark~4.5]{O-Cl}).  Under the notation in the preceding lemma, 
it is easy to see that the condition is equivalent to the fact that
\[
 \Lambda_k\cap\Lambda_j\ne\emptyset,\ \Lambda_j\not\subset\Lambda_k
 \text{ and }\Bigl(\exists
 \mu\in\Lambda_j\setminus\Lambda_k,\ \exists\mu'\in\Lambda_k
 \text{ such that } (\mu'-\mu)(k-j)>0
\Bigr)
\]
does not hold for any positive numbers $k$ and $j$ smaller or equal to $L$.
\end{rem}
\appendix 

\section{Infinitesimal Mackey's Tensor Product Theorem}\label{app:Mackey}
In this appendix
we explain {\it infinitesimal Mackey's tensor product theorem\/}
following the method given in ~\cite{Ma}.

Let $\mathfrak g$ be a finite dimensional Lie algebra over $\mathbb C$
and $\mathfrak p$ a subalgebra of $\mathfrak g$.
Let $V$ and $U$ be a $U(\mathfrak g)$-module 
and a $U(\mathfrak p)$-module, respectively.
We denote by $V|_{\mathfrak p}$ and
$\operatorname{Ind}_{\mathfrak p}^{\mathfrak g} U$
the restriction of the coefficient ring $U(\mathfrak g)$
to $U(\mathfrak p)$
and 
the induced representation $U(\mathfrak g) \otimes_{U(\mathfrak p)} U$
in the usual way.

\begin{thm}[infinitesimal Mackey's tensor product theorem]\label{thm:Mackey}
The map defined by
\begin{equation}\label{eq:mac1}
\begin{aligned}
U(\mathfrak g) \otimes_{U(\mathfrak p)}
(U \otimes_{\mathbb C} V|_{\mathfrak p}) &\rightarrow 
( U(\mathfrak g) \otimes_{U(\mathfrak p)} U )
\otimes_{\mathbb C} V,\\
D\otimes_{U(\mathfrak p)}(u\otimes_{\mathbb C} v) &\mapsto D\cdot [(1\otimes_{U(\mathfrak p)} u)\otimes_{\mathbb C} v]
\end{aligned}
\end{equation}
gives a canonical $U(\mathfrak g)$-module isomorphism
\begin{equation}\label{eq:mac2}
\operatorname{Ind}_{\mathfrak p}^{\mathfrak g}(U \otimes_{\mathbb C} V|_{\mathfrak p})
\simeq 
\left( \operatorname{Ind}_{\mathfrak p}^{\mathfrak g} U \right)
\otimes_{\mathbb C} V.
\end{equation}
\end{thm}

To prove this
we need two lemmas.

\begin{lem}\label{lem:ring_equiv}
Let $R$ be a ring and $R\text{-\sf Mod}$ the category of left $R$-modules.
For $M, N \in R\text{-\sf Mod}$
consider 
$F_M:\cdot \mapsto \Hom_R(M, \cdot)$ and
$F_N:\cdot \mapsto \Hom_R(N, \cdot)$,
which are functors from $R\text{-\sf Mod}$
to the category of abelian groups.
Suppose that $F_M$ and $F_N$ are naturally equivalent,
namely, there exists an assignment $A \mapsto \tau_A$
for each object $A \in R\text{-\sf Mod}$ of an isomorphism 
$\tau_A:\Hom_R(M, A)\rightarrow \Hom_R(N, A)$
such that 
$F_N(f) \circ \tau_A = \tau_B \circ F_M(f)$
for each $f \in \Hom_R(A, B)$.
Then $M \simeq N$ as $R$-modules.
\end{lem}
\begin{proof}
Put
$\varphi = \tau_N^{-1}(\operatorname{id}_N) \in \Hom_R(M, N)$
and
$\psi = \tau_M(\operatorname{id}_M) \in \Hom_R(N, M)$.
Then $\varphi\circ\psi = F_N(\varphi)(\psi)
=F_N(\varphi)\circ\tau_M(\operatorname{id}_M)
=\tau_N\circ F_M(\varphi)(\operatorname{id}_M)
=\tau_N(\varphi)=\operatorname{id}_N.$
Similarly $\psi\circ\varphi=\operatorname{id}_M.$
Hence $M \simeq N$.
\end{proof}

\begin{lem}\label{lem:3_modules}
Let $(\pi_i, V_i)\ (i=1, 2, 3)$ be $U(\mathfrak g)$-modules.
Consider $\Hom_{\mathbb C}(V_2, V_3)$
as a $U(\mathfrak g)$-module
by $X \Phi = \pi_3(X)\circ\Phi - \Phi\circ\pi_2(X)$ 
for $\Phi\in\Hom_{\mathbb C}(V_2, V_3)$
and $X \in \mathfrak g$.
Then naturally
\[ \Hom_{U(\mathfrak g)}(V_1 \otimes_{\mathbb C} V_2, V_3)
\simeq \Hom_{U(\mathfrak g)}(V_1, \Hom_{\mathbb C}(V_2, V_3)).
\]
\end{lem}
\begin{proof}
We have only to define the mapping $\varphi \mapsto \Phi$ 
from the left-hand side to the right-hand side
by $\left(\Phi(v_1)\right)(v_2)
=\varphi(v_1 \otimes v_2)$
for $v_1 \in V_1$ and $v_2 \in V_2$.
\end{proof}

\begin{proof}[Proof of\/ {\rm Theorem \ref{thm:Mackey}}]
Lemma \ref{lem:3_modules}
implies the following isomorphism
for a given $U(\mathfrak g)$-module $A$:
\[\begin{split}
\Hom_{U(\mathfrak g)}\left(
(U(\mathfrak g) \otimes_{U(\mathfrak p)} U )
\otimes_{\mathbb C} V, A
\right)
&\simeq
\Hom_{U(\mathfrak g)}\left(
U(\mathfrak g) \otimes_{U(\mathfrak p)} U,
\Hom_{\mathbb C}(V, A)
\right) \\
&\simeq
\Hom_{U(\mathfrak p)}\left(
U,
\Hom_{\mathbb C}(V|_{\mathfrak p}, A|_{\mathfrak p})
\right) \\
&\simeq
\Hom_{U(\mathfrak p)}\left(
U \otimes_{\mathbb C} V|_{\mathfrak p},
A|_{\mathfrak p}
\right) \\
&\simeq
\Hom_{U(\mathfrak g)}\left(
U(\mathfrak g) \otimes_{U(\mathfrak p)} 
(U \otimes_{\mathbb C} V|_{\mathfrak p}),
A
\right).
\end{split}\]
It gives a natural equivalence between
$F_{(U(\mathfrak g) \otimes_{U(\mathfrak p)} U )
\otimes_{\mathbb C} V}$ and 
$F_{U(\mathfrak g) \otimes_{U(\mathfrak p)} 
(U \otimes_{\mathbb C} V|_{\mathfrak p})}$
under the notation of Lemma \ref{lem:ring_equiv} with $R=U(\mathfrak g)$.
Hence by Lemma \ref{lem:ring_equiv},
we have \eqref{eq:mac2}.
It is easy to see the isomorphism is explicitly given by \eqref{eq:mac1}.
\end{proof}

\section{Undesirable Cases}\label{app:badcase}
In this appendix we give counter examples stated in Remark~\ref{rem:badcase}.
Let $\mathfrak g=\mathfrak{sl}_n$ and use the notation in \S\ref{sec:min} and \S\ref{sec:ideal}.
Suppose the Dynkin diagram of the fundamental system $\Psi=\{\alpha_1,\ldots,\alpha_{n-1}\}$ 
is the same as in Example~\ref{ex:A}. 
Let $\{\Lambda_1,\ldots,\Lambda_{n-1}\}$ be
the system of fundamental weights corresponding to $\Psi$.
Let $\pi$ be the irreducible representation of $\mathfrak g$
with lowest weight $\bar\pi=-m_1\Lambda_1-m_2\Lambda_2$.
Here $m_1$ and $m_2$ are %
positive integers. 
Then the multiplicity of the weight $\varpi':=\bar\pi+\alpha_1+\alpha_2\in\mathcal{W}(\pi)$
equals $2$.

Now take $\Theta=\Psi\setminus\{\alpha_2\}=\{\alpha_1,\alpha_3,\alpha_4,\ldots,\alpha_{n-1}\}$.
Since the multiplicity of the weight $\bar\pi+\alpha_2$
is $1$, both $\varpi'$ and $\bar\pi+\alpha_2$ belong to $\overline{\mathcal{W}}_\Theta(\pi)$.
On the other hand, by Remark~\ref{rem:lowweight}, the weight
$\varpi_{\alpha_{n-1}}:=\bar\pi+\alpha_2+\alpha_3+\cdots+\alpha_{n-2}$
is a unique extremal low weight of $\pi$
with respect to $\alpha_{n-1}$.
Note that 
$\{\varpi\in\overline{\mathcal{W}}_\Theta(\pi);\,
 \varpi\leq\varpi_{\alpha_{n-1}}\}
 =\{\bar\pi, \bar\pi+\alpha_2\}$
and the weight $\varpi'_{\alpha_{n-1}}:=\bar\pi+\alpha_2+\alpha_3+\cdots+\alpha_{n-1}$
satisfies 
$\varpi'_{\alpha_{n-1}}|_{\mathfrak a_\Theta}=\varpi'|_{\mathfrak a_\Theta}
 =(\bar\pi+\alpha_2)|_{\mathfrak a_\Theta} \ne \bar\pi|_{\mathfrak a_\Theta}$.
Moreover, it follows from Lemma~\ref{lem:order}
\begin{align*}
D_\pi(\varpi')-D_\pi(\bar\pi+\alpha_2)&=-\ang{\bar\pi+\alpha_2}{\alpha_1}
= \frac{m_1+1}2\ang{\alpha_1}{\alpha_1},\\
D_\pi(\varpi'_{\alpha_{n-1}})-D_\pi(\bar\pi+\alpha_2)&=
-\ang{\alpha_2}{\alpha_3}-\cdots-\ang{\alpha_{n-2}}{\alpha_{n-1}}
=\frac{n-3}2\ang{\alpha_1}{\alpha_1}.
\end{align*}
It shows the first factor of the function \eqref{eq:gapcond}
with ($\alpha, \varpi_\alpha)=(\alpha_{n-1}, \varpi_{\alpha_{n-1}})$
is identically zero
if $n=m_1+4$.

\end{document}